\RequirePackage{fix-cm}
\documentclass{svjour3}                     
\smartqed  
\usepackage{graphicx}
\usepackage{verbatim,xcolor,soul,relsize}
\usepackage{amsmath,lineno,amssymb,mathtools,url}

\newtheorem{dfn}{Definition}
\newtheorem{pro}[dfn]{Problem}
\newtheorem{thm}[dfn]{Theorem}
\newtheorem{exa}[dfn]{Example}
\newtheorem{rem}[dfn]{Remark}

\newtheorem{cor}[dfn]{Corollary}
\newtheorem{lem}[dfn]{Lemma}
\newtheorem{prop}[dfn]{Proposition}
\counterwithin{dfn}{section}

\newcommand{\C}{\mathbb C}
\newcommand{\R}{\mathbb R}
\newcommand{\Z}{\mathbb Z}

\newcommand{\al}{\alpha}
\newcommand{\be}{\beta}

\newcommand{\de}{\delta}

\newcommand{\La}{\Lambda}
\newcommand{\si}{\sigma}
\newcommand{\ep}{\varepsilon}

\newcommand{\inv}{\mathrm{Inv}}

\newcommand{\GL}{\mathrm{GL}}
\newcommand{\SO}{\mathrm{SO}}
\newcommand{\Or}{\mathrm{O}}

\newcommand{\PDD}{\mathrm{PDD}}
\newcommand{\RSD}{\mathrm{RSD}}
\newcommand{\AMD}{\mathrm{AMD}}

\newcommand{\sign}{\mathrm{sign}}

\newcommand{\DC}{\mathrm{DC}}
\newcommand{\DT}{\mathrm{DT}}
\newcommand{\QS}{\mathrm{QS}}
\newcommand{\VF}{\mathrm{VF}}
\newcommand{\CF}{\mathrm{CF}}
\newcommand{\RC}{\mathrm{RC}}
\newcommand{\PC}{\mathrm{PC}}
\newcommand{\RI}{\mathrm{RI}}
\newcommand{\PRF}{\mathrm{\overline{\RI}}}
\newcommand{\PI}{\mathrm{PI}}

\newcommand{\Oct}{\mathrm{Oct}}
\newcommand{\TC}{\mathrm{TC}}
\newcommand{\TP}{\mathrm{TP}}
\newcommand{\QT}{\mathrm{QT}}
\newcommand{\MS}{\mathrm{MS}}
\newcommand{\CM}{\mathrm{CM}}
\newcommand{\RM}{\mathrm{RM}}
\newcommand{\PM}{\mathrm{PM}}
\newcommand{\SIM}{\mathrm{SIM}}
\newcommand{\SSM}{\mathrm{SSM}}
\newcommand{\RIS}{\mathrm{RIS}}
\newcommand{\PIN}{\mathrm{PIN}}
\newcommand{\OSI}{\mathrm{OSI}}
\newcommand{\OSS}{\mathrm{OSS}}
\newcommand{\LIS}{\mathrm{LIS}}
\newcommand{\LSS}{\mathrm{LSS}}

\newcommand{\Obt}{\mathrm{Obt}}
\newcommand{\RB}{\mathrm{Red}}
\newcommand{\ws}{\hfill $\square$}
\newcommand{\bt}{\hfill $\blacktriangle$}
\newcommand{\bs}{\hfill $\blacksquare$}
\newcommand{\lra}{\leftrightarrow}
\newcommand{\Lra}{\Leftrightarrow}
\newcommand{\bd}{\partial}
\newcommand{\vl}{\,:\,}
\newcommand{\matfour}[4]{\left(\begin{array}{ccc}
#1 & #2 \\ #3 & #4 \end{array}\right)}

\begin{document}

\title{Mathematics of 2-dimensional lattices} 
\subtitle{Continuously parameterised spaces of all 2-dimensional lattices classified up to similarity, isometry, or rigid motion}


\author{Vitaliy Kurlin 
}


\institute{V.Kurlin \at 
Computer Science, University of Liverpool, UK \email{vitaliy.kurlin@liverpool.ac.uk}
}

\date{Received: date / Accepted: date}
\maketitle

\begin{abstract}
A periodic lattice in Euclidean space is the infinite set of all integer linear combinations of basis vectors.
Any lattice can be generated by infinitely many different bases.
This ambiguity was only partially resolved, but
standard reductions remained discontinuous under perturbations modelling crystal vibrations.

This paper completes a continuous classification of 2-dimensional lattices up to Euclidean isometry (or congruence), rigid motion (without reflections), and similarity (with uniform scaling).
The new homogeneous invariants allow easily computable metrics on lattices considered up to the equivalences above.
The metrics up to rigid motion are especially non-trivial and settle all remaining questions on (dis)continuity of lattice bases.
These metrics lead to real-valued chiral distances that continuously measure a lattice deviation from a higher-symmetry neighbour.

\keywords{Lattice \and rigid motion \and isometry \and invariant \and metric \and continuity}

\end{abstract}

\section{Motivations for a new continuous classification of lattices}
\label{sec:intro}

\begin{figure}[h]
\includegraphics[width=1.0\textwidth]{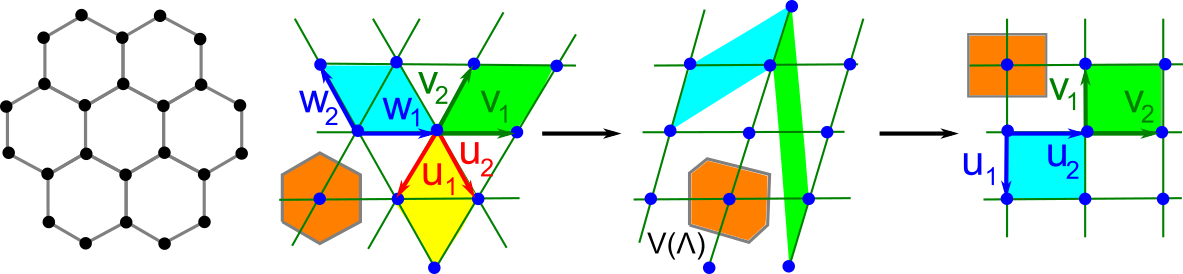}
\label{fig:hexagonal_graphene}
\caption{
\textbf{Left}: a 2-dimensional layer of graphene is formed by carbon atoms.
\textbf{Right}: 
one can generate a hexagonal lattice (as any other) by infintiely many bases and continuously deform into a rectangular lattice (or any other) whose bases $\{v_1,v_2\}$ and $\{u_1,u_2\}$ are related by an orientation-reversing map.
The yellow Voronoi domain $V(\La)$ of any point $p$ in a lattice $\La$ consists of all points $q\in\R^2$ that are non-strictly closer to $p$ than to other points of $\La-p$.
}
\end{figure}

A \emph{lattice} $\La\subset\R^n$ consists of all integer linear combinations of basis vectors $v_1,\dots,v_n$.
This basis spans a parallelepiped called a \emph{unit cell} $U\subset\R^n$.
A periodic point set is obtained as a union of translated copes $\La+p_i$ for finitely many $p_1,\dots,p_m\in U$.
Any crystal structure can be modelled as a periodic set whose points represent atomic centers.
For example, graphene is a 2-dimensional periodic set of carbon atoms based on a hexagonal lattice, see Fig.~\ref{fig:hexagonal_graphene}.
\medskip

Since crystal structures are determined in a rigid form, the most fundamental equivalence of their underlying lattices is a rigid motion.
Any \emph{rigid motion} in $\R^2$ is a composition of translations and rotations.
A more general \emph{isometry} includes mirror reflections and is sometimes called a \emph{congruence} in Euclidean geometry.
\medskip

In the language of Computer Science, the classification of lattices up to isometry is a binary classification problem deciding if lattices $\La,\La'$ are isometric, which can be denoted as $\La\cong\La'$.
If a certain approach wrongly indicates that isometric lattices $\La\cong\La'$ are in different isometry classes, this pair $\La,\La'$ is called a \emph{false negative}.
Many descriptors of crystals and their lattices allow false negatives by a simple comparison of lattice bases.
Any lattice can be represented by a reduced cell \cite{gruber1989reduced}, see Definition~\ref{dfn:reduced_cell} in section~\ref{sec:review}, which is unique up to isometry but this cell still has different bases as in Fig.~\ref{fig:hexagonal_graphene}. 
A descriptor without false negatives takes the same value on all isometric lattices and can be called an \emph{isometry invariant}.
\medskip

For example, the area of the unit cell $U$ spanned by any basis of a lattice $\La$ is an isometry invariant because a change of basis is realised by a $2\times 2$ matrix with determinant $\pm 1$, which preserves the absolute value of the area.
Such an invariant $I$ may allow \emph{false positives} $\La\not\cong\La'$ with $I(\La)=I(\La')$.
All lattices in Fig.~\ref{fig:hexagonal_graphene} have unit cells of the same area.
The area and many other invariants allow infinitely many false positives.
An invariant $I$ without false positives is called \emph{complete} and distinguishes all non-isometric lattices so that if $I(\La)=I(\La')$ then $\La\cong\La'$. 
\medskip

The traditional approach to deciding if lattices are isometric is to compare their conventional or reduced cells up to isometry.
Though this comparison theoretically gives a complete invariant, in practice all real crystal lattices are non-isometric because of inevitable noise in measurements. 
All atoms vibrate above the absolute zero temperature, hence any real lattice basis is always slightly perturbed.
The discontinuity of reduced bases under perturbations was experimentally known since 1980 \cite{andrews1980perturbation}, highlighted in \cite[section~1]{edels2021} and formally proved in \cite[Theorem~15]{widdowson2022average}.
\medskip

A more practically important goal is to find a complete invariant that is continuous under any perturbations of lattices.
Such a continuous and complete invariant will unambiguously parameterise the \emph{Lattice Isometry Space} ($\LIS$) consisting of infinitely many isometry classes of lattices.
For example, the latitude and longitude also continuously parameterise the surface of Earth. 
\medskip

The space $\LIS$ of isometry classes is continuous and connected because any two lattices can be joined by a continuous deformation of their bases as in Fig.~\ref{fig:hexagonal_graphene}.
Such deformation can be always visualised as a continuous path in the space $\LIS$, whose full geometry remained unknown even for 2-dimensional lattices.
\medskip

The main contribution is a full solution to the mapping problem below.

\begin{pro}[lattice mapping]
\label{pro:map}
Find a bijective and continuous invariant $I:\LIS\to\inv$ mapping the Lattice Isometry Space to a simpler space such that
\smallskip

\noindent
(\ref{pro:map}a) 
\emph{invariance} : $I(\La)$ is independent of a lattice basis and is preserved under isometry of $\R^2$, so $I$ has no false negatives : if $\La\cong\La'$ then $I(\La)=I(\La')$;
\smallskip

\noindent
(\ref{pro:map}b) 
\emph{completeness} : if $I(\La)=I(\La')$, then $\La,\La'$ are isometric, so $I$ has no false positives and defines a bijection (or a 1-1 map) $I:\LIS\to \inv=I(\LIS)$;
\smallskip

\noindent
(\ref{pro:map}c)
\emph{continuity} : the invariant 1-1 map $I:\LIS\to\inv$ is continuous in a suitable metric $d$ satisfying all axioms: (1) $d(\La,\La')=0$ \emph{if and only if} $\La\cong\La'$, (2) symmetry $d(\La,\La')=d(\La',\La)$, (3) triangle inequality $d(\La,\La')+d(\La',\La'')\geq d(\La,\La'')$;
\smallskip

\noindent
(\ref{pro:map}d)
\emph{computability} : the above metric $d(\La,\La')$ can be explicitly computed in a constant time from reduced bases of $\La,\La'$, see Definition~\ref{dfn:reduced_cell} in section~\ref{sec:review};
\smallskip

\noindent
(\ref{pro:map}e)
\emph{inverse design} : a basis of $\La$ can be explicitly reconstructed from $I(\La)$.
\bs 
\end{pro}
 
The metric axioms in~(\ref{pro:map}c) imply positivity due to $2d(\La,\La')=d(\La,\La')+d(\La',\La)\geq d(\La,\La)=0$. 
However, the identity of indiscernibles ($d(\La,\La')=0$ $\Lra$ $\La\cong\La'$) cannot be missed, otherwise even the zero function $d=0$ satisfies all other axioms. 
A binary answer to the isometry problem can provide only a discontinuous metric $d(\La,\La')$ equal to 1 or another positive number for any non-isometric lattices $\La\not\cong\La'$ even if $\La,\La'$ are nearly identical.
So the new continuity condition in (\ref{pro:map}c) makes Problem~\ref{pro:map} harder than a classification, especially up to rigid motion.
\medskip

\begin{figure}[h]
\includegraphics[width=1.0\textwidth]{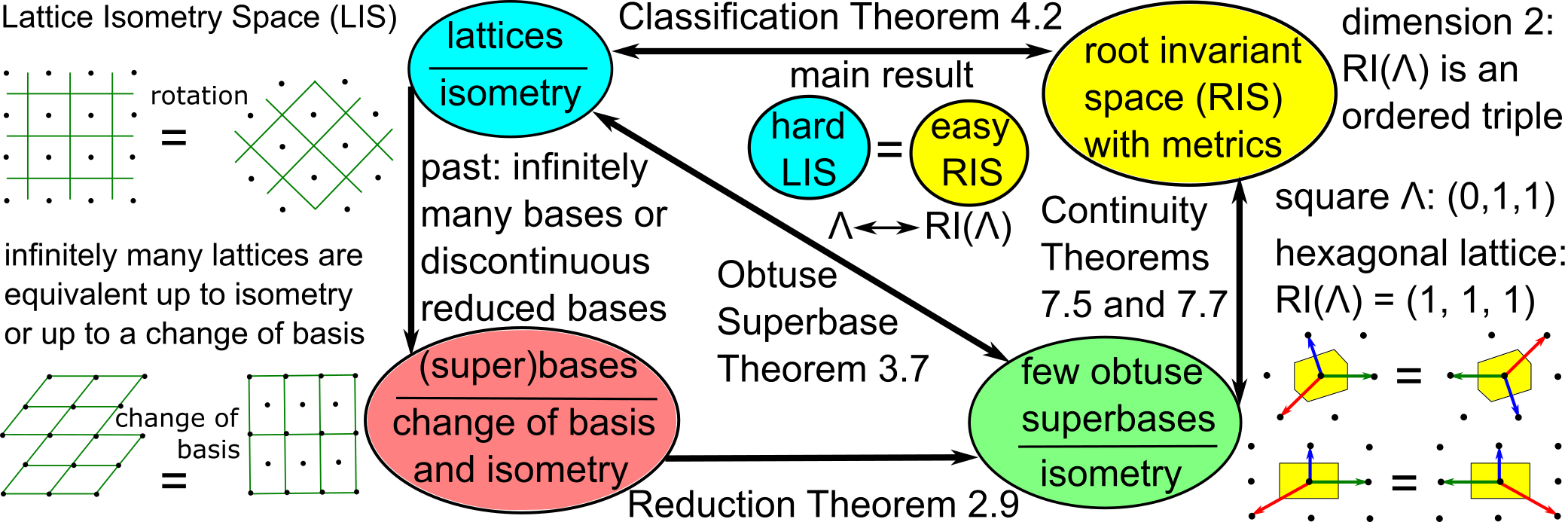}
\caption{
$\LIS(\R^2)$ is bijectively and continuously mapped to root invariants, which are ordered triples of square roots of scalar products of vectors of an obtuse superbase of a lattice $\La\subset\R^2$.}
\label{fig:lattice_classification}
\end{figure}

Fig.~\ref{fig:lattice_classification} summarises the past obstacles and a full solution to Problem~\ref{pro:map}.
The space $\inv$ will be the root invariant space $\RIS(\R^2)$ of ordered triples with continuous metrics.
Similar invariants will solve Problem~\ref{pro:map} up to three other equivalence relations: rigid motion, similarity and orientation-preserving similarity.

Lattices were previously represented by ambiguous or reduced bases, which are discontinuous under perturbations.
Most discrete invariants such as symmetry groups are also discontinuous and cut the Lattice Isometry Space ($\LIS$) into finitely many disjoint strata.
Delone \cite{delone1934mathematical}, later Conway and Sloane \cite{conway1992low} reduced ambiguity of lattice representations by using obtuse superbases.
Hence new continuous metrics and other structures on lattice spaces below are the next natural step.
\medskip

The inverse design in (\ref{pro:map}e) will raise Problem~\ref{pro:map} above metric geometry to define a richer structure of a vector space on $\LIS$.
It is easy to multiply any lattice by a fixed scalar, but a sum of any two lattices is harder to define in a meaningful way independent of lattice bases.
We will overcome this obstacle due to a linear structure on the root invariant space ($\RIS$) completely solving Problem~\ref{pro:map}.

\section{Main definitions and an overview of past work and new results}
\label{sec:review}

This section defines the main concepts and reviews past work on lattice comparisons, see the definition of an isometry and orientation in the appendix.
Any point $p$ in Euclidean space $\R^n$ can be represented by the vector from the origin $0\in\R^n$ to $p$.
This vector is also denoted by $p$,
An equal vector $p$ can be drawn at any initial point.
The \emph{Euclidean} distance between points $p,q\in\R^n$ is $|p-q|$.

\begin{dfn}[a lattice $\La$, a primitive unit cell $U(v_1,\dots,v_n)$]
\label{dfn:lattice_cell}
Let vectors $v_1,\dots,v_n$ form a linear {\em basis} in $\R^n$ so that 
any vector $v\in\R^n$ can be written as $v=\sum\limits_{i=1}^n c_i v_i$ for some real $c_i$, and if $\sum\limits_{i=1}^n c_i v_i=0$ then all $c_i=0$.
A {\em lattice} $\La$ in $\R^n$ consists of 
$\sum\limits_{i=1}^n c_i v_i$  with integer coefficients $c_i\in\Z$.
The parallelepiped $U(v_1,\dots,v_n)=\left\{ \sum\limits_{i=1}^n c_i v_i \vl c_i\in[0,1) \right\}$ is called a \emph{primitive unit cell} of $\La$.
\bs
\end{dfn}

The conditions $0\leq c_i<1$ on the coefficients $c_i$ above guarantee that  
the copies of unit cells $U(v_1,\dots,v_n)$ translated by all $v\in\La$ are disjoint and cover $\R^n$.

\begin{dfn}[orientation, isometry, rigid motion, similarity]
\label{dfn:isometry}
For a basis $v_1,\dots,v_n$ of $\R^n$, the \emph{signed volume} of $U(v_1,\dots,v_n)$ is the determinant of the $n\times n$ matrix with columns $v_1,\dots,v_n$.
The sign of this $\det(v_1,\dots,v_n)$ can be called an \emph{orientation} of the basis $v_1,\dots,v_n$. 
An \emph{isometry} is any map $f:\R^n\to\R^n$ such that $|f(p)-f(q)|=|p-q|$ for any $p,q\in\R^n$.
The unit cells $U(v_1,\dots,v_n)$ and $U(f(v_1),\dots,f(v_n))$ have non-zero volumes with equal absolute values.
If these volumes have equal signs, $f$ is \emph{orientation-preserving}, otherwise $f$ is \emph{orientation-reversing}.
Any orientation-preserving isometry $f$ is a composition of translations and rotations, and can be included into a continuous family of isometries $f_t$ (a \emph{rigid motion}), where $t\in[0,1]$, $f_0$ is the identity map and $f_1=f$.
A \emph{similarity} is a composition of isometry and uniform scaling $v\mapsto sv$ for a fixed scalar $s>0$. 
\bs
\end{dfn}

Any orientation-reversing isometry is a composition of a rigid motion and one reflection in a linear subspace of dimension $n-1$ (a line in $\R^2$). 
\medskip

Any lattice $\La$ can be generated by infinitely many bases or unit cells, see Fig.~\ref{fig:hexagonal_graphene}.
A standard approach to resolve this ambiguity is to consider a reduced basis below.
In $\R^3$, there are several ways to define a reduced basis \cite{gruber1989reduced}.
The most commonly used is Niggli's reduced cell \cite{niggli1928krystallographische}, whose 2-dimensional version is defined below.
\medskip

For vectors $v_1=(a_1,a_2)$ and $v_2=(b_1,b_2)$ in $\R^2$, let $\det(v_1,v_2)=a_1b_2-a_2b_1$ be the determinant of the matrix $\matfour{a_1}{b_1}{a_2}{b_2}$ with the columns $v_1,v_2$.

\begin{dfn}[reduced cell]
\label{dfn:reduced_cell}
For a lattice $\La\subset\R^2$ up to isometry, a basis and its unit cell $U(v_1,v_2)$ are called \emph{reduced} if $|v_1|\leq|v_2|$ and $-\frac{1}{2}v_1^2\leq v_1\cdot v_2\leq 0$.
Up to rigid motion, the conditions are weaker: $|v_1|\leq|v_2|$ and $-\frac{1}{2}v_1^2<v_1\cdot v_2\leq\frac{1}{2}v_1^2$, $\det(v_1,v_2)>0$, and the new \emph{special condition} : if $|v_1|=|v_2|$ then $v_1\cdot v_2\geq 0$.
\bs 
\end{dfn}

All bases marked in Fig.~\ref{fig:hexagonal_graphene} are reduced.
The condition $|v_1\cdot v_2|\leq\frac{1}{2}v_1^2$ in Definition~\ref{dfn:reduced_cell} geometrically means that $v_1,v_2$ are close to being orthogonal: the projection of $v_1$ to $v_2$ is between $\pm\frac{1}{2}|v_2|$.
The isometry conditions $|v_1|\leq|v_2|$ and $-\frac{1}{2}v_1^2\leq v_1\cdot v_2\leq 0$ in Definition~\ref{dfn:reduced_cell} coincide with
the conventional definition from \cite[section 9.2.2]{aroyo2013international} for type II (obtuse) cells in $\R^3$ if we choose $v_3$ to be very long and orthogonal to $v_1,v_2$.
If a basis $v_1,v_2$ satisfies the isometry conditions above, then so do three more bases obtained from $v_1,v_2$ by reflections in the lines parallel and orthogonal to $v_1$.
Up to rigid motion, we have maximum two mirror images of the same lattice $\La$, so the two extra bases are excluded due to $\det(v_1,v_2)>0$. 
If $\La$ is mirror-symmetric to itself due to $|v_1|=|v_2|$, one more basis is excluded by $v_1\cdot v_2\geq 0$.
Proposition~\ref{prop:reduced_bases}(a) proves uniqueness of a reduced basis in all cases. 

\begin{figure}[h]
\includegraphics[width=1.0\textwidth]{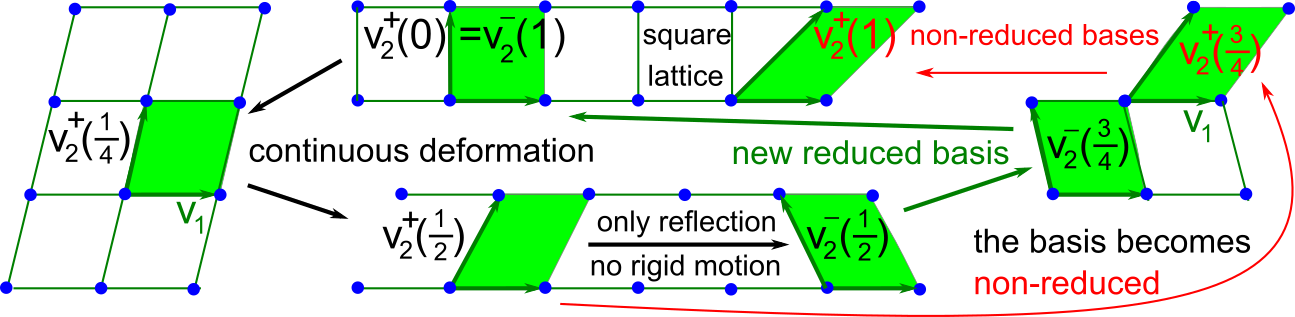}
\caption{
The basis $v_1=(1,0)$, $v_2^+(t)=(t,1)$ remains reduced only for $t\in[0,\frac{1}{2}]$ and at $t=\frac{1}{2}$ discontinuously changes to another reduced basis $v_1=(1,0)$, $v_2^-(t)=(t-1,1)$ for $t\in(\frac{1}{2},1]$.
}
\label{fig:basis_discontinuity}
\end{figure}

Fig.~\ref{fig:basis_discontinuity} shows the deformation of the lattice $\La(t)$ generated by $v_1=(1,0)$, $v_2^{+}(t)=(t,1)$ for $t\in[0,1]$.
The square lattice $\La(0)$ deforms to $\La(\frac{1}{2})$ in the bottom middle picture of Fig.~\ref{fig:basis_discontinuity}, where the previously reduced basis $v_1,v_2^+(t)=(t,1)$ becomes non-reduced for $t\in(\frac{1}{2},1]$ and substantially differs from the new reduced basis $v_1,v_2^-(t)=(t-1,1)$ at $t=\frac{1}{2}$.
In general, \cite[Theorem~15]{widdowson2022average} proved that there is no continuous reduction of bases if we compare them coordinate-wise.
Theorems~\ref{thm:OSI->RIS},~\ref{thm:RIS->OSI} and Corollary~\ref{cor:rect_discontinuity} will settle all (dis)continuity cases in $\R^2$.
\medskip

The book \cite{engel2004lattice} considered actions on lattices by groups with reflections.
In $\R^3$, crystallography classifies symmetry groups into 219 classes up to affine transformations including orientation-reversing maps, more often into 230 classes when orientation is preserved as by rigid motion of real crystals. 
So the classification of lattices up to rigid motion is more practically important than up to isometry.
\medskip

Another well-known cell of a lattice is the \emph{Voronoi domain} \cite{voronoi1908nouvelles}, also called the \emph{Wigner-Seitz cell}, \emph{Brillouin zone} or \emph{Dirichlet cell}.
We use the word \emph{domain} to avoid a confusion with a unit cell in Definition~\ref{dfn:lattice_cell}. 
Though the Voronoi domain can be defined for any point of a lattice, it suffices to consider only the origin $0$.

\begin{dfn}[Voronoi domain $V(\La)$]
\label{dfn:Voronoi_vectors}
The \emph{Voronoi domain} of a lattice $\La$ is the neighbourhood $V(\La)=\{p\in\R^n: |p|\leq|p-v| \text{ for any }v\in\La\}$ of the origin $0\in\La$ consisting of all points $p$ that are non-strictly closer to $0$ than to other points $v\in\La$.
A vector $v\in\La$ is called a \emph{Voronoi vector} if the bisector hyperspace $H(0,v)=\{p\in\R^n \vl p\cdot v=\frac{1}{2}v^2\}$ between 0 and $v$ intersects $V(\La)$.
If $V(\La)\cap H(0,v)$ is an $(n-1)$-dimensional face of $V(\La)$, then $v$ is called a \emph{strict} Voronoi vector. 
\bs
\end{dfn}

\begin{figure}[h]
\includegraphics[width=1.0\textwidth]{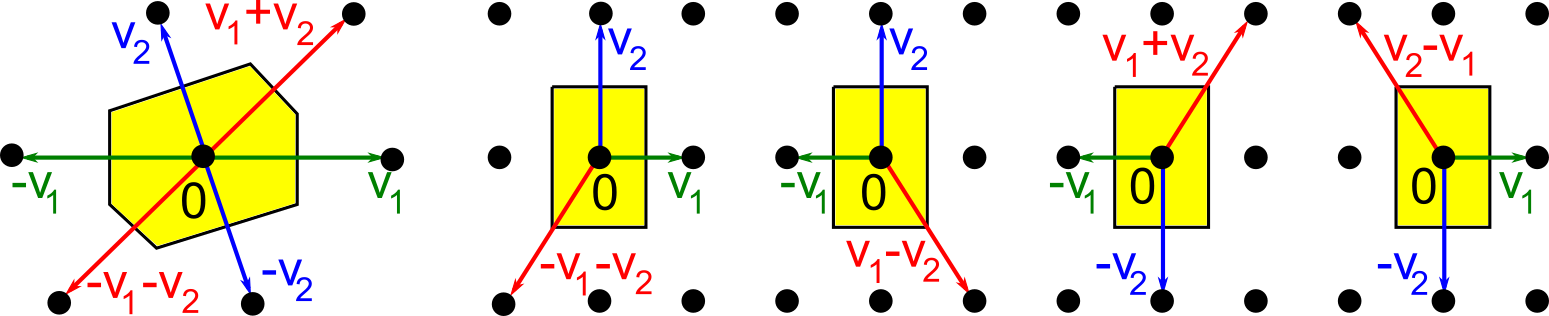}
\caption{\textbf{Left}: a generic lattice $\La\subset\R^2$ has a hexagonal Voronoi domain with an obtuse superbase $v_1,v_2,v_0=-v_1-v_2$, which is unique up to permutations and central symmetry.
\textbf{Other pictures}: two pairs of obtuse superbases (related by reflection) for a rectangular lattice.}
\label{fig:Voronoi2D}
\end{figure}

Fig.~\ref{fig:Voronoi2D} shows how the Voronoi domain $V(\La)$ can be obtained as the intersection of the closed half-spaces $S(0,v)=\{p\in\R^n \vl p\cdot v\leq\frac{1}{2}v^2\}$ whose boundaries $H(0,v)$ are bisectors between $0$ and all strict Voronoi vectors $v\in\La$.
A generic lattice $\La\subset\R^2$ has a hexagonal Vronoi domain $V(\La)$ with six Voronoi vectors.
\medskip

Any lattice is determined by its Voronoi domain by Lemma~\ref{lem:Voronoi_classification} in the appendix.
However, the combinatorial structure of $V(\La)$ is discontinuous under perturbations.
Almost any perturbation of a rectangular basis in $\R^2$ gives a non-rectangular basis generating a lattice whose Voronoi domain $V(\La)$ is hexagonal, not rectangular.
Hence any integer-valued descriptors of $V(\La)$ such as the numbers of vertices or edges are always discontinuous and unsuitable for continuous quantification of similarities between arbitrary crystals or periodic point sets.
\medskip

Optimal geometric matching of Voronoi domains with a shared centre led \cite{mosca2020voronoi} to two continuous metrics (up to orientation-preserving isometry and similarity) on lattices.
The minimisation over infinitely many rotations was implemented in \cite{mosca2020voronoi} by sampling and gave approximate algorithms for these metrics.
The complete invariant isoset \cite{anosova2021isometry} for periodic point sets in $\R^n$ has a continuous metric that can be approximated \cite{anosova2021introduction} with a factor $O(n)$.  
The metric on invariant density functions \cite{edels2021} required a minimisation over $\R$, so far without approximation guarantees.
\medskip

Lemma~\ref{lem:Voronoi_vectors} shows how to find all Voronoi vectors of any lattice $\La\subset\R^n$. 
The doubled lattice is $2\La=\{2v \vl v\in\La\}$.
Vectors $u,v\in\La$ are called \emph{$2\La$-equivalent} if $u-v\in 2\La$.
Then any vector $v\in\La$ generates its $2\La$-class $v+2\La=\{v+2u \vl u\in\La\}$, which is $2\La$ translated by $v$ and containing $-v$.
All classes of $2\La$-equivalent vectors form the quotient space $\La/2\La$.
Any 1-dimensional lattice $\La$ generated by a vector $v$ has the quotient $\La/2\La$ consisting of only two classes $\La$ and $v+\La$.  
 
\begin{lem}[criterion for Voronoi vectors {\cite[Theorem~2]{conway1992low}, \cite[p.~477]{conway2013sphere}}]
\label{lem:Voronoi_vectors}
For any lattice $\La\subset\R^n$, a non-zero vector $v\in\La$ is a Voronoi vector of $\La$ if and only if $v$ is a shortest vector in its $2\La$-class $v+2\La$.
Also, $v$ is a strict Voronoi vector if and only if $\pm v$ are the only shortest vectors in the $2\La$-class $v+2\La$.
\bt
\end{lem}

Appendix~\ref{sec:proofs} includes detailed proofs of key past results such as Lemma~\ref{lem:Voronoi_vectors}.
Any lattice $\La\subset\R^2$ generated by $v_1,v_2$ has $\La/2\La=\{v_1,v_2,v_1+v_2\}+\La$.
Notice that the vectors $v_1\pm v_2$ belong to the same $2\La$-class.
Assume that $v_1,v_2$ are not longer than $v_1+v_2$, which holds if $\angle(v_1,v_2)\in[60^{\circ},120^{\circ}]$.
If $v_1+v_2$ is shorter than $v_1-v_2$ as in Fig.~\ref{fig:Voronoi2D}~(left), then $\La$ has three pairs of strict Voronoi vectors $\pm v_1,\pm v_2,\pm(v_1+v_2)$.
If $v_1\pm v_2$ have the same length, the unit cell spanned by $v_1,v_2$ degenerates to a rectangle, $\La$ has four non-strict Voronoi vectors $\pm v_1\pm v_2$. 
\medskip

The triple of vector pairs $\pm v_1,\pm v_2,\mp(v_1+v_2)$ in Fig.~\ref{fig:Voronoi2D} motivates the concept of a superbase with the extra vector $v_0=-v_1-v_2$, which extends to any dimension $n$ by setting $v_0 = -\sum\limits_{i=1}^n v_n$.
For dimensions 2 and 3, Theorem~\ref{thm:reduction} will prove that any lattice has an obtuse superbase of vectors whose pairwise scalar products are non-positive and are called \emph{Selling parameters} \cite{selling1874ueber}.
For any superbase in $\R^n$, the negated parameters $p_{ij}=-v_i\cdot v_j$ can be interpreted as conorms of lattice characters, functions $\chi: \La\to\{\pm 1\}$ satisfying $\chi(u+v)=\chi(u)\chi(v)$), see \cite[Theorem~6]{conway1992low}.
So $p_{ij}$ will be defined as \emph{conorms} only for an obtuse superbase below. 

\begin{dfn}[obtuse superbase, conorms $p_{ij}$]
\label{dfn:conorms}
For any basis $v_1,\dots,v_n$ in $\R^n$, the \emph{superbase} $v_0,v_1,\dots,v_n$ includes the vector $v_0=-\sum\limits_{i=1}^n v_i$.
The \emph{conorms} $p_{ij}=-v_i\cdot v_j$ are the negative scalar products of the vectors above. 
The superbase is \emph{obtuse} if all conorms $p_{ij}\geq 0$, so all angles between vectors $v_i,v_j$ are non-acute for distinct indices $i,j\in\{0,1,\dots,n\}$.
The superbase is \emph{strict} if all $p_{ij}>0$.
\bs
\end{dfn}

\cite[formula (1)]{conway1992low} has a typo initially defining $p_{ij}$ as exact Selling parameters, but later Theorems 3,~7,~8 use the non-negative conorms $p_{ij}=-v_i\cdot v_j\geq 0$.
\medskip

The indices of a conorm $p_{ij}$ are distinct and unordered.
We set $p_{ij}=p_{ji}$ for all $i,j$.
For $n=1$, the 1-dimensional lattice generated by a vector $v_1$ has the obtuse superbase consisting of the two vectors $v_0=-v_1$ and $v_1$, so the only conorm $p_{01}=-v_0\cdot v_1=v_1^2$ is the squared length of $v_1$.
Any superbase of $\R^n$ has $\dfrac{n(n+1)}{2}$ conorms $p_{ij}$, for example, three conorms $p_{01},p_{02},p_{12}$ in dimension 2.

\begin{dfn}[partial sums $v_S$, vonorms $v_S^2$]
\label{dfn:vonorms}
Let a lattice $\La\subset\R^n$ have a superbase $B=\{v_0,v_1,\dots,v_n\}$. 
For any proper subset $S\subset\{0,1,\dots,n\}$ of indices,
 consider its complement $\bar S=\{0,1,\dots,n\}-S$ and the \emph{partial sum} $v_S=\sum\limits_{i\in S} v_i$ whose squared lengths $v_S^2$ are called the \emph{vonorms} of $B$ and can be expressed as $v_S^2=(\sum\limits_{i\in S} v_i)(-\sum\limits_{j\in\bar S}v_j)=-\sum\limits_{i\in S,j\in\bar S}v_{j}\cdot v_j=\sum\limits_{i\in S,j\in\bar S}p_{ij}$.
For $n=2$, we get
$$
v_0^2=p_{01}+p_{02},\qquad
v_1^2=p_{01}+p_{12},\qquad 
v_2^2=p_{02}+p_{12}.
\leqno{(\ref{dfn:vonorms}a)}$$
The above formulae allow us to express the conorms via vonorms as follows
$$
p_{12}=\dfrac{1}{2}(v_1^2+v_2^2-v_0^2),\quad
p_{01}=\dfrac{1}{2}(v_0^2+v_1^2-v_2^2),\quad
p_{02}=\dfrac{1}{2}(v_0^2+v_2^2-v_1^2).
\leqno{(\ref{dfn:vonorms}b)}$$
So $p_{ij}=\dfrac{1}{2}(v_i^2+v_j^2-v_k^2)$ for distinct $i,j\in\{0,1,2\}$ and $k=\{0,1,2\}-\{i,j\}$.
\bs
\end{dfn}

Lemma~\ref{lem:partial_sums} will later help to prove that a lattice is uniquely determined up to isometry by an obtuse superbase, hence by its vonorms or, equivalently, conorms.

\begin{lem}[Voronoi vectors $v_S$ {\cite[Theorem~3]{conway1992low}}]
\label{lem:partial_sums}
For any obtuse superbase $v_0,v_1,\dots,v_n$ of a lattice, all partial sums $v_S$ from Definition~\ref{dfn:vonorms} split into $2^n-1$ symmetric pairs $v_S=-v_{\bar S}$, which are Voronoi vectors representing distinct $2\La$-classes in $\La/2\La$.
All Voronoi vectors $v_S$ are strict if and only if all $p_{ij}>0$.
\bt
\end{lem}

By Conway and Sloane \cite[section~2]{conway1992low}, any lattice $\La\subset\R^n$ that has an obtuse superbase is called a \emph{lattice of Voronoi's first kind}.
It turns out that any lattice in dimensions 2 and 3 is of Voronoi's first kind by Theorem~\ref{thm:reduction}.

\begin{thm}[reduction to an obtuse superbase]
\label{thm:reduction}
Any lattice $\La$ in dimensions $2$ and $3$ has an obtuse superbase $\{v_0,v_1,\dots,v_n\}$ so that $v_0=-\sum\limits_{i=1}^n v_i$ and all conorms $p_{ij}=-v_i\cdot v_j\geq 0$ for all distinct indices $i,j\in\{0,1,\dots,n\}$.
\bt
\end{thm}

Conway and Sloane in \cite[section~7]{conway1992low} attempted to prove Theorem~\ref{thm:reduction} for $n=3$ by example, which is corrected in \cite{kurlin2022complete}.
Appendix~\ref{sec:proofs} proves Theorem~\ref{thm:reduction} for $n=2$.
Finding an obtuse superbase is related to solving the shortest vector problem in a lattice.
The latter problem is NP-hard \cite{ajtai1998shortest} in $\R^n$, see the great review in \cite{nguyen2010lll}.
\medskip

The following result implies that the space of all lattices and general periodic point sets in any $\R^n$ is continuous, path-connected in the language of topology.
Due to Proposition~\ref{prop:cont_deform}, if we call any periodic structures equivalent (or similar) when they differ up to any small perturbation of points, then any two point sets become equivalent by the transitivity axiom: if $A_1\sim A_2\sim\dots\sim A_k$ then $A_1\sim A_k$. 

\begin{prop}
\label{prop:cont_deform}
Any periodic point sets $\La+M=\{v+p \mid v\in \La,\, p\in M\}$, where $\La\subset\R^n$ is a lattice, $M$ is a finite set (motif) of points in a unit cell of $\La$, can be deformed into each other so that coordinates of all points change continuously.
\bs
\end{prop}
\begin{proof}
Starting from any basis $v_1,\dots,v_n$, continuously rotate $v_2,\dots,v_n$ to make all basis vectors pairwise orthogonal.
If given periodic sets have different numbers $m_1\neq m_2$ of points in their unit cells, we can enlarge their cells in the direction of $v_1$ by factors $m_2,m_1$ so that both sets have the same number of $m=m_1m_2$ points in their cells.  
If the coordinates of points $p\in M$ remain constant in the moving basis $v_1,\dots,v_n$, they change continuously in a fixed basis of $\R^n$.
Then we can continuously elongate $v_1,\dots,v_n$ to a get a sufficiently large unit cubic cell.
After two periodic point sets are put in this `gas state' in a common large cube, continuously move all points from one motif into any other configuration without collisions.
A composition of the above movements connects any periodic sets.
\ws
\end{proof}

Any lattice $\La\subset\R^2$ with a basis $v_1,v_2$ defines the \emph{positive quadratic form} 
$$Q(x,y)=(xv_1+yv_2)^2=q_{11}x^2+2q_{12}xy+q_{22}y^2\geq 0 
\text{ for all }x,y\in\R,$$ 
where $q_{11}=v_1^2$, $q_{22}=v_2^2$, $q_{12}=v_1\cdot v_2$.
Changing the basis $v_1,v_2$ is equivalent to replacing $x,y$ by the linear combinations of the coordinates of the vector $xv_1+yv_2$ in a new basis.
Conversely, any positive quadratic form $Q(x,y)$ can be written as a sum of two squares $(a_1x+b_1y)^2+(a_2x+b_2y)^2$, see \cite[Theorem~2 on p.~116]{delone1975bravais}, and defines the lattice with the basis $v_1=(a_1,a_2)$, $v_2=(b_1,b_2)$.
\medskip

In 1773 Lagrange \cite{lagrange1773recherches} proved that any positive quadratic form can be rewritten so that $0<q_{11}\leq q_{22}$ and $-q_{11}\leq 2q_{12}\leq 0$. 
The resulting \emph{reduced} basis $v_1,v_2$ satisfies $0<v_1^2\leq v_2^2$ and $-v_1^2\leq 2v_1\cdot v_2\leq 0$ without the new special conditions in Definition~\ref{dfn:reduced_cell}.
Then the lattice generated by $v_1=(1,0),v_2^{\pm}(\frac{1}{2})=(\pm\frac{1}{2},1)$ in Fig.~\ref{fig:basis_discontinuity} (bottom middle) is represented by two Lagrange-reduced forms $x^2\pm xy+\frac{5}{4}y^2$.
\medskip

To resolve the basis ambiguity, Delone defined the parameters \cite[section~29]{delone1938geometry} equal to the conorms from Definition~\ref{dfn:conorms}.
We use the notations $p_{ij}$ from \cite{conway1992low}:
$$\begin{array}{l}
p_{01}=q_{11}+q_{12}=v_1^2+v_1\cdot v_2=v_1\cdot(v_1+v_2)=-v_0\cdot v_1,\\ 
p_{02}=q_{22}+q_{12}=v_2^2+v_1\cdot v_2=v_2\cdot(v_1+v_2)=-v_0\cdot v_2,\\
p_{12}=-q_{12}=-v_1\cdot v_2. \end{array}$$ 
The quadratic form becomes a sum of squares: $Q_\La=p_{01}x^2+p_{22}y^2+p_{12}(x-y)^2$.
The inequalities for $q_{ij}$ are equivalent to the simple ordering $0\leq p_{12}\leq p_{01}\leq p_{02}$, which Definition~\ref{dfn:RI} will use to introduce a more convenient  root invariant.
\medskip

The isometry classification in 
(\ref{pro:map}ab) can be interpreted via group actions, see \cite{engel2004lattice} and \cite{zhilinskii2016introduction}.
Let $\mathcal{B}_n$ be the space of all linear bases in $\R^n$.
Up to a change of basis, all lattices in $\R^n$ form the $n^2$-dimensional orbit space $\mathcal{L}_n=\mathcal{B}_n/\GL_n(\R)$, see \cite[formula (1.37) on p.~34]{engel2004lattice}.
Up to orthogonal maps from the group $\Or_n(\R)$, the orbit space of lattices can be identified with the cone $\mathcal{C}_+(\mathcal{Q}_n)=\mathcal{B}_n/O_n(\R)$ of positive quadratic forms, where $\mathcal{Q}_n$ denotes the space of real symmetric $n\times n$ matrices, see \cite[formula (1.67) on p.~41]{engel2004lattice}.
The Lattice Isometry Space $\LIS(\R^n)$ was called the space of \emph{intrinsic} lattices $\mathcal{L}_n^o=\mathcal{C}_+(\mathcal{Q}_n)/\GL_n(\Z)$ in \cite[formula (1.70) on p.~42]{engel2004lattice}.
\medskip

The past approach to uniquely identify an intrinsic lattice (isometry class), say for $n=2$,  was to choose a fundamental domain of the action of $\GL_2(\Z)$ on the cone $\mathcal{C}_+(\mathcal{Q}_2)$.
This choice is equivalent to a choice of a reduced basis, which can be discontinuous.
Mirror reflections of any lattice $\La$ correspond to quadratic forms $q_{11}x^2\pm 2q_{12}xy+q_{22}y^2$ that differ by a sign of $q_{12}$.
To distinguish mirror images of lattices, Definition~\ref{dfn:sign} will introduce $\sign(\La)$.
Then continuous deformations of lattices become continuous paths in a space of invariants, see Remark~\ref{rem:group_action}.

\begin{figure}[h]
\includegraphics[width=\textwidth]{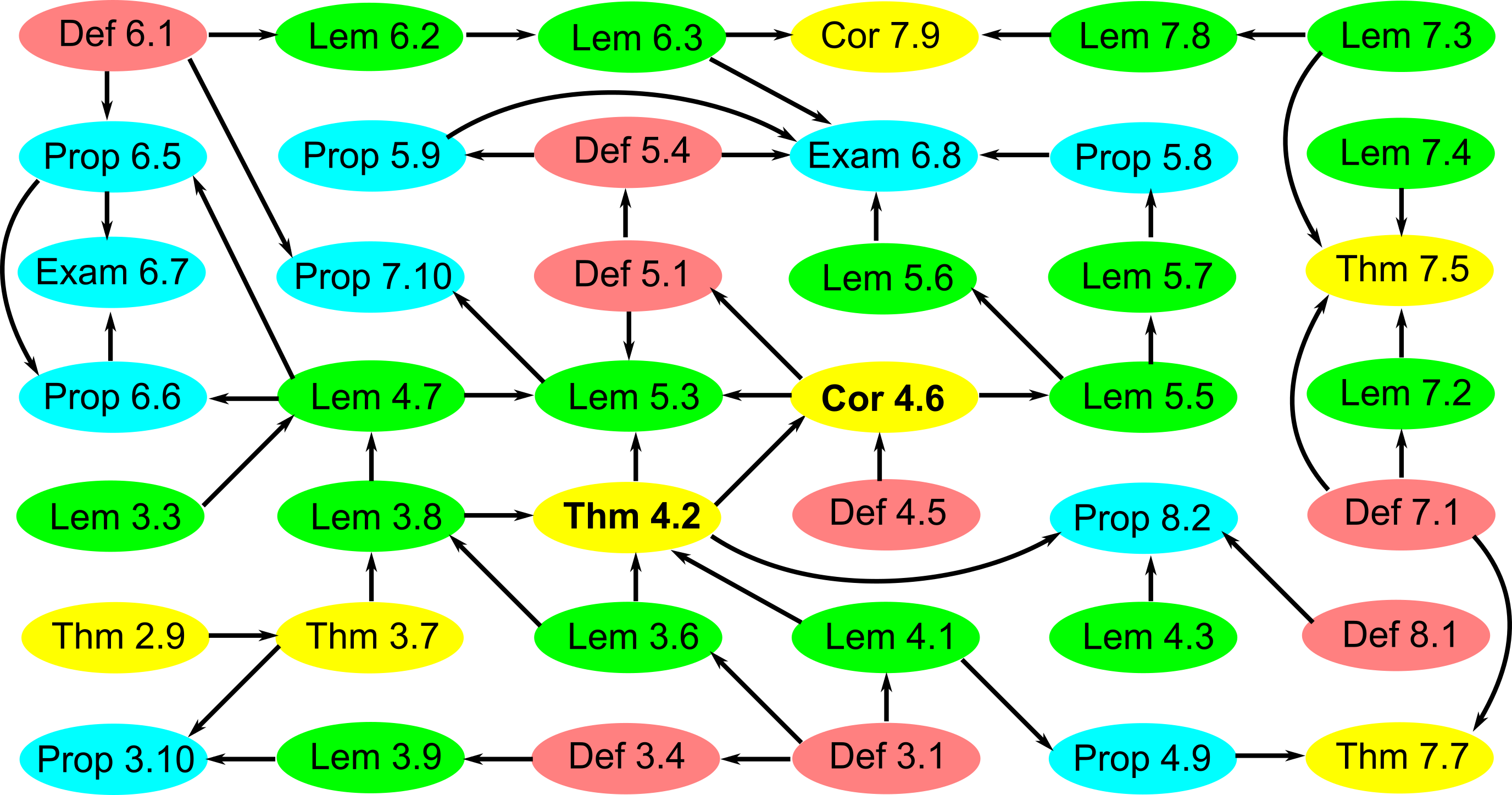}
\caption{Major logical connections between new definitions, auxiliary lemmas and main results.
}
\label{fig:lattices2d_logic}
\end{figure}

Proposition~\ref{prop:reduced_bases} establishes a 1-1 correspondence between obtuse superbases and reduced bases.
The latter bases are common in crystallography and implemented by many fast algorithms \cite{aroyo2011crystallography}.
So our lattice input will be any obtuse superbase.
Fig.~\ref{fig:lattices2d_logic} shows logical flows from key concepts to new contributions.
\medskip

The main results are complete classifications in Theorem~\ref{thm:classification2d}, Corollary~\ref{cor:similar_lattices2d}, and metrics on lattice invariants in Definitions~\ref{dfn:RM},~\ref{dfn:RMo}.
Continuity of invariants in Theorems~\ref{thm:OSI->RIS},~\ref{thm:RIS->OSI} convert the Lattice Isometry Space $\LIS(\R^2)$ into a continuously parameterised map solving Problem~\ref{pro:map}.
Definition~\ref{dfn:RC} extends the binary chirality to a real-valued deviation of a lattice from a higher-symmetry neighbour.
\medskip

Petitjean \cite{petitjean2003chirality} comprehensively described past approaches to quantify chirality of bounded objects such as a rigid molecule.
The most rigorous approach is to use a metric between such rigid objects. 
However, even for the simplest case of a finite set of points, the Hausdorff-like distances between finite sets required approximate minimisations over infinitely many rotations.
Definition~\ref{dfn:RC} will introduce chiral distances for 2D lattices, which are easily computable by Propositions~\ref{prop:RC},~\ref{prop:PC}.

\section{Isometry invariants of an obtuse superbase of a 2-dimensional lattice}
\label{sec:forms2d}

Definition~\ref{dfn:RI} introduces voforms $\VF$ and coforms $\CF$, which are triangular cycles whose three nodes are marked by vonorms and conorms, respectively.
We start from any obtuse superbase $B$ of a lattice $\La\subset\R^2$ to define $\VF$, $\CF$, and a root invariant $\RI$. Lemma~\ref{lem:lattice_invariants}(a) will justify that $\RI$ depends only on $\La$, not on $B$.

\begin{figure}[h]
\includegraphics[height=35mm]{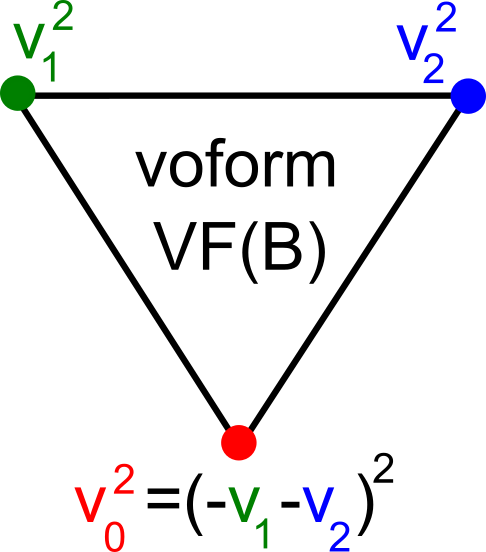}
\includegraphics[height=35mm]{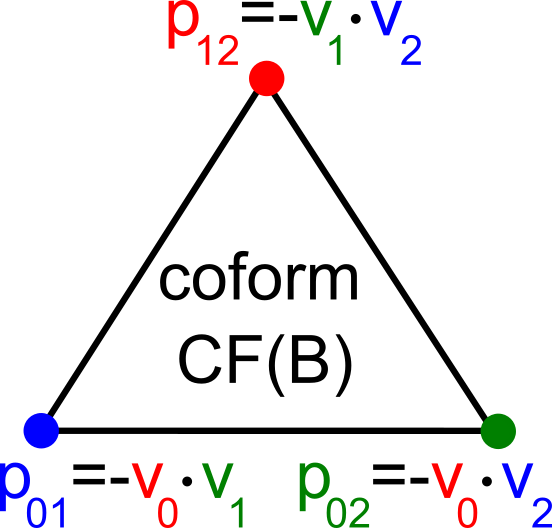}
\hspace*{2mm}
\includegraphics[height=35mm]{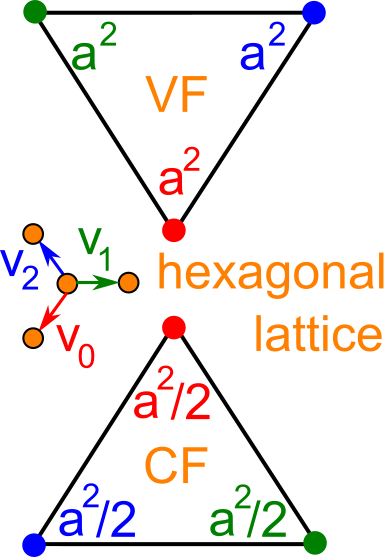}
\hspace*{2mm}
\includegraphics[height=35mm]{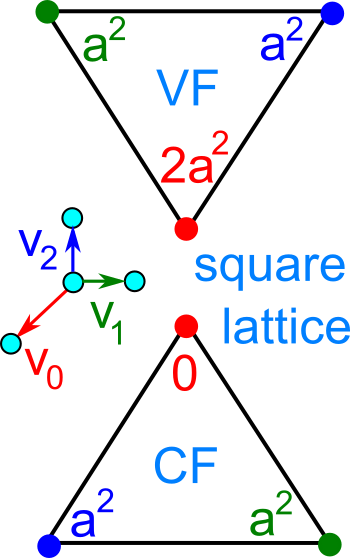}
\caption{\textbf{1st picture}: a voform $\VF(B)$ of a 2D lattice with an obtuse superbase $B=\{v_0,v_1,v_2\}$.
\textbf{2nd picture}: nodes of a coform $\CF(B)$ are marked by conorms $p_{ij}$.
\textbf{3rd and 4th pictures}: $\VF$ and $\CF$ of the hexagonal and square lattice with a minimum inter-point distance $a$.}
\label{fig:forms2d}
\end{figure}

\begin{dfn}[voform $\VF$, coform $\CF$, ordered \emph{root invariant} $\RI$]
\label{dfn:RI}
For any ordered obtuse superbase $B$ in $\R^2$, the \emph{voform} $\VF(B)$ is the cycle on three nodes marked by the vonorms $v_0^2,v_1^2,v_2^2$, see Fig.~\ref{fig:forms2d}.
The \emph{coform} $\CF(B)$ is the cycle on three nodes marked by the conorms $p_{12},p_{02},p_{01}$.
Since all conorms $p_{ij}\geq 0$, we can define the \emph{root products} $r_{ij}=\sqrt{p_{ij}}$.
The \emph{root invariant} $\RI(B)$ is obtained by writing the three root products $r_{12},r_{01},r_{02}$ in increasing order.
\bs
\end{dfn}

The ordering $r_{12}\leq r_{01}\leq r_{02}$ is equivalent to $v_1^2\leq v_2^2\leq v_0^2$ by formulae~(\ref{dfn:vonorms}a).
Root products have the same units as original coordinates of basis vectors, for example, Angstroms: $1\AA=10^{-10}$m.
The ordered root invariant $\RI(B)$ is more convenient than $\VF(B)$ and $\CF(B)$, which depend on an order of vectors of $B$.

\begin{figure}[h]
\includegraphics[width=1.0\textwidth]{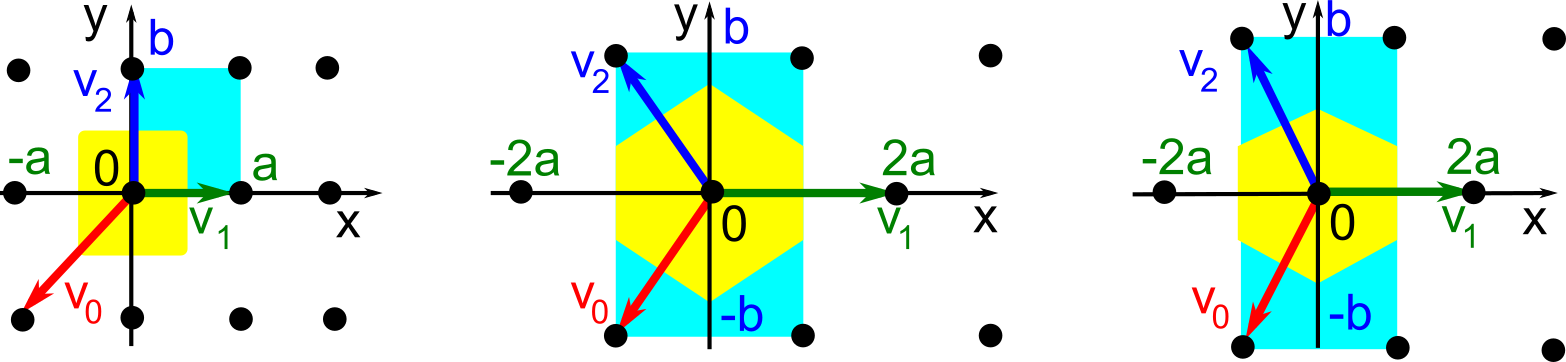}
\caption{
\textbf{Left}: $\La$ has a rectangular cell and obtuse superbase $B$ with $v_1=(a,0)$, $v_2=(0,b)$, $v_0=(-a,-b)$, see Example~\ref{exa:achiral_lattices} and Lemma~\ref{lem:achiral_lattices}.
Other lattices $\La$ have a rectangular cell $2a\times 2b$ and an obtuse superbase $B$ with $v_1=(2a,0)$, $v_2=(-a,b)$, $v_0=(-a,-b)$.
\textbf{Middle}: $\RI(B)=(\sqrt{b^2-a^2}, a\sqrt{2},a\sqrt{2})$, $a\leq b\leq a\sqrt{3}$.
\textbf{Right}: $\RI(B)=(a\sqrt{2},a\sqrt{2},\sqrt{b^2-a^2})$, $a\sqrt{3}\leq b$.
}
\label{fig:achiral_lattices}
\end{figure}

\begin{exa}
\label{exa:achiral_lattices}
\textbf{(a)} 
A lattice $\La$ with a rectangular cell of sides $a\leq b$ has an obtuse superbase $B$ with $v_1=(a,0)$, $v_2=(0,b)$, $v_0=(-a,-b)$, and $\RI(B)=(0,a,b)$.
\medskip

\noindent
\textbf{(b)} 
For any lattice $\La\subset\R^2$ whose Voronoi domain $V(\La)$ is a mirror-symmetric hexagon, assume that the $x$-axis is its line of symmetry.
Since $V(\La)$ is centrally symmetric with respect to the origin $0$, the $y$-axis is also its line of symmetry, see Fig.~\ref{fig:achiral_lattices}.
Then $\La$ has the centred rectangular (non-primitive) cell with sides $2a\leq 2b$.
The obtuse superbase $B$ with $v_1=(2a,0)$, $v_2=(-a,b)$, $v_0=(-a,-b)$ has $\RI(B)=(a\sqrt{2},a\sqrt{2},\sqrt{b^2-a^2})$ for $b\geq a\sqrt{3}$.
For $a\leq b<a\sqrt{3}$, we should swap $r_{02}=\sqrt{b^2-a^2}$ with $r_{12}=a\sqrt{2}$ to get an ordered root invariant $\RI(B)$.
\bs
\end{exa}

A lattice $\La\subset\R^n$ that can be mapped to itself by a mirror reflection with respect to a $(n-1)$-dimensional hyperspace can be called \emph{mirror-symmetric} or \emph{achiral}.
Since a mirror reflection of any lattice $\La\subset\R^2$ with respect to a line $L\subset\R^2$ can be realised by a rotation in $\R^3$ around $L$ through $180^{\circ}$, the term \emph{achiral} sometimes applies to all 2D lattices and becomes non-trivial only for 3D lattices.
This paper for 2D lattices uses the clearer adjective \emph{mirror-symmetric}.
 
\begin{lem}[root invariants of mirror-symmetric lattices $\La\subset\R^2$]
\label{lem:achiral_lattices}
An obtuse superbase $B$ generates a \emph{mirror-symmetric} lattice $\La(B)$ \emph{if and only if} 
\smallskip

\noindent
(\ref{lem:achiral_lattices}a) 
the root invariant $\RI(B)$ contains a zero value and $\La(B)$ is rectangular, or
\smallskip

\noindent
(\ref{lem:achiral_lattices}b) 
$\RI(B)$ has equal root products and the Voronoi domain of $\La(B)$ is a square or a hexagon whose symmetry group has two orthogonal axes of symmetry.
\bt
\end{lem}
\begin{proof}
The part \emph{if} $\Leftarrow$.
Let $\RI(B)$ include a zero, which should be the first root product, say $0=r_{12}=\sqrt{-v_1\cdot v_2}$.
The vectors $v_1,v_2$ of the superbase $B$ are orthogonal and generate a rectangular lattice, which is mirror-symmetric.
If $\RI(B)$ has two equal root products, say $r_{01}=r_{02}$, the conorms are also equal: $p_{01}=p_{02}$.
Formulae~(\ref{dfn:vonorms}a) imply that $v_1^2=p_{01}+p_{12}=p_{02}+p_{12}=v_2^2$.
The vectors $v_1,v_2$ have equal lengths and can be swapped ($v_1\lra v_2$) by the reflection in the bisector $L$ between $v_1,v_2$, which preserves $v_0=-v_1-v_2$, so $\La(B)$ is mirror-symmetric.
\medskip

The part \emph{only if} $\Rightarrow$.
If $\La(B)$ is mirror-symmetric, then so is its Voronoi domain $V(\La)$.
If $V(\La)$ is a rectangle or a mirror-symmetric hexagon as in Fig.~\ref{fig:achiral_lattices}, $\RI(B)$ computed in Example~\ref{exa:achiral_lattices} contains either a zero or two equal root products. 
\ws
\end{proof}


\begin{dfn}[$\sign(B)$, the oriented root invariant $\RI^o(B)$]
\label{dfn:sign}
If an obtuse superbase $B$ generates a mirror-symmetric lattice, set $\sign(B)=0$.
Else all vectors of $B$ have different lengths and angles not equal to $90^{\circ}$ by Lemma~\ref{lem:achiral_lattices}.
Let $v_1,v_2$ be the shortest vectors of $B$ so that $|v_1|<|v_2|$.
Then $\sign(B)=\pm 1$ is the sign of the determinant $\det(v_1,v_2)$ of the matrix with the columns $v_1,v_2$.
The \emph{oriented root invariant} $\RI^o(B)$ is obtained by adding $\sign(B)$ as a superscript to $\RI(B)$.
\bs
\end{dfn}

\begin{figure}[h]
\label{fig:root_forms2d_reflection}
\caption{The lattices $\La,\La'$ are mirror reflections of each other and have oriented root invariants $\RI^o=(\sqrt{3},\sqrt{6},\sqrt{7})_{\pm}$ with opposite signs introduced in Definition~\ref{dfn:sign}, see Example~\ref{exa:signs}.} 
\includegraphics[width=\textwidth]{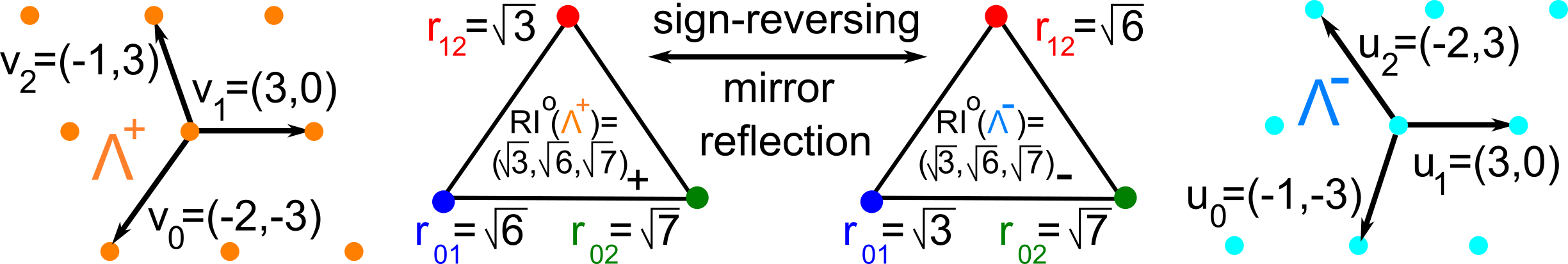}
\end{figure}

If $\sign(B)=0$, this zero superscript in $\RI^o(B)$ can be skipped for simplicity, so $\RI^o(B)=\RI(B)$ in this case.
Theorem~\ref{thm:isometric_superbases} will show that $\sign(B)$ can be considered as an invariant of a lattice $\La$ up to orientation-preserving similarity.
\medskip

In Definition~\ref{dfn:sign} the determinant $\det(v_1,v_2)$ is the signed area of the unit cell $U(v_1,v_2)$ equal to $|v_1|\cdot|v_2|\sin\angle(v_1,v_2)$, where the angle is measured from $v_1$ to $v_2$ in the anticlockwise direction around the origin $0\in\R^2$.
For a strict obtuse superbase $B$, all angles between its basis vectors are strictly obtuse.
Then $\sign(B)=+1$ if $\angle(v_1,v_2)$ is in the positive range $(90^\circ,180^\circ)$, else $\sign(B)=-1$.

\begin{exa}[signs of lattices]
\label{exa:signs}
The lattice $\La^+$ in the first picture of Fig.~\ref{fig:root_forms2d_reflection} has the obtuse superbase $B$ with $v_1=(3,0)$, $v_2=(-1,3)$, $v_0=(-2,-3)$ of lengths $3,\sqrt{10},\sqrt{13}$, respectively, so $\La^+$ is not mirror-symmetric.
Since $v_1,v_2$ are the two shortest vectors of $B^+$ and $\det(v_1,v_2)=\det\matfour{3}{-1}{0}{3}>0$, we get $\sign(B^+)=+1$.
The anticlockwise angle is $\angle(v_1,v_2)=180^\circ-\arcsin\frac{3}{\sqrt{10}}\approx 108^\circ$.
\medskip

The lattice $\La^-$ in the last picture of Fig.~\ref{fig:root_forms2d_reflection} is obtained from $\La^+$ by a mirror reflection and has the obtuse superbase $B^-$ with $u_1=v_1$, $u_2=(-2,3)$, $u_0=(-1,-3)$ of lengths
$3,\sqrt{13},\sqrt{10}$, respectively, so $\La^-$ is not mirror-symmetric. 
Since $u_1,u_0$ are the shortest vectors, $\det(u_1,u_0)=\det\matfour{3}{-1}{0}{-3}<0$, we get $\sign(B^-)=-1$.
The anticlockwise angle is $\angle(u_1,u_0)=\arcsin\frac{3}{\sqrt{10}}-180^\circ\approx -108^\circ$.
\bs
\end{exa}

\begin{lem}[$\RI$ invariance]
\label{lem:RI_invariance} 
For an unordered obtuse superbase $B$ in $\R^2$, any isometry preserves $\RI(B)$.
Any rigid motion preserves $\sign(B)$ and $\RI^o(B)$.
\bt
\end{lem}
\begin{proof}
Any isometry of an ordered obtuse superbase $B$ preserves the lengths and scalar products of the ordered vectors, so $\RI(B)$ is unchanged.
Any re-ordering of vectors of $B$ permutes conorms. 
$\RI(B)$ is unique due to ordered root products. 
\medskip

If a lattice is mirror-symmetric, then so is its image under any rigid motion in $\R^2$, hence $\sign(B)=0$ is preserved.
If $B$ generates a non-mirror symmetric lattice, $B$ has unique shortest vectors $v_1,v_2$.
A rigid motion acts on $v_1,v_2$ as a special orthogonal matrix with determinant 1, hence preserving $\det(v_1,v_2)$, $\sign(B)$.
\ws 
\end{proof}

Theorem~\ref{thm:isometric_superbases}  below is crucial for a complete classification of 2D lattices in Theorem~\ref{thm:classification2d} and Corollary~\ref{cor:similar_lattices2d}. 
Theorem~\ref{thm:isometric_superbases} highlights that mirror-symmetric lattices have more options for obtuse superbases up to rigid motion.
The same rectangular lattice can have two obtuse bases with $v_1=(1,0)$, $v_2=(0,\pm 2)$, which are related by reflection in the $x$-axis, not by rigid motion.
This symmetry-related ambiguity is much harder to resolve for 3D lattices even up to isometry, see \cite{kurlin2022complete}.

\begin{thm}[isometric obtuse superbases]
\label{thm:isometric_superbases}  
Any lattices $\La,\La'\subset\R^2$ are isometric if and only if any obtuse superbases of $\La,\La'$ are isometric.
If $\La,\La'$ are not rectangular, the same conclusion holds for rigid motion instead of isometry.
Any rectangular (non-square) lattice has two obtuse superbases related by reflection.
\bt
\end{thm}
\begin{proof}
Part \emph{if} ($\Leftarrow$): any isometry between obtuse superbases of $\La,\La'$ linearly extends to an isometry $\La\to\La'$.
Part \emph{only if} ($\Rightarrow$) means that any obtuse superbase of $\La$ is unique up to isometry.
By Lemma~\ref{lem:partial_sums} for $n=2$, if a lattice $\La$ has a strict obtuse superbase $B=\{v_0,v_1,v_2\}$, the Voronoi vectors of $\La$ are the pairs of opposite partial sums $\pm v_0,\pm v_1,\pm v_2$, see Fig.~\ref{fig:Voronoi2D}~(left). 
Hence $B$ is uniquely determined by the strict Voronoi vectors up to a sign.
So $B$ is one of only two obtuse superbases $\pm\{v_0,v_1,v_2\}$ related by central symmetry or rotation through $180^{\circ}$ around $0$.
Hence $\La$ has a unique obtuse superbase up to rigid motion. 
\medskip

If a superbase of $\La$ is non-strict, one conorm vanishes, say $p_{12}=0$.
Then $v_1,v_2$ span a rectangular unit cell and $\La$ has four non-strict Voronoi vectors $\pm v_1\pm v_2$ with all possible combinations of signs.
Hence $\La$ has four obtuse superbases 
$\{v_1,v_2,-v_1-v_2\}$, 
$\{-v_1,-v_2,v_1+v_2\}$,
$\{-v_1,v_2,v_1-v_2\}$, 
$\{v_1,-v_2,v_2-v_1\}$, see Fig.~\ref{fig:Voronoi2D}.
The first two (and the last two) superbases are obtained from each other by rotation through $180^\circ$ around the origin. 
Unless the lattice is square, the resulting two classes of superbases are related by reflection, not by rigid motion.
\ws
\end{proof}

\begin{lem}[lattice invariants]
\label{lem:lattice_invariants}
\textbf{(a)}
For any obtuse superbase $B$ of a lattice $\La\subset\R^2$,
 the root  invariant $\RI(B)$ is an isometry invariant of $\La$ and can be denoted by $\RI(\La)$.
Similarly, $\RI^o(\La)$ and $\sign(\La)$ are invariants of a lattice $\La$ up to rigid motion and orientation-preserving similarity, respectively.
\medskip

\noindent
\textbf{(b)}
A lattice $\La\subset\R^2$ is mirror-symmetric if and only if $\sign(\La)=0$.
\bt
\end{lem}
\begin{proof}
\textbf{(a)}
An obtuse superbase $B$ of any lattice $\La$ is unique up to isometry by Theorem~\ref{thm:isometric_superbases}.
Lemma~\ref{lem:RI_invariance} implies that the root invariant $\RI$ is an isometry invariant of $\La$, independent of any obtuse superbase $B$, hence can be denoted by $\RI(\La)$.
\medskip

Since an obtuse superbase $B$ of any non-mirror-symmetric lattice $\La$ is unique up to rigid motion by part \textbf{(a)}, Lemma~\ref{lem:RI_invariance} implies that $\sign(B)$ and $\RI^o(B)$ are invariant up to rigid motion, hence can be denoted by $\sign(\La)$ and $\RI^o(\La)$, respectively.
If $\La$ is mirror-symmetric, then any rigid motion preserves $\sign(\La)=0$ as well as $\RI(\La)$.
So $\RI^o(\La)$ is invariant up to rigid motion for all $\La\subset\R^2$.
\medskip

Any orientation-preserving similarity is a composition of a rigid motion and a uniform scaling (or a dilation) of all vectors by a factor $s>0$.
This similarity preserves any symmetries of the lattice $\La$ and multiplies the determinant $\det(v_1,v_2)$ from Definition~\ref{dfn:sign} by $s^2>0$, hence preserving $\sign(\La)$. 
\medskip

\noindent
\textbf{(b)}
By Definition~\ref{dfn:sign} any mirror-symmetric lattice has $\sign(\La)=0$.
Any basis $v_1,v_2$ of a non-mirror symmetric lattice $\La$ has $\det(v_1,v_2)\neq 0$, so $\sign(\La)\neq 0$.
\ws
\end{proof}

For any lattices $\La^{\pm}$ related by reflection, their unoriented root invariants $\RI(\La^{\pm})$ are identical, while 
$\RI^o(\La^{\pm})$ differ by sign.
Lemma~\ref{lem:signs} computes $\sign(\La)$ from any obtuse superbase whose first vector can be assumed to be $v_1=(1,0)$.

\begin{lem}[geometry of signs]
\label{lem:signs}
\textbf{(a)}
Up to orientation-preserving similarity, any lattice $\La\subset\R^2$ has an obtuse superbase with $v_1=(1,0)$ and $v_2=(x,y)$, where $(x,y)$ belongs to the region $\Obt=\{-1\leq x\leq 0<y,\; x^2+x+y^2\geq 0\}$. 
Then $\sign(\La)$ from Definition~\ref{dfn:sign} is determined by $(x,y)$ in Fig.~\ref{fig:signs}, see also Table~\ref{tab:signs}.
\medskip

\begin{figure}[h]
\includegraphics[width=\textwidth]{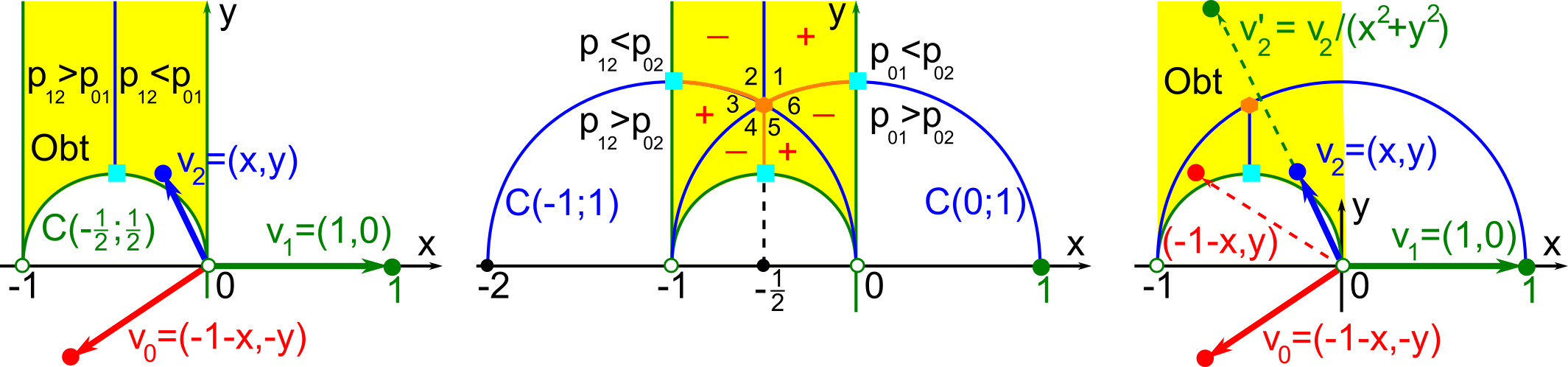}
\caption{\textbf{Left}: if $\La$ has an obtuse superbase with $v_1=(1,0)$, $v_2=(x,y)$, then $\sign(\La)$ is determined by $(x,y)\in\Obt$ above the circle $C(-\frac{1}{2};\frac{1}{2})=\{(x,y)\in\R^2\mid x^2+x+y^2=0\}$.
\textbf{Middle}: $\Obt$ splits into six subregions by the vertical line $x=-\frac{1}{2}$ and the circles $C(0;1)=\{(x,y)\in\R^2\mid x^2+y^2=1\}$ and $C(-1;1)=\{(x,y)\in\R^2\mid x^2+2x+y^2=0\}$, see Lemma~\ref{lem:signs}(a) and Table~\ref{tab:signs}.
\textbf{Right}: re-ordering and re-scaling vectors of an obtuse superbase is realised by the symmetries acting on $v_2=(x,y)$ within the yellow region $\Obt$, see Lemma~\ref{lem:signs}(b).
} 
\label{fig:signs}
\end{figure}

\noindent
\textbf{(b)}
Up to similarity, any lattice $\La\subset\R^2$ with an obtuse superbase $v_0,v_1,v_2$ can be represented by up to six points $(x,y)$ in the subregions of $\Obt$.
Swapping $v_0\lra v_2$ is realised by the reflection in the line $x=-\frac{1}{2}$, so $v_2=(x,y)\mapsto(-1-x,y)$.
Swapping and re-scaling the vectors $v_1\lra v_2$ is realised by the inversion with respect to the circle $x^2+y^2=1$ so that $v_2\mapsto v_1\mapsto \dfrac{v_2}{x^2+y^2}$, see 
Fig.~\ref{fig:signs}~(right).
\bt
\end{lem}


\begin{table}[h]
\caption{The sign of a lattice $\La\subset\R^2$ can be found from an obtuse superbase  with $v_1=(1,0)$, $v_2=(x,y)$, see Lemma~\ref{lem:signs}(a), Fig.~\ref{fig:signs}.
If any inequality becomes equality, then $\sign(\La)=0$.}
\medskip

\begin{tabular}{ccll}    
$k$ & $\sign(\La)$ & conditions on $v_2=(x,y)$ in the $k$-th subregion in Fig.~\ref{fig:signs}  & $p_{ij}$ inequalities \\
\hline

1 & $+$  & $-\frac{1}{2}<x<0$, \hspace*{4mm} $x^2+y^2>1$ & $p_{12}<p_{01}<p_{02}$ \\

2 & $-$  & $-1<x<-\frac{1}{2}$, \hspace*{27mm} 
$x^2+2x+y^2>0$ & $p_{01}<p_{12}<p_{02}$ \\

3 & $+$  & $-1<x<-\frac{1}{2}$, \hspace*{2mm} 
$x^2+y^2>1$, \hspace*{7mm}
$x^2+2x+y^2<0$ & $p_{01}<p_{02}<p_{12}$ \\

4 & $-$  & $-1<x<-\frac{1}{2}$, \hspace*{2mm} $x^2+y^2<1$, \hspace*{8mm}$x^2+x+y^2>0$ & $p_{02}<p_{01}<p_{12}$ \\

5 & $+$  & $-\frac{1}{2}<x<0$, \hspace*{4mm} 
$x^2+x+y^2>0$, \hspace*{2mm} 
$x^2+2x+y^2<0$ & $p_{02}<p_{12}<p_{01}$ \\

6 & $-$  & $-\frac{1}{2}<x<0$, \hspace*{4mm} $x^2+y^2<1$, \hspace*{8mm}$x^2+2x+y^2>0$ & $p_{12}<p_{02}<p_{01}$
\end{tabular}
\label{tab:signs}
\end{table}

\begin{proof}
\textbf{(a)}
Let $B=\{v_0,v_1,v_2\}$ be an obtuse superbase of $\La$.
Any point $p\in\La$ can be translated to the origin.
Then a suitable rotation puts the basis vector $v_1$ along the positive $x$-axis so that $v_1=(s,0)$ for $s>0$.
The uniform scaling by the factor $s$, maps $v_1$ to $(1,0)$.
Since both vectors $v_0,v_2$ have non-acute angles with $v_1$, they should have non-positive $x$-coordinates.
Since the vectors $v_0,v_2$ have a non-acute angle, one of them should be in the second quadrant $\{x\leq 0<y\}$.
Since we can swap $v_0,v_2$ without affecting $\La$, we can assume that $v_2=(x,y)$ for $x\leq 0<y$.
Then $v_0=(-1-x,-y)$.
The ordered superbase 
$B=\{v_0,v_1,v_2\}$ has the conorms 
$$p_{12}=-v_1\cdot v_2=-x\geq 0,\quad
p_{01}=-v_0\cdot v_1=1+x,\quad
p_{02}=-v_0\cdot v_2=x^2+x+y^2.$$
Since all conorms should be non-negative, we need that $0\leq p_{01}=1+x$, $x\geq -1$.
Also $0\leq p_{02}=x^2+x+y^2=(x+\frac{1}{2})^2-\frac{1}{4}+y^2$, so $(x+\frac{1}{2})^2+y^2\geq\frac{1}{4}$.
The endpoint $(x,y)$ of $v_2$ should be in the vertical strip $\{-1\leq x\leq 0\}$ non-strictly above the green circle $C(-\frac{1}{2};\frac{1}{2})$ with the centre $(-\frac{1}{2},0)$ and radius $\frac{1}{2}$ in Fig.~\ref{fig:signs}.
The yellow region $\Obt$ of allowed endpoints $(x,y)$ of $v_2$ in Fig.~\ref{fig:signs} is bounded by the vertical lines $x=0$, $x=-1$ and the green circle $C(-\frac{1}{2};\frac{1}{2})$. 
All boundary points represent all rectangular lattices.
For example, the points $(x,y)=(0,1)$ and $(x,y)=(-1,1)$ in the vertical boundaries represent the same square lattice.
For $(x,y)=(-\frac{1}{2},\frac{1}{2})$ in the green circle $C(-\frac{1}{2};\frac{1}{2})$, the vectors $v_0=(-\frac{1}{2},-\frac{1}{2})$ and $v_2=(-\frac{1}{2},\frac{1}{2})$ span a square unit cell with edge-length $\frac{1}{\sqrt{2}}$.
Now we split the yellow region into three pairs of symmetric subregions according to inequalities between three conorms.
\medskip

The inequality $p_{12}<p_{01}$ is equivalent to $-x<1+x$, $x>-\frac{1}{2}$, see Table~\ref{tab:sign_regions}.
The inequality $p_{01}<p_{02}$ is equivalent to $1+x<x^2+x+y^2$, $x^2+y^2>1$, so the point $(x,y)$ is above the circle $C(0;1)$ with the centre $(0,0)$ and radius $1$ in Fig.~\ref{fig:signs}.
The inequality $p_{12}<p_{02}$ is equivalent to $-x<x^2+x+y^2$, $(x+1)^2+y^2>1$, so the point $(x,y)$ is above the circle $C(-1;1)$ with the centre $(-1,0)$ and radius $1$.
\medskip

\begin{table}
\caption{Inequalities between conorms are interpreted in terms of endpoints $(x,y)$ of a vector $v_2$ complementing $v_1=(1,0)$ in an obtuse superbase $\{v_0,v_1,v_2\}$, see Lemma~\ref{lem:signs}.}
\medskip

\begin{tabular}{lll}      
$p_{ij}$ inequality & condition on $(x,y)$ & subregion within the yellow region $\Obt$ in Fig.~\ref{fig:signs} \\
\hline

$p_{02}\geq 0$ & $x^2+x+y^2\geq 0$ & the region $\Obt$ is non-strictly above $C(-\frac{1}{2};\frac{1}{2})$ \\

$p_{12}<p_{01}$ & $-\frac{1}{2}<x<0$ & the right hand side vertical strip of the region $\Obt$ \\

$p_{01}<p_{02}$ & $x^2+y^2>1$ & the subregion in $\Obt$ above the circle $C(0;1)$ \\

$p_{12}<p_{02}$ & $x+2x+y^2>0$ & the subregion in $\Obt$ above the circle $C(-1;1)$
\end{tabular}
\label{tab:sign_regions}
\end{table}

The inequalities on $p_{ij}$ from Table~\ref{tab:sign_regions} justify that the region $\Obt$ splits into six subregions split by the vertical line $x=-\frac{1}{2}$ and two circles $C(-1;0)$ and $C(0;1)$.
Each subregion is defined by one of six possible orderings of the conorms $p_{12},p_{01},p_{02}$, see the last column of Table~\ref{tab:signs}.
To check the signs in the second column of Table~\ref{tab:signs}, notice that if $p_{ij}$ is a minimal conorm, the formula $p_{ij}=\frac{1}{2}(v_i^2+v_j^2-v_k^2)$ from~(\ref{dfn:vonorms}b) implies that $v_i,v_j$ are the shortest of three vectors $v_0,v_1,v_2$.
\medskip

For example, Table~\ref{tab:signs} says that $v_1=(1,0)$ and $v_2$ are the two shortest vectors in the cases of the first and last rows.
In the first row, $v_2=(x,y)$ has the length $|v_2|=\sqrt{x^2+y^2}>1=|v_1|$, hence by Definition~\ref{dfn:sign} $\sign(\La)$ equals the sign of $\det(v_1,v_2)=y>0$.
In the last row, $v_2=(x,y)$ has the length $|v_2|=\sqrt{x^2+y^2}<1=|v_1|$, hence by Definition~\ref{dfn:sign} $\sign(\La)$ equals the sign of $\det(v_2,v_1)=-y<0$.
The signs in the remaining four rows of Table~\ref{tab:signs} are similarly checked.
\medskip

If any of the strict inequalities above becomes equality, we get a point either on the boundary of $\Obt$ (representing all rectangular lattices in $|R^2$) or in one of the lines $x=-\frac{1}{2}$ or the circles $C(0;1)$ and $C(-1;1)$. 
These internal curves contain points $(x,y)$ representing centred rectangular lattices.
For instance, the triple intersection of the internal curves at $(x,y)=(-\frac{1}{2},\frac{\sqrt{3}}{2})$ represents all hexagonal lattices.
All these lattices are mirror-symmetric and have $\sign(\La)=0$.
\medskip

\noindent
\textbf{(b)}
Any two vectors of an obtuse superbase $B=\{v_0,v_1,v_2\}$ can be mapped by similarity to $(1,0)$ and $(x,y)$.
Each of the resulting six pairs $(x,y)$ belongs to one of the six subregions marked by $k=1,2,3,4,5,6$ in the middle picture of Fig.~\ref{fig:signs}.
It suffices to understand the action of two transpositions $v_0\lra v_2$ and $v_1\lra v_2$.
\medskip

When we swap $v_2=(x,y)$ and $v_0=(-1-x,-y)$, while keeping $v_1=(x,y)$ fixed, we reflect the lattice $\La(B)$ generated by $B$ in the $x$-axis so that $v_1=(x,y)$ has an obtuse angle to the image $(-1-x,y)$ of $v_0$.
This new vector $(-1-x,y)$ plays the role of $v_2$ in the reflected lattice and is symmetric to $v_2=(x,y)$ in the vertical line $x=-\frac{1}{2}$ within the yellow region $\Obt$ in Fig.~\ref{fig:signs}~(right).     
\medskip

When we swap $v_1=(1,0)$ and $v_2=(x,y)$, the second vector is divided by its length $|v_2|=\sqrt{x^2+y^2}$.
Hence the first vector $v_1$ maps to the vector that is parallel to $v_2=(x,y)$ and has the length $1/\sqrt{x^2+y^2}$.
This new vector $v'_2=v_2/(x^2+y^2)$ plays the role of $v_2$ and is obtained from $v_2=(x,y)$ by the inversion with respect to the circle $x^2+y^2=1$.
The inversion keeps all points on $x^2+y^2=1$ fixed, maps the $y$-axis $x=0$ to itself, swaps the half-line $\{x=0, \; y>0\}$ with the upper half-circle $\{x^2+x+y^2=0,\; y>0\}$.
Compositions of the symmetry in $x=-\frac{1}{2}$ and this inversion generate up to six images of $(x,y)$ in the six subregions of $\Obt$, though the point $(x,y)=(-\frac{1}{2},\frac{\sqrt{3}}{2})$ representing all hexagonal lattices is fixed.
\ws
\end{proof}

\begin{prop}[reduced bases]
\label{prop:reduced_bases}
\textbf{(a)}
Up to isometry in $\R^2$, all reduced bases $v_1,v_2$ from Definition~\ref{dfn:reduced_cell} are in a 1-1 correspondence with all obtuse superbases $B=\{v_0,v_1,v_2\}$ such that $|v_1|\leq|v_2|\leq|v_0|$.
Up to isometry, any lattice $\La\subset\R^2$ has a unique reduced basis specified by the conditions of Definition~\ref{dfn:reduced_cell}.
\medskip

\noindent
\textbf{(b)}
Up to rigid motion, any lattice has a unique reduced basis in Definition~\ref{dfn:reduced_cell}.
\bt
\end{prop}

\begin{figure}[h]
\includegraphics[width=\textwidth]{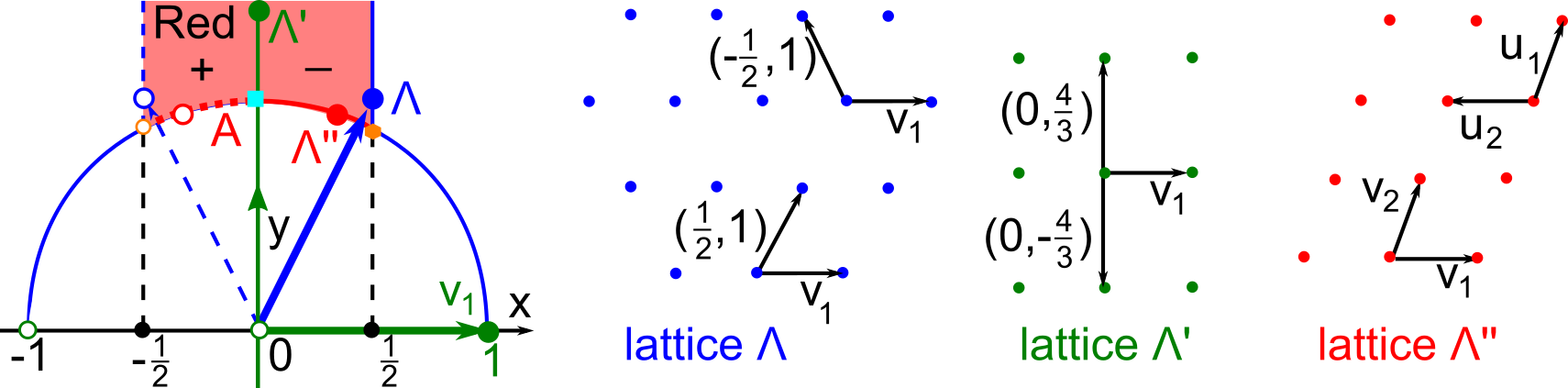}
\caption{\textbf{Left}: any reduced basis in Definition~\ref{dfn:reduced_cell} up to orientation-preserving similarity maps to $v_1=(1,0)$ and $v_2=(x,y)\in\RB$ from Proposition~\ref{prop:reduced_bases}.
\textbf{Right}: for each of the lattices $\La,\La',\La''$ represented by small blue, green, red  circles/disks on the right, the conditions of Definition~\ref{dfn:reduced_cell} choose one reduced basis among two bases that differ up to rigid motion.
}
\label{fig:reduced_bases}
\end{figure}

\begin{proof}
\textbf{(a)}
Up to similarity by Lemma~\ref{lem:signs}(b), any lattice $\La$ has an obtuse superbase $B=\{v_0,v_1,v_2\}$ with $v_1=(1,0)$ and $v_2=(x,y)$, where the point $(x,y)$ belongs to the yellow region $\Obt$ in Fig.~\ref{fig:signs}.
By Lemma~\ref{lem:signs}(b) the six permutations of $v_0,v_1,v_2$ are realised by internal symmetries of $\Obt$, so we may assume that $|v_1|\leq|v_2|\leq|v_0|$.
The equivalent inequalities in conorms $p_{12}\leq p_{01}\leq p_{02}$ define the 1st subregion of $\Obt$ in Fig.~\ref{fig:signs}(middle), which coincides with the closure of the region $\overline{\RB_+}=\{(x,y)\in\R^2\mid x^2+y^2\geq 1, -\frac{1}{2}\leq x\leq 0<y\}$
in Fig.~\ref{fig:reduced_bases}~(left).
\medskip

Due to uniqueness of $B$ up to isometry by Theorem~\ref{thm:isometric_superbases}, the position of $v_2=(x,y)\in\overline{\RB_+}$ is unique for $\La$.
Up to uniform scaling and reflection $y\lra -y$, the closure $\overline{\RB_+}$ is defined by the same conditions $v_1^2\leq v_2^2=x^2+y^2$ and $-\dfrac{1}{2}\leq\dfrac{v_1\cdot v_2}{v_1^2}=x\leq 0$ as a reduced basis whose uniqueness up to isometry follows now.
\medskip

\noindent
\textbf{(b)}
If orientation should be preserved, Theorem~\ref{thm:isometric_superbases} proves the uniqueness of the obtuse superbase $B=\{v_0,v_1,v_2\}$ from part (a) up to rigid motion for any non-rectangular lattice $\La$. 
Cyclic permutations of $v_0,v_1,v_2$  allow us to assume that $v_1$ is the shortest vector.
The equivalent condition on conorms says that $p_{02}$ is the largest, hence $v_2=(x,y)$ belongs to the first two subregions of $\Obt$ in Fig.~\ref{fig:signs}.
\medskip

For the first open subregion with $\sign(B)=+1$, the conditions $v_1^2<v_2^2=x^2+y^2$, $-\frac{1}{2}<x=\dfrac{v_1\cdot v_2}{v_1^2}<0<y=\det(v_1,v_2)$ of Definition~\ref{dfn:reduced_cell}  specify the interior $\overline{\RB_+}$, so the basis $\{v_1,v_2\}$ is reduced.
For the second open subregion with $\sign(B)=-1$, another basis $\{v_1,v_1+v_2\}$ is reduced by Definition~\ref{dfn:reduced_cell} because the shifted point $v_1+v_2=(x+1,y)$ belongs to the interior of the right-half region in Fig.~\ref{fig:reduced_bases}~(left):
$\overline{\RB_-}=\{(x,y)\in\R^2\mid x^2+y^2\geq 1,\, 0<x<\frac{1}{2},\, 0<y\}$.
\medskip

It remains to consider singular cases.
We include the common boundary line $x=\dfrac{v_1\cdot v_2}{v_1^2}=-\dfrac{1}{2}$ and the boundary round arcs of both subregions represent mirror-symmetric lattices with a unique (up to rigid motion) obtuse superbase.  
We exclude these boundaries from the first region, include them into the second region and shift by $x\mapsto x+1$ so that the unique reduced basis $v_1,v_1+v_2$ satisfies the conditions $\dfrac{v_1\cdot v_2}{v_1^2}=\dfrac{1}{2}$ and $v_1\cdot v_2\geq 0$ for $|v_1|=|v_2|$ in Definition~\ref{dfn:reduced_cell}.
\medskip

The above boundaries in Fig.~\ref{fig:reduced_bases}~(left) include the blue and red points representing the basis vectors $v_2=(x,y)$ of the lattices $\La,\La''$, respectively.
\medskip

The final boundary lines $x=-1$ and $x=0$ represent rectangular lattices with a unit cell $a\times b$ for $0<a<b$ and two obtuse superbases related by reflection, not by rigid motion, for example $v_1=(a,0)$, $v_2=(0,\pm b)$.
In this case there is no 1-1 correspondence between obtuse superbases and reduced bases.
Definition~\ref{dfn:reduced_cell} selects a unique reduced basis $v_1=(a,0)$, $v_2=(0,b)$ due to $\det(v_1,v_2)>0$.
\ws
\end{proof}

The dotted arc $A$ in Fig.~\ref{fig:reduced_bases}~(left) should be indeed excluded from \cite[Fig.~1.3 on p.~85]{engel2004lattice}, otherwise the lattice $\La''$ in the last picture of Fig.~\ref{fig:reduced_bases} has two potential reduced bases $v_1=(1,0)$, $v_2=(\frac{1}{3},\frac{2\sqrt{2}}{3})$ and $u_1=(\frac{1}{3},\frac{2\sqrt{2}}{3})$ and $u_2=(-1,0)$. 
Indeed, all basis vectors have length 1 and the second basis can be rotated to $u'_1=(1,0)$, $u'_2=(-\frac{1}{3},\frac{2\sqrt{2}}{3})\in A$.
The second basis $(u_1,u_2)$ of the same lattice $\La''$ is related to $(v_1,v_2)$ by a reflection, but not by rigid motion.
So the region $\RB$ with the excluded left boundary for $x<0$ contains a unique vector $v_2=(x,y)$ of a reduced basis up to orientation-preserving similarity.
Forgetting about the uniform scaling, we get uniqueness of a reduced basis up to rigid motion.
\medskip

The region $\RB$ in Fig.~\ref{fig:reduced_bases}~(left) is a fundamental domain of all bases by the action of $\SO(\R^2)\times\R_+$ and $\GL_2(\Z)$ in the sense that any lattice up to orientation-preserving similarity can be represented by a unique point $(x,y)\in\RB$.
$\RB$ or any other half-open fundamental domain of a group action suffers from \emph{discontinuity on boundary} when close lattices are represented by distant bases. 
For each of the lattices $\La,\La''$ in Fig.~\ref{fig:reduced_bases}, a slight perturbation of the non-reduced basis makes it reduced but distant from the initial reduced basis up to rigid motion.
The discontinuity above can be resolved by identifying boundary points of  of $\RB$ by the reflection $x\lra-x$.
Section~\ref{sec:classification2d} will describe a simpler way to continuously parameterise lattices up to orientation-preserving similarity in Corollary~\ref{cor:similar_lattices2d}.
 
\section{Complete classifications of 2D lattices up to isometry and similarity}
\label{sec:classification2d}

Lemma~\ref{lem:lattice_invariants} showed that $\RI(\La),\RI^o(\La)$ are invariants of lattices up to isometry and rigid motion, respectively. 
To prove completeness of the invariants in Theorem~\ref{thm:classification2d}, Lemma~\ref{lem:superbase_reconstruction} reconstructs an obtuse superbase of $\La$.
Corollary~\ref{cor:similar_lattices2d} will classify lattices up to similarity by projected invariants introduced in Definition~\ref{dfn:PF}.

\begin{lem}[superbase reconstruction]
\label{lem:superbase_reconstruction}
An obtuse superbase $B=\{v_0,v_1,v_2\}$ of a lattice $\La\subset\R^2$ can be uniquely reconstructed up to isometry and up to rigid motion from its root invariant $\RI(\La)$ and its oriented root invariant $\RI^o(\La)$, respectively.
If $\RI(\La)=(r_{12},r_{01},r_{02})$, the basis vectors $v_1,v_2$ are determined by 
$$|v_1|=\sqrt{r_{12}^2+r_{01}^2},\quad 
|v_2|=\sqrt{r_{12}^2+r_{02}^2},\quad
\cos\angle(v_1,v_2)=\dfrac{-r_{12}^2}{\sqrt{r_{12}^2+r_{01}^2}\sqrt{r_{12}^2+r_{02}^2}},
$$ and span a primitive unit cell of the area $A(\La)=\sqrt{r_{12}^2 r_{01}^2 + r_{12}^2 r_{02}^2 + r_{01}^2 r_{02}^2}$.
\bt
\end{lem}
\begin{proof}
Assuming that a root invariant $\RI(\La)$ is ordered as $r_{12}\leq r_{01}\leq r_{02}$, we will build an obtuse superbase $\{v_0,v_1,v_2\}$ such that $r_{ij}=\sqrt{-v_i\cdot v_j}$ for any distinct $i,j\in\{0,1,2\}$.
Find the lengths from (\ref{dfn:vonorms}a): $|v_i|=\sqrt{p_{12}+p_{0i}}=\sqrt{r_{12}^2+r_{0i}^2}$ for $i=1,2$.
Using $v_1\cdot v_2=-r_{12}^2$, the anticlockwise angle has 
$\cos\angle(v_1,v_2)=\dfrac{v_1\cdot v_2}{|v_1|\cdot|v_2|}=\dfrac{-r_{12}^2}{\sqrt{r_{12}^2+r_{01}^2}\sqrt{r_{12}^2+r_{02}^2}}$.
The unit cell $U(v_1,v_2)$ has the area 
$$A(\La)=|v_1|\cdot |v_2|\sin\al=|v_1|\cdot |v_2|\sqrt{1-\cos^2\al}=\sqrt{|v_1|^2|v_2|^2-(v_1\cdot v_2)^2}=$$
$$=\sqrt{(r_{12}^2+r_{01}^2)(r_{12}^2+r_{02}^2)-r_{12}^2}
=\sqrt{r_{12}^2 r_{01}^2 + r_{12}^2 r_{02}^2 + r_{01}^2 r_{02}^2}.$$
Up to rigid motion, the length $|v_1|$ is enough to fix the vector $v_1$ along the positive $x$-axis.
The length $|v_2|$ and $\cos\angle(v_1,v_2)$ determine the position of $v_2$ relative to the fixed vector $v_1$ up to reflection in the $x$-axis.
Up to isometry or if $\sign(\La)=0$ (when $\La$ is mirror-symmetric), the above options for $v_2$ are not important.
If $\sign(\La)=+1$, then we choose $v_2$ in the upper half-plane above the $x$-axis so that $\angle(v_1,v_2)\in(90^\circ,180^\circ)$, otherwise we put $v_2$ into the lower half-plane.
\medskip

Finally, $v_0=-v_1-v_2$ and the reconstructed ordered obtuse superbase $B=\{v_0,v_1,v_2\}$ is unique up to isometry and up to rigid motion by Theorem~\ref{thm:isometric_superbases}.
\ws
\end{proof}

\begin{thm}[isometric 2D lattices $\lra$ root invariants]
\label{thm:classification2d}
Any lattices $\La,\La'\subset\R^2$ are isometric if and only if their root invariants coincide: $\RI(\La)=\RI(\La')$.
Any lattices $\La,\La'$ are related by rigid motion
if and only if $\RI^o(\La)=\RI^o(\La')$.
\bt 
\end{thm}
\begin{proof}
The part \emph{only if} ($\Rightarrow$) follows from Lemma~\ref{lem:lattice_invariants}(a) saying that $\RI(\La),\RI^o(\La)$ are invariant under isometry and rigid motion, respectively.
The part \emph{if} ($\Leftarrow$) follows from Lemma~\ref{lem:superbase_reconstruction} reconstructing a superbase from $\RI(\La)$ or $\RI^o(\La)$. 
\ws
\end{proof}

The above classification helps prove that some other isometry invariants of lattices are also complete and continuous.
By (\ref{dfn:vonorms}ab) the voform $\VF=(v_0^2,v_1^2,v_2^2)$ and coform $\CF=(p_{12},p_{01},p_{02})$ are both complete if considered
up to $3!$ permutations.
The root invariant $\RI$ is a uniquely ordered version of $\CF$ and deserves its own name.
The square roots $r_{ij}=\sqrt{p_{ij}}$ have original units of vector coordinates.
\medskip

The oriented part of Theorem~\ref{thm:classification2d} didn't appear in the past to the best of our knowledge.
Conway and Sloane studied 2D lattices in \cite[section~6]{conway1992low} only up to general isometry including reflections.
Here is the closest formal claim from \cite{conway1992low}.

\begin{lem}[geometry of vonorms {\cite[Theorem~7]{conway1992low}}]
\label{lem:VF}
For any lattice $\La\subset\R^2$, the vonorms $v_0^2,v_1^2,v_2^2$ are squared lengths of three shortest Voronoi vectors. 
\bt
\end{lem}

\begin{figure}[h]
\includegraphics[height=34mm]{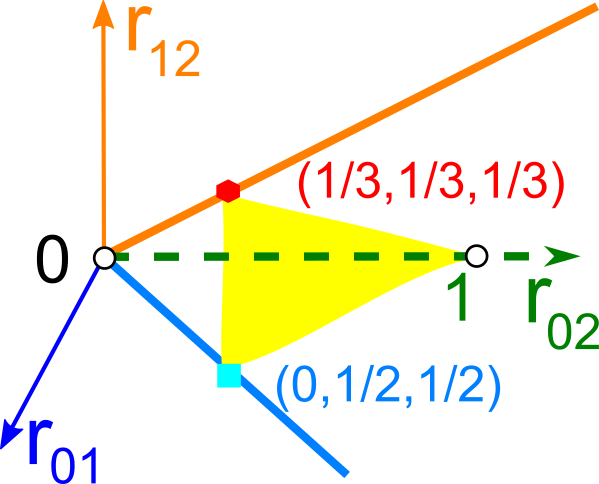}\hspace*{2mm}
\includegraphics[height=34mm]{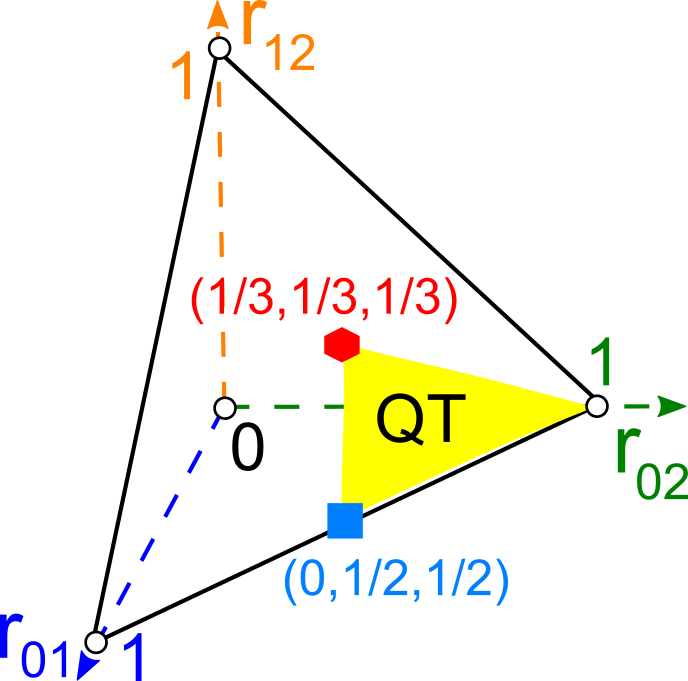}\hspace*{2mm}
\includegraphics[height=34mm]{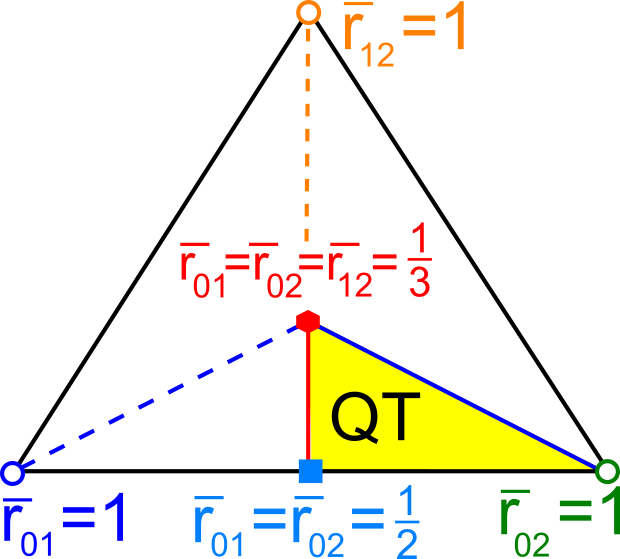}
\caption{\textbf{Left}: the triangular cone $\TC=\{(r_{12},r_{01},r_{02})\in\R^3 \mid 0\leq r_{12}\leq r_{01}\leq r_{02}\neq 0\}$ represents the space $\RIS$ of all root invariants of 2D lattices, see Definition~\ref{dfn:TC}. 
\textbf{Middle}: 
$\TC$ projects to the quotient triangle $\QT=\TC\cap\{r_{12}+r_{01}+r_{02}=1\}$ representing the space $\LSS$ of 2D lattices up to similarity, see Corollary~\ref{cor:similar_lattices2d}.
\textbf{Right}: the quotient triangle $\QT$ can be parameterised by $x=\bar r_{02}-\bar r_{01}\in[0,1)$ and $y=3\bar r_{12}\in[0,1]$, see $\QT$ also in Fig.~\ref{fig:QT+QS}.
}
\label{fig:TC}
\end{figure}

\begin{figure}[h]
\includegraphics[height=55mm]{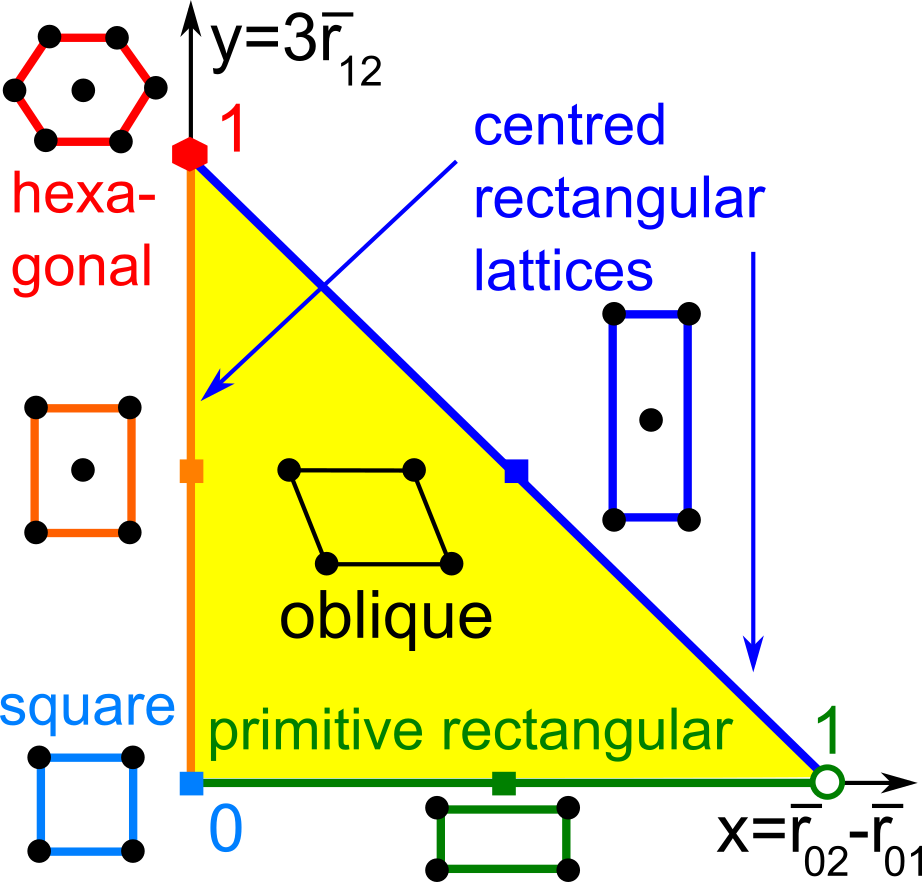}
\hspace*{3mm}
\includegraphics[height=55mm]{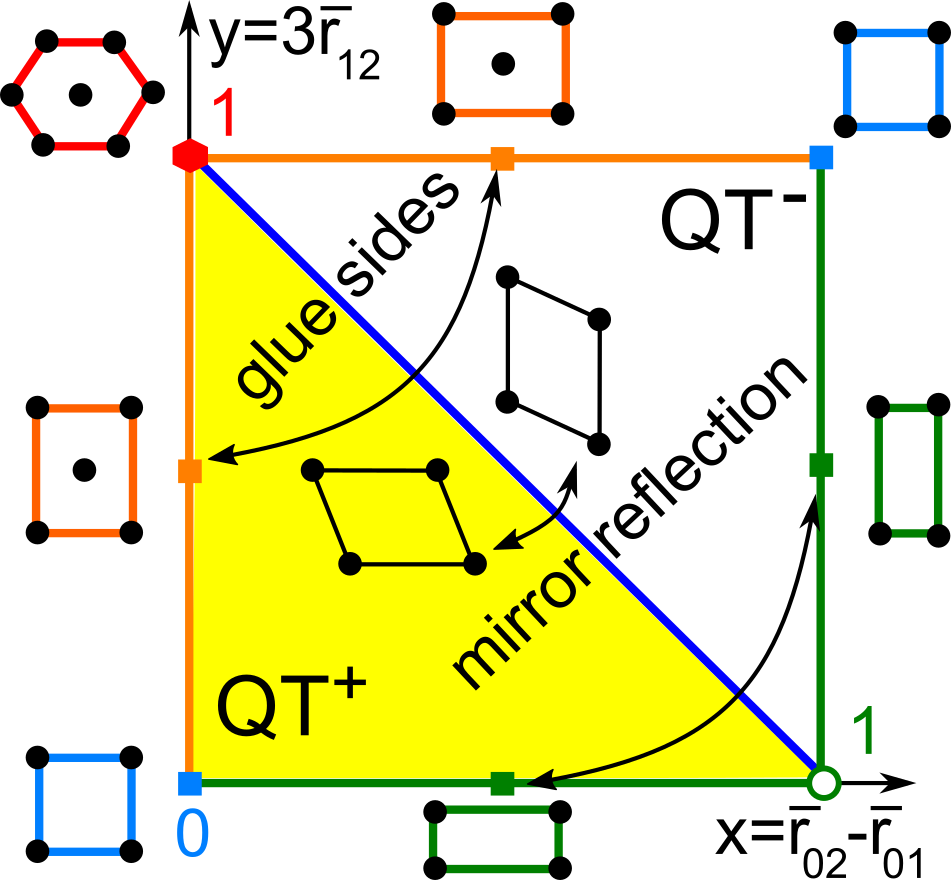}
\caption{\textbf{Left}:
all projected invariants $\PI(\La)$ of lattices $\La\subset\R^2$ live in the quotient triangle $\QT$ from Fig.~\ref{fig:TC}, which is parameterised by $x=\bar r_{02}-\bar r_{01}\in[0,1)$ and $y=3\bar r_{12}\in[0,1]$.
\textbf{Right}: mirror reflections $\La^{\pm}$ of any non-mirror-symmetric lattice can be represented by a pair of points in the quotient square $\QS=\QT^+\cup\QT^-$ symmetric in the diagonal $x+y=1$.}
\label{fig:QT+QS}
\end{figure}

Theorem~\ref{thm:classification2d} and Lemma~\ref{lem:VF} imply that, after taking square roots of vonorms, the ordered lengths, say $|v_1|\leq|v_2|\leq|v_0|$, form a complete invariant that should satisfy the triangle inequality $|v_1|+|v_2|\geq|v_0|$.
This inequality is the only disadvantage of the complete invariant $|v_1|\leq|v_2|\leq|v_0|$ in comparison with ordered root products $r_{12}\leq r_{01}\leq r_{02}$, which are easier to visualise in Fig.~\ref{fig:TC},~\ref{fig:QT+QS}.
\medskip

Classification Theorem~\ref{thm:classification2d} says that all isometry classes of lattices $\La\subset\R^2$ are in a 1-1 correspondence with all ordered triples $0\leq r_{12}\leq r_{01}\leq r_{02}$ of root products in $\RI(\La)$.
Only the smallest root product $r_{12}$ can be zero, two others $r_{01}\leq r_{02}$ should be positive, otherwise $v_1^2=r_{12}^2+r_{01}^2=0$ by formulae (\ref{dfn:vonorms}a).
\medskip

We explicitly describe the set of all possible root invariants, which will be later converted into metric spaces with continuous metrics in Definitions~\ref{dfn:RM} and~\ref{dfn:RMo}.

\begin{dfn}[triangular cone $\TC$]
\label{dfn:TC}
All root invariants $\RI(\La)=(r_{12},r_{01},r_{02})$ of lattices $\La\subset\R^2$ live in the \emph{triangular cone} $\TC=\{0\leq r_{12}\leq r_{01}\leq r_{02}\}$ within the octant $\Oct=[0,+\infty)^3$ excluding the axes in the coordinates $r_{12},r_{01},r_{02}$, see Fig.~\ref{fig:TC}~(left).
The boundary $\bd(\TC)$ of the cone $\TC$ consists of root invariants of all mirror-symmetric lattices from Lemma~\ref{lem:achiral_lattices}: 
the bisector planes $\{r_{01}=r_{02}\}$ and $\{r_{12}=r_{01}\}$ within $\TC$. 
The orange line $\{0<r_{12}=r_{01}=r_{02}\}\subset\bd(\TC)$ in Fig.~\ref{fig:TC}~(left) consists of root invariants of hexagonal lattices with a minimum inter-point distance $r_{12}\sqrt{2}$.
The blue line $\{r_{12}=0<r_{01}=r_{02}\}\subset\bd(\TC)$ consists of root invariants of square lattices with a minimum inter-point distance $r_{01}$.
\bs
\end{dfn}

To classify lattices up to similarity, it is convenient to scale them by the \emph{size} $\si(\La)=r_{12}+r_{01}+r_{02}$.
This sum is a simpler uniform measure of size than (say) the unit cell area $A(\La)$ from Lemma~\ref{lem:superbase_reconstruction}, which can be small even for long cells. 

\begin{dfn}[projected invariants $\PI(\La)$ and $\PI^o(\La)$]
\label{dfn:PF}
The \emph{triangular projection} $\TP:\TC\to\{r_{12}+r_{01}+r_{02}=1\}$ divides each coordinate by the \emph{size} $\si(\La)=r_{12}+r_{01}+r_{02}$ and gives $\PRF(\La)=(\bar r_{12},\bar r_{01},\bar r_{02})=\dfrac{(r_{12},r_{01},r_{02})}{r_{12}+r_{01}+r_{02}}$ in $\TC\cap\{r_{12}+r_{01}+r_{02}=1\}$.
Then we map $(\bar r_{12},\bar r_{01},\bar r_{02})$ to the \emph{projected invariant} $\PI(\La)=(x,y)$ with $x=\bar r_{02}-\bar r_{01}\in[0,1)$ and $y=3\bar r_{12}\in[0,1]$ in the \emph{quotient triangle} $\QT=\{(x,y)\in\R^2\mid 0\leq x<1,\, 0\leq y\leq 1,\, x+y\leq 1\}$, see Fig.~\ref{fig:QT+QS}.
\medskip

All oriented root invariants $\RI^o(\La)$ live in the \emph{doubled cone} $\DC$ that is the union of two triangular cones $\TC^{\pm}$, where we identify any two boundary points representing the same root invariant $\RI(\La)$ with $\sign(\La)=0$.
The \emph{oriented projected invariant} $\PI^o(\La)=(x,y)^{\pm}$ is $\PI(\La)$ with the superscript from $\sign(\La)$.
\bs
\end{dfn}

The inequality $1\geq x+y=(\bar r_{02}-\bar r_{01})+3\bar r_{12}$ follows after multiplying both sides by the size $\si(\La)$, because $r_{12}+r_{01}+r_{02}\geq (r_{02}-r_{01})+3r_{12}$ becomes $r_{01}\geq r_{12}$.
\medskip

The set of oriented projected invariants $\PI^o$ is visualised in Fig.~\ref{fig:QT+QS} (right) as the \emph{quotient square} $\QS$ obtained by gluing the quotient triangle $\QT^+$ with its mirror image $\QT^-$.
The boundaries of both triangles excluding the vertex $(x,y)=(1,0)$ are glued by the diagonal reflection $(x,y)\lra(1-y,1-x)$.
Any pair of points $(x,y)\in\QT^+$ and $(1-y,1-x)\in\QT^-$ in Fig.~\ref{fig:QT+QS}~(right) represent mirror images of a lattice up to similarity, see Corollary~\ref{cor:similar_lattices2d}.
So $\QS$ is a topological sphere without a single point and will be parameterised by geographic-style coordinates in \cite{bright2021geographic}.
\medskip

Following Fig.~\ref{fig:forms2d}, any square lattice 
has a root invariant $\RI=(0,a,a)$, so its projected invariant $\PI=(0,0)$ is at the bottom left vertex of $\QT$ in Fig.~\ref{fig:QT+QS}~(left), identified with top right vertex of $\QS$ in Fig.~\ref{fig:QT+QS}~(right).
Any hexagonal lattice has a root invariant $\RI=(a,a,a)$, so its projected invariant $\PI=(0,1)$ is at the top left vertex of $\QT$ in Fig.~\ref{fig:QT+QS}~(left), identified with bottom right vertex of $\QS$.
\medskip

By Example~\ref{exa:achiral_lattices}(a) any rectangular lattice has $\RI=(0,a,b)$ for $a<b$, hence its projected invariant $\PI=(\frac{b-a}{a+b},0)$ belongs to the bottom edge of $\QT$ identified with the top edge of $\QS$.
By Example~\ref{exa:achiral_lattices}(b) any lattice with a mirror-symmetric Voronoi domain has $\RI$ with 0 or two equal root products.
Such lattices have a rhombic unit cell and form the centred rectangular Bravais class.
Their projected invariants belong to the vertical edges and diagonal of $\QS$ in Fig.~\ref{fig:QT+QS}~(right).
The companion paper \cite{bright2021geographic} discusses Bravais classes of 2-dimensional lattices in detail.
\medskip

In the theory of complex functions, any lattice $\La\subset\R^2$ can be considered as a subgroup of the complex plane $\C$ whose quotient $\C/\La$ is a torus.
By the Riemann mapping theorem any compact Riemann surface of genus 1 is conformally equivalent (holomorphically homeomorphic) to the quotient $\C/\La$ for some lattice $\La$, see \cite[Section 5.3]{jost2013compact}.
Such tori $\C/\La$ and $\C/\La'$ are \emph{conformally equivalent} if and only if $\La,\La'$ are similar, see \cite[Theorem~6.1.4]{jones1987complex}.
The spaces $\LSS(\R^2)$ and $\LSS^o(\R^2)$ of all lattices $\La\subset\C=\R^2$ up to similarity and orientation-preserving similarity are the quotient triangle $\QT$ and square $\QS$, respectively, see Fig.~\ref{fig:QT+QS}.

\begin{cor}[similar lattices $\lra$ projected invariants $\PI$]
\label{cor:similar_lattices2d}
Lattices $\La,\La'\subset\R^2$ are \emph{similar} (related by an isometry composed with a uniform scaling) 
if and only if their projected invariants are equal: $\PI(\La)=\PI(\La')$.
The lattices $\La,\La'$ are related by an orientation-preserving similarity if and only if $\PI^o(\La)=\PI^o(\La')$.
\bt
\end{cor}
\begin{proof}
follows from Theorem~\ref{thm:classification2d} because a uniform scaling of all basis vectors $v_i\mapsto sv_i$ by a factor $s>0$ multiplies all root products $r_{ij}=\sqrt{-v_i\cdot v_j}$ by $s$, which is neutralised by the triangular projection $\TP$ from Definition~\ref{dfn:TC}.
\ws
\end{proof}

\begin{lem}[criteria of mirror-symmetric lattices in $\R^2$]
\label{lem:sign0}
A lattice $\La$ in $\R^2$ is mirror-symmetric if and only if one of the following equivalent conditions holds:
$\sign(\La)=0$ or $\RI(\La)\in\bd\TC$ or $\PI(\La)\in\bd\QT$.
So the boundaries of the triangular cone $\TC$ and the quotient triangle $\QT$ consist of root invariants and projected invariants, respectively, of all mirror-symmetric lattices $\La\subset\R^2$.
\bt
\end{lem}
\begin{proof}
By Lemma~\ref{lem:lattice_invariants} a lattice $\La$ is mirror-symmetric if and only if $\sign(\La)=0$.
By Lemma~\ref{lem:achiral_lattices} any mirror-symmetric lattice $\La$ has $\RI(\La)$ with 0 (rectangular lattice) or with two equal root products (centred rectangular lattices).
The last conditions on $\RI$ define the boundary $\bd\TC$ of the triangular cone in Fig.~\ref{fig:TC} or, equivalently, the projected invariant $\PI(\La)$ belongs to the boundary of $\QT$ in Fig.~\ref{fig:QT+QS}~(left). 
\ws
\end{proof}

\begin{rem}[lattices via group actions]
\label{rem:group_action}
Another parameterisation of the Lattice Similarity Space $\LSS(\R^2)$ can be obtained from a fundamental domain of the action of $\GL_2(\Z)\times\R_+^\times$ on the cone $\mathcal{C}_+(\mathcal{Q}_2)$ of positive quadratic forms.
Recall that any lattice $\La\subset\R^2$ with a basis $v_1,v_2$ defines the positive quadratic form
$$Q_\La(x,y)=(xv_1+yv_2)^2=v_1^2x^2+2v_1v_2xy+v_2^2y^2=q_{11}x^2+2q_{12}xy+q_{22}y^2\geq 0$$ 
whose positivity for all $(x,y)\in\R^2-0$ means that $q_{12}^2<q_{11}q_{22}$.
The cone $\mathcal{C}_+(\mathcal{Q}_2)$ of all positive quadratic forms projects to the unit disk $\xi^2+\eta^2<1$ parameterised by $\xi=\dfrac{q_{22}-q_{11}}{q_{11}+q_{22}}$ and $\eta=\dfrac{-2q_{12}}{q_{11}+q_{22}}$.
Indeed, the positivity 
condition $q_{12}^2<q_{11}q_{22}$ for the form $Q_\La(x,y)$ is equivalent to $\xi^2+\eta^2<1$ in the coordinates above.
\medskip

A form $Q_\La$ is called \emph{reduced} if $0\leq -2q_{12}\leq q_{11}\leq q_{22}$ and $q_{11}>0$, see \cite[formula~(1.130) on p.~75]{engel2004lattice}.
The above conditions define the fundamental domain $T=\{0\leq\xi<1,\; 0\leq\eta\leq\frac{1}{2},\; \xi+2\eta\leq 1\}$, see  \cite[Fig.~8.1]{zhilinskii2016introduction}. 
This non-isosceles triangle is one of infinitely many triangular domains within the disk $\xi^2+\eta^2<1$ in \cite[Fig.~1.2 on p.~82]{engel2004lattice} or \cite[Fig.~6.2]{zhilinskii2016introduction}.
Choosing one triangular domain is equivalent to choosing a reduced basis up to isometry, not up to rigid motion.
For example, the mirror-symmetric bases of $v_1=(1,0)$, $v_2^{\pm}=(\pm\frac{1}{2},1)$ have the quadratic forms $x^2\pm xy+\frac{5}{4}y^2$ represented by $(\xi,\eta)=(\frac{1}{9},\mp\frac{4}{9})$.
\medskip

The above ambiguity or discontinuity up to rigid motion is now resolved by $\sign(\La)$ in the twice larger space $\LSS^o(\R^2)$ visualised as the quotient square $\QS$, see
Example~\ref{exa:RF_deformation}.
The coefficients of $Q_\La$ can be used for parameterising $\LIS(\R^2)$: 
$$q_{11}=v_1^2=r_{12}^2+r_{01}^2,\quad
q_{22}=v_2^2=r_{12}^2+r_{02}^2,\quad
q_{12}=v_1\cdot v_2=-r_{12}^2.$$
The triple $(v_1^2,v_1\cdot v_2,v_2^2)$ is called the \emph{metric tensor} whose disadvantages are squared units and non-homogeneity (one scalar product, two squared lengths).
\medskip

Another complete invariant is an ordered Voronoi form $v_1^2\leq v_2^2\leq v_0^2$ or the lengths $|v_1|\leq|v_2|\leq|v_0|$ of three shortest Voronoi vectors from Lemma~\ref{lem:VF}.
However, this invariant doesn't extend even to dimension $n=3$ due to a 6-parameter family of pairs of non-isometric lattices $\La_1\not\cong\La_2$ that have the same lengths of seven shortest Voronoi vectors in $\R^3$, see \cite{kurlin2022complete}.
The above reasons justify the choice of homogeneous coordinates $r_{ij}$, which easily extend to higher dimensions. 
\bs
\end{rem}

The projected invariant $\PI=(x,y)$ obtained from $\RI$ is preferable to the coordinates $(\xi,\eta)$, which define a non-isosceles triangle, while the isosceles quotient triangle $\QT$ will lead to easier formulae for metrics in the next section.
Since the metric tensor $(v_1\cdot v_2,v_1^2,v_2^2)=(-q_{12},q_{11},q_{22})$ and its 3-dimensional analogue are more familiar to crystallographers, we will rephrase key results from sections~\ref{sec:metrics}-\ref{sec:chiral_distances} by using these non-homogeneous cooordinates in the companion paper \cite{bright2021geographic}.

\begin{prop}[inverse design of 2D lattices]
\label{prop:inverse_design}
For $\si>0$ and any point $(x,y)$ in the quotient triangle $\QT$, there is a unique (up to isometry) lattice $\La$ with the projected invariant $\PI(\La)=(x,y)$ and size $\si=r_{12}+r_{01}+r_{02}$. 
Then
$$\RI(\La)=(r_{12},r_{01},r_{02})=\left(\frac{\si}{3} y,\; \frac{\si}{6}(3-3x-y),\; \frac{\si}{6}(3+3x-y)\right).
\leqno{(\ref{prop:inverse_design}a)}$$ 
If $(x,y)$ is in the interior of $\QT$, the invariant $\RI$ defines a pair of lattices $\La^{\pm}$ that have opposite signs and 
unique (up to isometry) reduced basis vectors $v_1,v_2$ with the lengths $|v_1|=\sqrt{r_{12}^2+r_{01}^2}$, $|v_2|=\sqrt{r_{12}^2+r_{02}^2}$ and the anticlockwise angle
$(\ref{prop:inverse_design}b)\qquad \angle(v_1,v_2)=\arccos\dfrac{-4y^2}{\sqrt{(9x^2+5y^2-6y+9)^2-36x^2(3-y)^2}}$.
\bt
\end{prop}
\begin{proof}
In Definition~\ref{dfn:PF} the projected invariant $\PI(\La)=(x,y)$ is obtained from the coordinates $(\bar r_{12},\bar r_{01},\bar r_{02})$ of $\PRF(\La)$ satisfying the equations 
$\left\{\begin{array}{l}
x=\bar r_{02}-\bar r_{01}, \\
y=3\bar r_{12}, \\
\bar r_{12}+\bar r_{01}+\bar r_{02}=1.
\end{array}\right.$
The solution is $\PRF(\La)=(\bar r_{12},\bar r_{01},\bar r_{02})=(\frac{y}{3},\frac{1}{2}-\frac{x}{2}-\frac{y}{6},\frac{1}{2}+\frac{x}{2}-\frac{y}{6})$.
Multiplying all coordinates by the size $\si=r_{12}+r_{01}+r_{02}$ gives the root invariant in (\ref{prop:inverse_design}a). 
\medskip

By Proposition~\ref{prop:reduced_bases} a (unique up to isometry) reduced basis $v_1,v_2$ consist of two shortest vectors of an obtuse superbase.
Lemma~\ref{lem:superbase_reconstruction} implies the formulae for $|v_1|,|v_2|$.
The angle formula from Lemma~\ref{lem:superbase_reconstruction} can be expressed in $x,y$ as follows:
$$\cos\angle(v_1,v_2)=\dfrac{-r_{12}^2}{|v_1|\cdot|v_2|}
=\dfrac{-y^2/9}{\frac{1}{6}\sqrt{(2y)^2+(3-3x-y)^2}\frac{1}{6}\sqrt{(2y)^2+(3+3x-y)^2}}=$$
$$=\dfrac{-4y^2}{\sqrt{4y^2+(3-y)^2+9x^2+6x(3-y)}\sqrt{4y^2+(3-y)^2+9x^2-6x(3-y)}}=$$
$=\dfrac{-4y^2}{\sqrt{(9x^2+5y^2-6y+9)+6x(3-y)}\sqrt{(9x^2+5y^2-6y+9)-6x(3-y)}}$.
\ws
\end{proof}

Example~\ref{exa:inverse_design} shows the power of Proposition~\ref{prop:inverse_design} based on Theorem~\ref{thm:classification2d} and Corollary~\ref{cor:similar_lattices2d} for inverse design by sampling the square $\QS$ at interesting places.
\medskip

Fig.~\ref{fig:QS+DC}~(right) visualises the doubled cone $\DC$ of oriented root invariants $\RI^o$ from Definition~\ref{dfn:sign} 
 by uniting the triangular cone $\TC=\{0\leq r_{12}\leq r_{01}\leq r_{02}\}$ with its mirror reflection in the vertical plane $\{r_{01}=r_{02}\}$ including the $r_{12}$-axis.
\medskip

The lattice $L_0$ with $\RI=(1,1,4)$ is represented by two boundary points of $\DC$ identified by $(r_{01},r_{02})\lra(r_{02},r_{01})$.
The lattices $L_{\infty}^{\pm}$ with the root invariant $\RI=(r_{12},r_{01},r_{02})=(1,4,7)$ are represented by $(1,4,7)$ and its mirror image $(1,7,4)$ in $\DC$ related by the reflection in the vertical bisector plane $r_{01}=r_{02}$ containing the root invariants of $\La_4,\La_6$.
The superscript shows $\sign(L_{\infty}^{\pm})=\pm 1$. 
\medskip
  
The lattice notations $L_0,L_2^{\pm},L_{\infty}^{\pm}$ slightly differ from $\La_4,\La_6$ in Fig.~\ref{fig:QS+DC}, because the companion paper \cite{bright2021geographic} considers other lattices $\La_0,\La_2^{\pm},\La_{\infty}^{\pm}$ at similar positions in the quotient triangle $\QT$ but with different coordinates $(\xi,\eta)$ from Remark~\ref{rem:group_action}.

\begin{exa}[inverse design of 2D lattices]
\label{exa:inverse_design}
We will inversely design the lattices $\La_4,\La_6,L_0,L_2^{\pm},L_{\infty}^\pm$, see their visualised invariants in Fig.~\ref{fig:QS+DC}~(right).
\medskip

\noindent
\emph{($\pmb{\La_4}$)}
We design the square lattice $\La_4$ starting from its projected invariant  at the origin $\PI(\La_4)=(0,0)\in\QT$, which is identified with the top right vertex $(1,1)\in\QS$ in Fig.~\ref{fig:QS+DC}~(left).
Formula~(\ref{prop:inverse_design}a) for the size $\si(\La_4)=2$ (only to get simplest integers) gives $\RI(\La_4)=(0,1,1)$.
An obtuse superbase $\{v_0,v_1,v_2\}$ can be reconstructed by Lemma~\ref{lem:superbase_reconstruction}.
The vonorms are $v_1^2=v_2^2=0^2+1^2=1$, $v_0^2=1^2+1^2=2$.
We can choose the standard obtuse superbase $v_1=(1,0)$, $v_2=(0,1)$, $v_0=(-1,-1)$.
\medskip

\noindent
\emph{($\pmb{\La_6}$)}
We design the hexagonal lattice $\La_6$ starting from the projected invariant at the top left vertex $\PI(\La_6)=(0,1)\in\QT$, which is identified with the bottom right vertex $(1,0)\in\QS$ in Fig.~\ref{fig:QS+DC}~(left).
Formula~(\ref{prop:inverse_design}a) for the size $\si(\La_6)=3$ (only to get simplest integers) gives $\RI(\La_6)=(1,1,1)$.
To reconstruct an obtuse superbase $\{v_0,v_1,v_2\}$ by Lemma~\ref{lem:superbase_reconstruction}, find the vonorms $v_1^2=v_2^2=v_0^2=1^2+1^2=2$.
Formula~(\ref{prop:inverse_design}b) gives the angle $\angle(v_1,v_2)=\arccos
\frac{-4}{\sqrt{(5-6+9)^2}}=
\arccos\left(-\frac{1}{2}\right)=120^\circ$.
We can choose the superbase $v_1=(\sqrt{2},0)$, $v_2=(-\frac{1}{\sqrt{2}},\frac{\sqrt{3}}{\sqrt{2}})$,
$v_0=(-\frac{1}{\sqrt{2}},-\frac{\sqrt{3}}{\sqrt{2}})$.
\medskip

\begin{figure}[h]
\includegraphics[height=45mm]{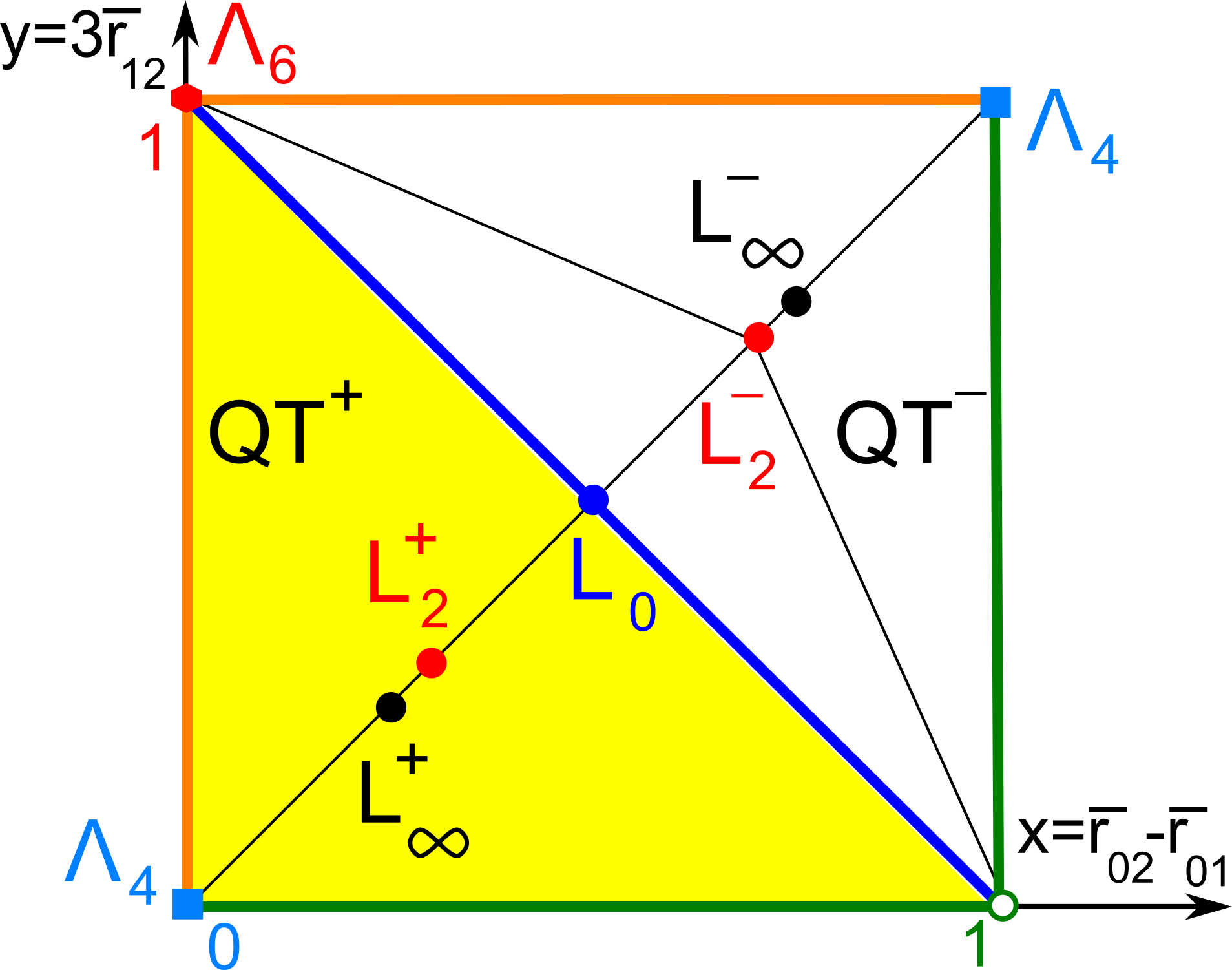}
\hspace*{2mm}
\includegraphics[height=45mm]{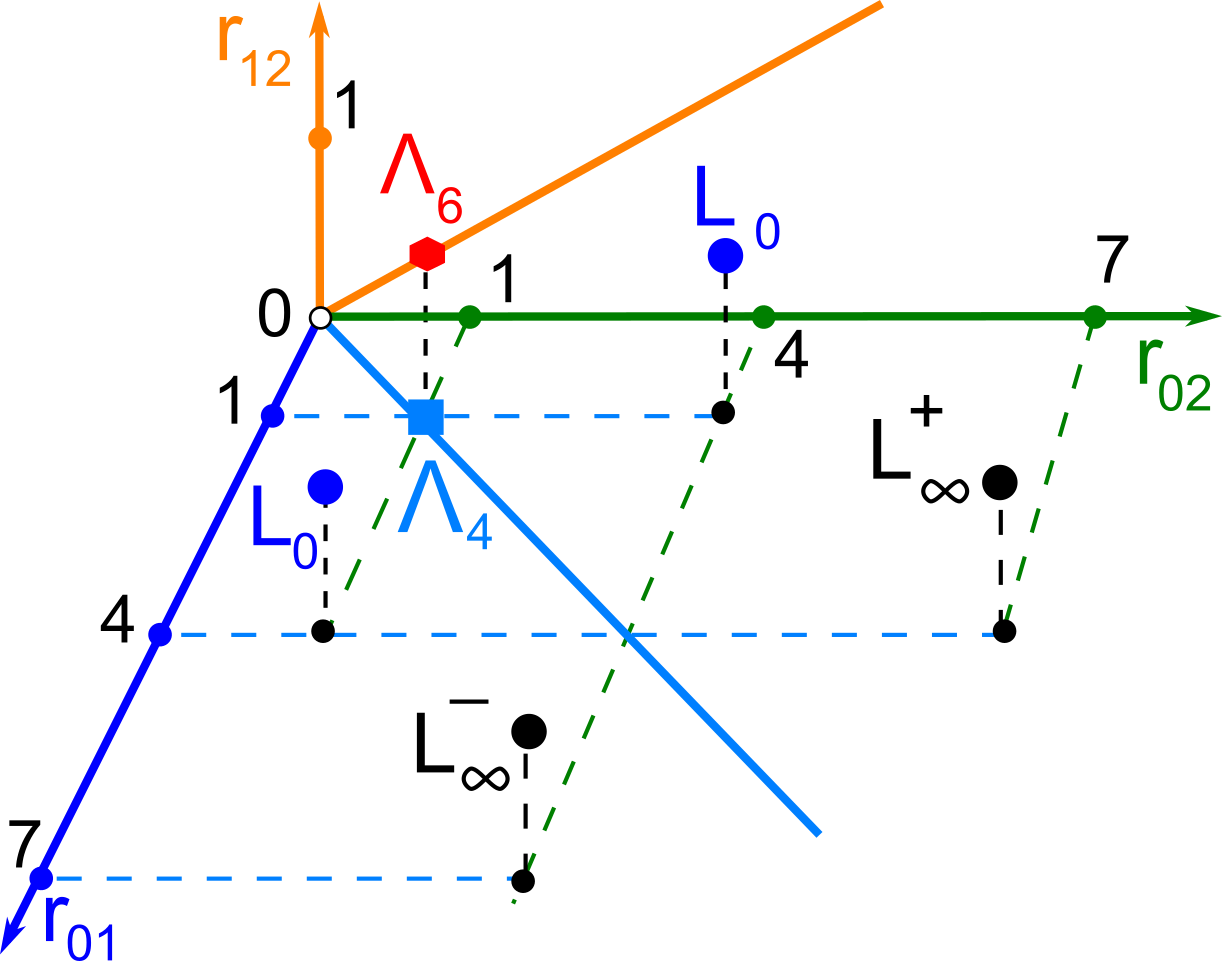}
\caption{\textbf{Left}: $\QS=\QT^+\cup\QT^-$ includes mirror-symmetric lattices $\La_4,\La_6,L_0$ and non-mirror-symmetric lattices $L_{\infty}^{\pm}$, see Example~\ref{exa:RM} and Table~\ref{tab:forms460inf} later.
\textbf{Right}: the doubled cone $\DC$ is visualised as $\{0\leq r_{12}\leq\min\{r_{01},r_{02}\}>0\}$ bounded by the planes $\{r_{12}=0\}$, $\{r_{12}=r_{01}\}$, $\{r_{12}=r_{02}\}$ with the identifications $(r_{12},r_{01},r_{02})\lra (r_{12},r_{02},r_{01})$ on the boundary $\bd\DC$.
}
\label{fig:QS+DC}
\end{figure}

\noindent
\emph{($\pmb{L_0}$)}
We inversely design the lattice $L_0$ in Fig.~\ref{fig:QS+DC} starting from $\PI(L_0)=(x,y)$ at the centre $(\frac{1}{2},\frac{1}{2})\in\QS$.
Formula~(\ref{prop:inverse_design}a) for the size $\si(L_0)=6$ (only to get simplest integers) gives $\RI(L_0)=(1,1,4)$.
To reconstruct an obtuse superbase $\{v_0,v_1,v_2\}$ by Lemma~\ref{lem:superbase_reconstruction}, find the vonorms $v_1^2=1^2+1^2=2$, $v_0^2=v_2^2=1^2+4^2=17$.
Formula~(\ref{prop:inverse_design}b) gives the angle
$\angle(v_1,v_2)=\arccos
\dfrac{-4}{\sqrt{(\frac{9}{4}+\frac{5}{4}-3+9)^2-9(\frac{5}{4})^2}}
=\arccos(-\frac{1}{\sqrt{34}})\approx 99.9^\circ$.
We can choose the following superbase, see Fig.~\ref{fig:DT}:  
$v_1=(\sqrt{2},0)$, 
$v_2=|v_2|(\cos\angle(v_1,v_2),\sin\angle(v_1,v_2))
=(-\frac{1}{\sqrt{2}},\frac{\sqrt{33}}{\sqrt{2}})$,
$v_0=(-\frac{1}{\sqrt{2}},-\frac{\sqrt{33}}{\sqrt{2}})$.
\medskip

\noindent
\emph{($\pmb{L_{2}}$)}
We inversely design the lattice $L_2$ in Fig.~\ref{fig:QS+DC} 
starting from their projected invariants $\PI(L_2)=(\frac{1}{2+\sqrt{2}},\frac{1}{2+\sqrt{2}})$, which will maximise the chiral distance $\PC[D_2]$ in Theorem~\ref{prop:PC}(a).
Formula~(\ref{prop:inverse_design}a) for the size $\si(L_2)=6$ (only to simplify the root invariant) gives $\RI(L_{2})=(2-\sqrt{2},2\sqrt{2}-1,5-\sqrt{2})$.
Since all root products are non-zero and distinct, by Lemma~\ref{lem:achiral_lattices} there is a pair of lattices $L_{2}^\pm$ with $\sign(L_{2}^\pm)=\pm 1$.
The lattices $L_{2}^\pm$ are related by reflection, not by rigid motion.
\medskip

To reconstruct an obtuse superbase $\{v_0,v_1,v_2\}$ of $L_{2}^\pm$ by Lemma~\ref{lem:superbase_reconstruction}, find \\
$v_0^2=(2\sqrt{2}-1)^2+(5-\sqrt{2})^2=(9-4\sqrt{2})+(27-10\sqrt{2})=36-14\sqrt{2}\approx 16.2,$ \\
$v_1^2=(2-\sqrt{2})^2+(2\sqrt{2}-1)^2=(6-4\sqrt{2})+(9-4\sqrt{2})=15-8\sqrt{2}\approx 3.7,$ \\
$v_2^2=(2-\sqrt{2})^2+(5-\sqrt{2})^2=(6-4\sqrt{2})+(27-10\sqrt{2})=33-14\sqrt{2}\approx 13.2,$  \\
and the anticlockwise angle
$\angle(v_1,v_2)=\arccos
\dfrac{-r_{12}^2}{|v_1|\cdot|v_2|}\approx 92.8^\circ$.
Then $L_2^\pm$ have the following obtuse superbases in
Fig.~\ref{fig:DT}:  
$v_1=(\sqrt{15-8\sqrt{2}},0)\approx(1.9,0)$, 
$v_2=|v_2|(\cos\angle(v_1,v_2),\sin\angle(v_1,v_2))
\approx(-0.18,3.63)$,
$v_0
\approx(-1.72,-3.63)$.
\medskip

\begin{figure}[h]
\includegraphics[width=\textwidth]{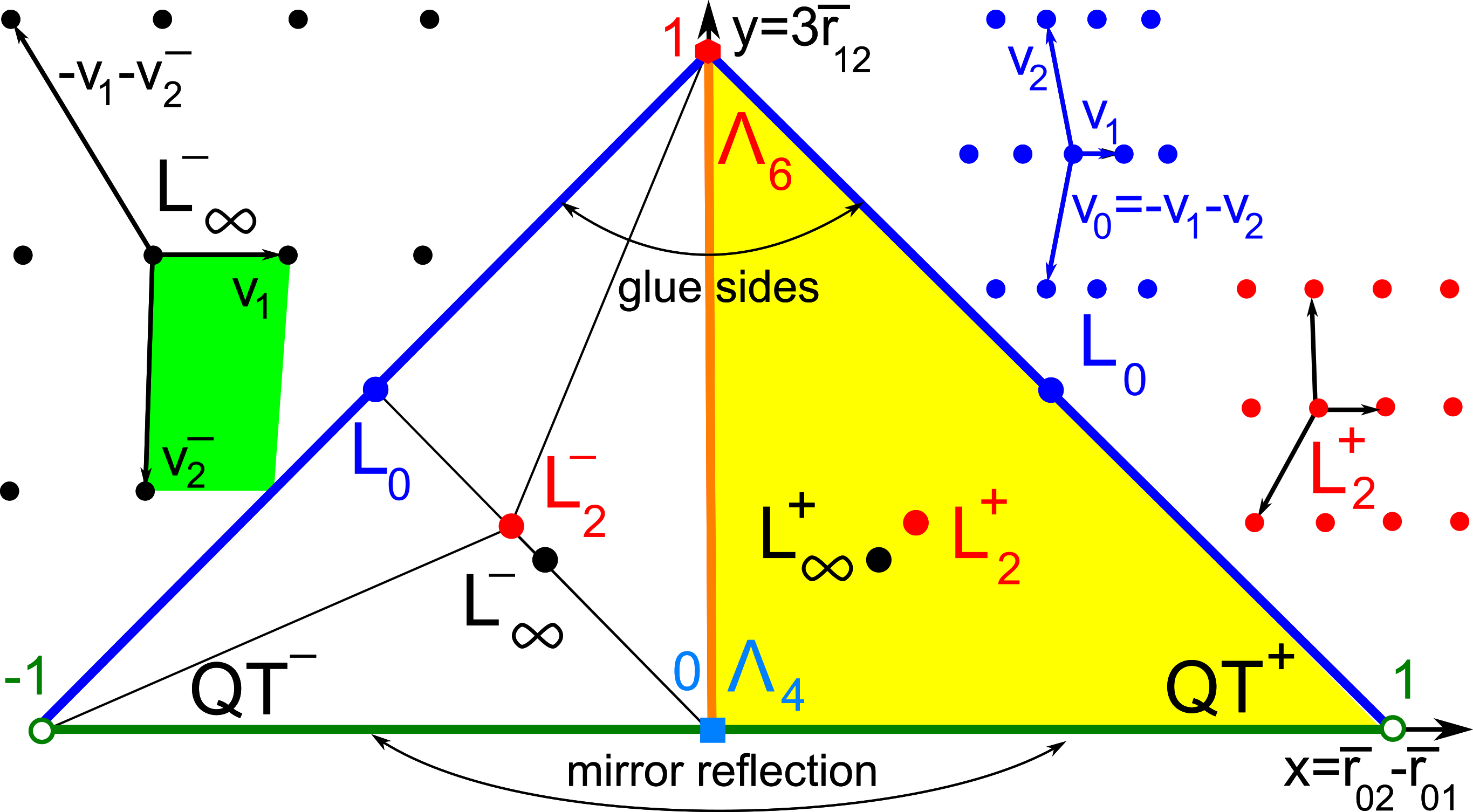}
\caption{
The doubled cone $\DC$ in Fig.~\ref{fig:QS+DC} (right) projects to the doubled triangle $\DT$ parameterised by $x\in(-1,1)$, $y\in[0,1]$ and obtained by gluing two copies $\QT^{\pm}$ of the quotient triangle along vertical sides instead of hypotenuses as in $\QS$, see Example~\ref{exa:inverse_design} and Table~\ref{tab:forms460inf}.
}
\label{fig:DT}
\end{figure}

\noindent
\emph{($\pmb{L_{\infty}}$)}
We inversely design the lattice $L_{\infty}$ in Fig.~\ref{fig:DT} starting from $\PI(L_{\infty})=(x,y)$ at the mid-point $(\frac{1}{4},\frac{1}{4})$ of the segment between $\PI(\La_4),\PI(L_0)\in\QT$.
Formula~(\ref{prop:inverse_design}a) for the size $\si(L_{\infty})=12$ (only to simplify the root invariant) gives 
 $\RI(L_{\infty})=(1,4,7)$.
Since all root products are non-zero and distinct, by Lemma~\ref{lem:achiral_lattices} there is a pair of lattices $L_{\infty}^\pm$ of opposite signs $\sign(L_{\infty}^\pm)=\pm 1$.
\medskip

To reconstruct an obtuse superbase $\{v_0,v_1,v_2\}$ of $L_{\infty}^\pm$ by Lemma~\ref{lem:superbase_reconstruction}, find the vonorms $v_0^2=4^2+7^2=65$, $v_1^2=1^2+4^2=17$, $v_2^2=1^2+7^2=50$, and the anticlockwise angle
$\angle(v_1,v_2)=\arccos
\dfrac{-r_{12}^2}{|v_1|\cdot|v_2|}
=\arccos(-\frac{1}{\sqrt{850}})\approx 92^\circ$.
Then $L_{\infty}^\pm$ have the following obtuse superbases in
Fig.~\ref{fig:DT}:  
$v_1=(\sqrt{17},0)\approx(4.12,0)$, 
$$v_2^{\pm}=|v_2|(\cos\angle(v_1,v_2),\sin\angle(v_1,v_2))
=\left(-\frac{1}{\sqrt{17}},\pm\frac{\sqrt{849}}{\sqrt{17}}\right)\approx(-0.24,\pm7.1),$$
$v_0^\pm=-v_1-v_2^{\pm}=(-\frac{16}{\sqrt{17}},\mp\frac{\sqrt{849}}{\sqrt{17}})\approx(-3.88,\mp7.1)$, see all forms in
Table~\ref{tab:forms460inf}.
\bs
\end{exa}

\begin{table}[h]
\caption{Various forms of the lattices computed in Example~\ref{exa:inverse_design} and shown Fig.~\ref{fig:QS+DC} and~\ref{fig:DT}. }
\medskip

\hspace*{-2mm}
\begin{tabular}{l|ccccc}      
$\La$ & 
$\La_4$ & 
$\La_6$ & 
$L_0$ & 
$L_2^{\pm}$ & 
$L_{\infty}^{\pm}$ \\
\hline

$\si(\La)$ & 2 & 3 & 6 & 6 & 12 \\

$\PI(\La)$ 
& $(0,0)$
& $(0,1)$
& $\left(\dfrac{1}{2},\dfrac{1}{2}\right)$
& $\left(\dfrac{1}{2+\sqrt{2}},\dfrac{1}{2+\sqrt{2}}\right)$
& $\left(\dfrac{1}{4},\dfrac{1}{4}\right)$ \\


$\RI^o(\La)$  & 
(0,1,1) & 
(1,1,1) & 
(1,1,4) & 
$(2-\sqrt{2},2\sqrt{2}-1,5-\sqrt{2})^{\pm}$ & 
$(1,4,7)^{\pm}$ \\

$\VF(\La)$  & 
(2,1,1) & 
(2,2,2) & (17,2,17) & 
$(15-8\sqrt{2},33-14\sqrt{2},36-14\sqrt{2})$ & 
(65,17,50)
\end{tabular}
\label{tab:forms460inf}
\end{table}

\section{Metrics on spaces of lattices up to isometry, rigid motion, similarity}
\label{sec:metrics}

All lattices $\La\subset\R^2$ are uniquely represented up to isometry and similarity by their invariants $\RI\in\TC$ and $\PI\in\QT$, respectively.
Then any metric $d$ on the triangular cone $\TC\subset\R^3$ or the quotient triangle $\QT\subset\R^2$ gives rise to a metric in Definition~\ref{dfn:RM} on the spaces $\LIS$ and $\LSS$, respectively. 
The oriented case in Definition~\ref{dfn:RMo} will be harder because of identifications on the boundary $\bd\TC$.

\begin{dfn}[root metrics $\RM$, projected metrics $\PM$]
\label{dfn:RM}
Any metric $d$ on $\R^3$ defines the \emph{root metric} 
$\RM(\La_1,\La_2)=d(\RI(\La_1),\RI(\La_2))$ on lattices $\La_1,\La_2\subset\R^2$  up to isometry.
The \emph{Root Invariant Space} $\RIS=(\TC,d)$ is the triangular cone with a fixed metric $d$.
If we use the Minkowski norm $M_q(v)=||v||_q=(\sum\limits_{i=1}^n |x_i|^q)^{1/q}$ of a vector $v=(x_1,\dots,x_n)\in\R^n$  for any real $q\in[1,+\infty]$, the root metric is denoted by $\RM_q(\La_1,\La_2)=||\RI(\La_1)-\RI(\La_2)||_q$.
The limit case $q=+\infty$ uses 
$||v||_{\infty}=\max\limits_{i=1,\dots,n}|x_i|$.
The \emph{projected metric} $\PM(\La_1,\La_2)=d(\PI(\La_1),\PI(\La_2))$ is on lattices up to similarity for any metric $d$ on $\R^2$.
The space of \emph{projected invariants} $\PIN=(\QT,d)$ is the quotient triangle with a metric $d$.
The notation $\PM_q(\La_1,\La_2)=||\PI(\La_1)-\PI(\La_2)||_q$ includes a parameter $q\in[1,+\infty]$ of $M_q$.
\bs
\end{dfn}

The Minkowski distance $M_q$ for $q=2$ is Euclidean.
The root metric $\RM_q$ can take any large values in original units of vector coordinates such as Angstroms.
The projected metric $\PM_q$ is unitless and the space $\PIN=(\QT,d)$ is bounded.

\begin{table}[h]
\caption{Metrics $\RM_q$ and $\PM_q$ for the lattices from Example~\ref{exa:RM} and shown Fig.~\ref{fig:QS+DC} and~\ref{fig:DT}. }
\medskip

\begin{tabular}{|l|cccc|}  
\hline    
$\RM_{\infty}$ & $\La_4$ & $\La_6$ & $L_0$ & $L_{\infty}^{\pm}$\\
\hline

$\RI(\La_4)=(0,1,1)$ & $0$ & $1$ & $3$ & $6$ \\

$\RI(\La_6)=(1,1,1)$ & $1$ & $0$ & $3$ & $6$ \\

$\RI(L_0)=(1,1,4)$ & $3$ & $3$ & $0$ & $3$  \\

$\RI(L_{\infty}^{\pm})=(1,4,7)$ & $6$ & $6$ & $3$ & $0$  \\
\hline
\end{tabular}
\begin{tabular}{|l|cccc|}
\hline      
$\PM_{\infty}$ & $\La_4$ & $\La_6$ & $L_0$ & $L_{\infty}^{\pm}$\\
\hline

$\PI(\La_4)=(0,0)$ & $0$ & $1$ & $\frac{1}{2}$ & $\frac{1}{4}$ \\

$\PI(\La_6)=(0,1)$ & $1$ & $0$ & $\frac{1}{2}$ & $\frac{3}{4}$ \\

$\PI(L_0)=(\frac{1}{2},\frac{1}{2})$ & $\frac{1}{2}$ & $\frac{1}{2}$ & $0$ & $\frac{1}{4}$ \\

$\PI(L_{\infty}^{\pm})=(\frac{1}{4},\frac{1}{4})$ & $\frac{1}{4}$ & $\frac{3}{4}$ & $\frac{1}{4}$ & $0$  \\
\hline
\end{tabular}
\medskip

\begin{tabular}{|l|cccc|}  
\hline    
$\RM_q$ for $q\in[1,+\infty)$ & $\La_4$ & $\La_6$ & $L_0$ & $L_{\infty}^{\pm}$\\
\hline

$\RI(\La_4)=(0,1,1)$ & $0$ & $1$ & $(1+3^q)^{1/q}$ & $(1+3^q+6^q)^{1/q}$ \\

$\RI(\La_6)=(1,1,1)$ & $1$ & $0$ & $3$ & $(3^q+6^q)^{1/q}$ \\

$\RI(L_0)=(1,1,4)$ & $(1+3^q)^{1/q}$ & $3$ & $0$ & $3\cdot 2^{1/q}$  \\

$\RI(L_{\infty}^{\pm})=(1,4,7)$ & $(1+3^q+6^q)^{1/q}$ & $(3^q+6^q)^{1/q}$ & $3\cdot 2^{1/q}$ & $0$  \\
\hline
\end{tabular}
\medskip

\begin{tabular}{|l|cccc|}
\hline      
$\PM_q$ for $q\in[1,+\infty)$ & $\La_4$ & $\La_6$ & $L_0$ & $L_{\infty}^{\pm}$\\
\hline

$\PI(\La_4)=(0,0)$ & $0$ & $1$ & $2^{(1/q)-1}$ & $2^{(1/q)-2}$ \\

$\PI(\La_6)=(0,1)$ & $1$ & $0$ & $2^{(1/q)-1}$ & $\frac{1}{4}(1+3^q)^{1/q}$ \\

$\PI(L_0)=(\frac{1}{2},\frac{1}{2})$ & $2^{(1/q)-1}$ & $2^{(1/q)-1}$ & $0$ & $2^{(1/q)-2}$  \\

$\PI(L_{\infty}^{\pm})=(\frac{1}{4},\frac{1}{4})$ & $2^{(1/q)-2}$ & $\frac{1}{4}(1+3^q)^{1/q}$ & $2^{(1/q)-2}$ & $0$  \\
\hline
\end{tabular}
\label{tab:RM}
\end{table}

\begin{exa}[metrics $\RM_q,\PM_q$]
\label{exa:RM}
Table~\ref{tab:RM} summarises metric computations for the lattices $\La_4,\La_6,L_0,L_{\infty}^{\pm}$, which  
were inversely designed in Example~\ref{exa:inverse_design}.
\bs
\end{exa}

\begin{lem}[metric axioms for $\RM,\PM$]
\label{lem:RM}
\textbf{(a)}
Any metrics $\RM$ and $\PM$ from Definition~\ref{dfn:RM} satisfy all metric axioms in (\ref{pro:map}c) on the Lattice Isometry Space $\LIS(\R^2)$ and the Lattice Similarity Space $\LSS(\R^2)$, respectively.
\bt
\end{lem}
\begin{proof}
The metric axioms for $\RM,\PM$ from Definition~\ref{dfn:RM} follow from the same axioms for an underlying metric $d$.
Only the first axiom is non-trivial: by the first axiom for $d$ we know that $\RM(\La_1,\La_2)=d(\RI(\La_1),\RI(\La_2))=0$ if and only if $\RI(\La_1)=\RI(\La_2)$.
Now Theorem~\ref{thm:classification2d} says that $\RI(\La_1)=\RI(\La_2)$ is equivalent to $\La_1,\La_2$ being isometric.
Corollary~\ref{cor:similar_lattices2d} classifying lattices up to similarity by projected invariants similarly justifies the first axiom for $\PM(\La_1,\La_2)$. 
\ws
\end{proof}

Since the mirror images $L_{\infty}^{\pm}$ have the same root invariant $\RI(L_{\infty}^{\pm})=(1,4,7)$, for any lattice $\La$, the distances $\RM(\La,L_{\infty}^{\pm})$ and $\PM(\La,L_{\infty}^{\pm})$ are independent of $\sign(L_{\infty}^{\pm})=\pm 1$.
Any mirror images $\La^{\pm}$ have $\RM(\La^+,\La^-)=0=\PM(\La^+,\La^-)$ because $\La^{\pm}$ are isometric to each other.
The metric $\RM$ from Definition~\ref{dfn:RM} is well-defined only for lattices up to any isometry including reflections.
\medskip

\begin{figure}[h]
\includegraphics[width=1.0\textwidth]{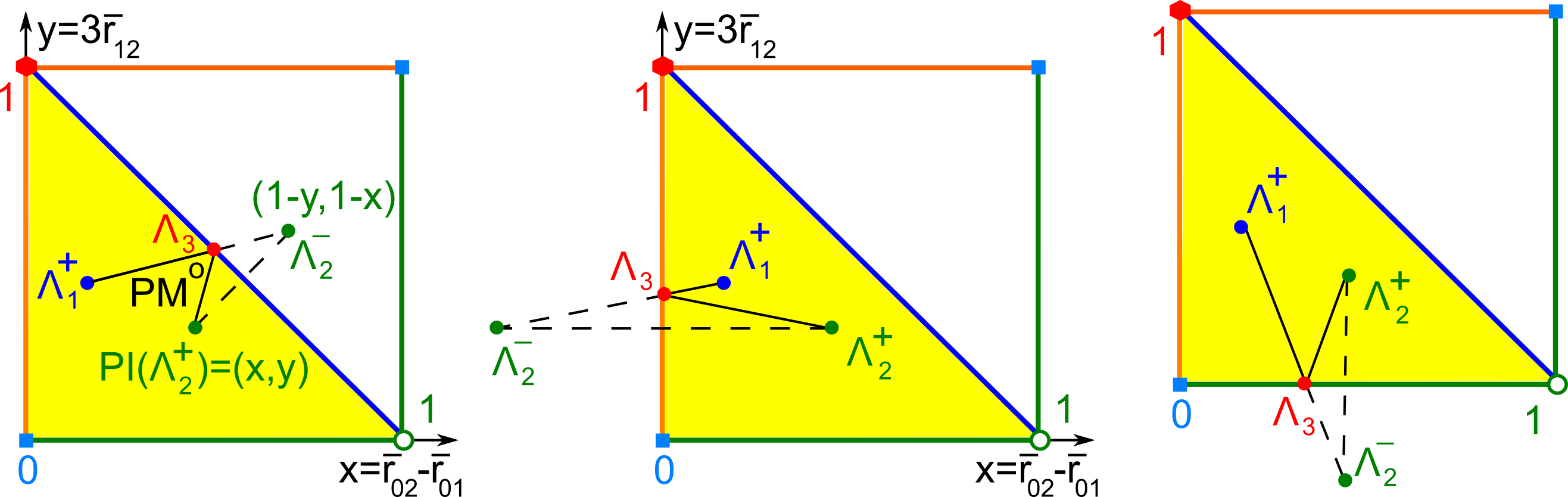}
\caption{By Definition~\ref{dfn:RMo}, the projected metric $\PM_2(\La_1^+,\La_2^-)$ is the minimum sum $\PM_2(\La_1^+,\La_3)+\PM_2(\La_3\La_2^-)$ achieved in the left image, see computations in Proposition~\ref{prop:RMq}.}
\label{fig:QS+RMo}
\end{figure}

Definition~\ref{dfn:RMo} introduces the metric $\RM^o$ on lattices up to rigid motion so that $\RM^o(\La^+,\La^-)>0$ on mirror images of a non-mirror-symmetric lattice, see Fig.~\ref{fig:QS+RMo}.

\begin{dfn}[orientation-aware metrics $\RM^o,\PM^o$]
\label{dfn:RMo}
For lattices $\La_1,\La_2\subset\R^2$ with $\sign(\La_1)\sign(\La_2)\geq 0$, the \emph{orientation-aware} root metric is $\RM^o(\La_1,\La_2)=\RM(\La_1,\La_2)$ as in Definition~\ref{dfn:RM}.
If any lattices $\La_1,\La_2$ have opposite signs, set $\RM^o(\La_1,\La_2)=\inf\limits_{\sign(\La_3)=0}(\RM(\La_1,\La_3)+\RM(\La_2,\La_3))$.
The \emph{orientation-based} metric $\PM^o(\La_1,\La_2)$ is defined by the same formula, where we replace $\RM$ by $\PM$.
\bs
\end{dfn}

The infimum in $\RM^o(\La_1,\La_2)$ is the greatest lower bound as in \cite[Lemma~I.5.24]{bridson2013metric} defining a metric on a union of metric spaces glued by isometries.
Theoretically, this bound may not be achieved over a non-compact domain.
When using a Minkowski base metric $M_q$, Propositions~\ref{prop:RMq}-\ref{prop:PMq} explicitly compute $\RM_q^o,\PM_q^o$ for $q=2,+\infty$, so the infimum in Definition~\ref{dfn:RMo} can be replaced by a minimum in practice.
\medskip

The \emph{oriented} root invariant space $\RIS^o$ and the space of \emph{oriented} projected invariants $\PIN^o$ can be defined similarly to $\RIS$ and $\PIN$ in Definition~\ref{dfn:RM} as the doubled cone $\DC$ and quotient square $\QS$ with any metrics from Definition~\ref{dfn:RMo}.

\begin{lem}[metric axioms for $\RM^o,\PM^o$]
\label{lem:RMo}
Any root metric $\RM^o$ and projected metric $\PM^o$ from Definition~\ref{dfn:RMo} satisfy all metric axioms in~(\ref{pro:map}c)
 on the Lattice Isometry Space $\LIS^o$ and the Lattice Similarity Space $\LSS^o$, respectively.
\bt
\end{lem}
\begin{proof}
Lemma~\ref{lem:RM} implies metric axioms in all cases when involved lattices have the same sign or sign 0.
For example, the metrics $\RM^o,\PM^o$ vanish only if the lattices $\La_1,\La_2$ in question have equal forms $\RI^o,\PI^o$, respectively, so $\La_1,\La_2$ are isometric or similar by Theorem~\ref{thm:classification2d} or Corollary~\ref{cor:similar_lattices2d}, respectively.
The symmetry axiom for any lattices $\La,\La'$ with opposite signs follows from the symmetry of the sum in Definition~\ref{dfn:RMo}: $\RM^o(\La_1,\La_2)=\min\limits_{\sign(\La_3)=0}(\RM(\La_1,\La_3)+\RM(\La_2,\La_3))$.
\medskip

Without loss of generality it suffices to prove the required triangle inequality
$\RM^o(\La_1,\La_2)+\RM^o(\La_2,\La_3)\geq \RM^o(\La_1,\La_3)$ in the following two cases below.
\medskip

\noindent
\emph{Case} $\sign(\La_1)\geq 0$ and $\sign(\La_2)\geq 0>\sign(\La_3)$.
Then $\RM^o(\La_1,\La_2)=\RM(\La_1,\La_2)$ is the root metric without minimisation from Definition~\ref{dfn:RM}.
Let $\La'$ be some achiral lattice minimising $\RM^o(\La_2,\La_3)=\min\limits_{\sign(\La')=0}(\RM(\La_2,\La')+\RM(\La_3,\La'))$.
Then 
$$\RM^o(\La_1,\La_2)+\RM^o(\La_2,\La_3)=\RM(\La_1,\La_2)+\RM(\La_2,\La')+\RM(\La_3,\La')\geq$$ 
$$\geq \RM(\La_1,\La')+\RM(\La_3,\La')\geq\min\limits_{\sign(\La')=0}(\RM(\La_1,\La')+\RM(\La_3,\La'))=\RM^o(\La_1,\La_3),$$ where we used the triangle inequality 
for the root metric $\RM$ and $\La_1,\La',\La_2$.
\medskip

\noindent
\emph{Case} $\sign(\La_1)\geq 0$ and $\sign(\La_3)\geq 0>\sign(\La_2)$.
Then $\RM^o(\La_1,\La_3)=\RM(\La_1,\La_3)$ is the root metric without minimisation from Definition~\ref{dfn:RM}.
Let $\La',\La''$ be mirror-symmetric lattices minimising 
$\RM^o(\La_1,\La_2)=\min\limits_{\sign(\La')=0}(\RM(\La_1,\La')+\RM(\La_2,\La'))$ and 
$\RM^o(\La_2,\La_3)=\min\limits_{\sign(\La'')=0}(\RM(\La_2,\La'')+\RM(\La_3,\La''))$, respectively.
Then 
$$\RM^o(\La_1,\La_2)+\RM^o(\La_2,\La_3)=\RM(\La_1,\La')+\RM(\La_2,\La')+\RM(\La_2,\La'')+\RM(\La_3,\La'')$$
$$\geq \RM(\La_1,\La')+\RM(\La',\La'')+\RM(\La_3,\La'')\geq \RM(\La_1,\La_3),\text{ where} $$
we used the triangle inequality for $\RM$ and the lattices $\La_2,\La',\La''$ with non-positive signs, then for the lattices $\La_1,\La',\La'',\La_3$, which have only non-negative signs.
\ws
\end{proof}

Lemma~\ref{lem:reversed_signs} speeds up computations in the oriented case, see Example~\ref{exa:RMo}.

\begin{lem}[reversed signs]
\label{lem:reversed_signs}
If lattices $\La_1^{\pm},\La_2^{\pm}\subset\R^2$ have specified signs, then 
$\RM^o(\La_1^+,\La_2^-)=\RM^o(\La_1^-,\La_2^+)$ and
$\PM^o(\La_1^+,\La_2^-)=\PM^o(\La_1^-,\La_2^+)$.
\bt
\end{lem}
\begin{proof}
By Definition~\ref{dfn:RMo}, for any base distance $d$ on $\R^3$, when minimising over mirror-symmetric lattices $\La_3$ with $\sign(\La_3)=0$, the metrics $\RM^o$ are computed for lattices that have one zero sign and one non-zero sign.
Hence $\RM^o(\La_1^{\pm},\La_3)$ can be replaced by the simpler metric $\RM(\La_1,\La_3)=d(\RI(\La_1),\RI(\La_3))$ depending only on the unoriented root invariants $\RI(\La_1)$ and $\RI(\La_3)$ without a sign.
After that the metric $\RM$ can lifted back to the lattices $\La_1^-,\La_3^+$ with reversed signs:
$$\RM^o(\La_1^+,\La_2^-)
=\min\limits_{\sign(\La_3)=0}(\RM(\La_1^+,\La_3)+\RM(\La_2^-,\La_3))=$$
$$=\min\limits_{\sign(\La_3)=0}(d(\RI(\La_1),\RI(\La_3))+d(\RI(\La_2),\RI(\La_3))=$$
$$=\min\limits_{\sign(\La_3)=0}(\RM(\La_1^-,\La_3)+\RM(\La_2^+,\La_3))=\RM^o(\La_1^-,\La_2^+).$$
The proof for the projected metric $\PM^o$ is similar to the above arguments. 
\ws
\end{proof}

Lemma~\ref{lem:max_sum_moduli} will help compute $\RM_q^o,\PM_q^o$ for $q=2,+\infty$ in Propositions~\ref{prop:RMq},~\ref{prop:PMq}.

\begin{lem}[maximum sum of moduli]
\label{lem:max_sum_moduli}
For any real numbers $a,b,c,d\geq 0$, 
the maximum sum $S(x)=\max\{|a-x|+|x-b|,|c-x|+|x-d|\}$ over $x\geq 0$ has the minimum value 
$\MS(a,b,c,d)=\max\{|a-b|,|c-d|,\frac{1}{2}|a+b-c-d|\}$.
\bt
\end{lem}
\begin{proof}
The maximum sum $S(x)$ and the formula for $\MS(a,b,c,d)$ are invariant under permutations $a\lra b$, $c\lra d$ and $(a,b)\lra(c,d)$.
Without loss of generality one can assume that
$a\leq b$ and $a\leq c\leq d$.
Then Fig.~\ref{fig:max_sum_moduli1} shows the graphs of $y=|a-x|+|x-b|$, $y=|c-x|+|x-d|$, $y=S(x)$ in green, blue, red, respectively.

\begin{figure}[h]
\includegraphics[width=1.0\textwidth]{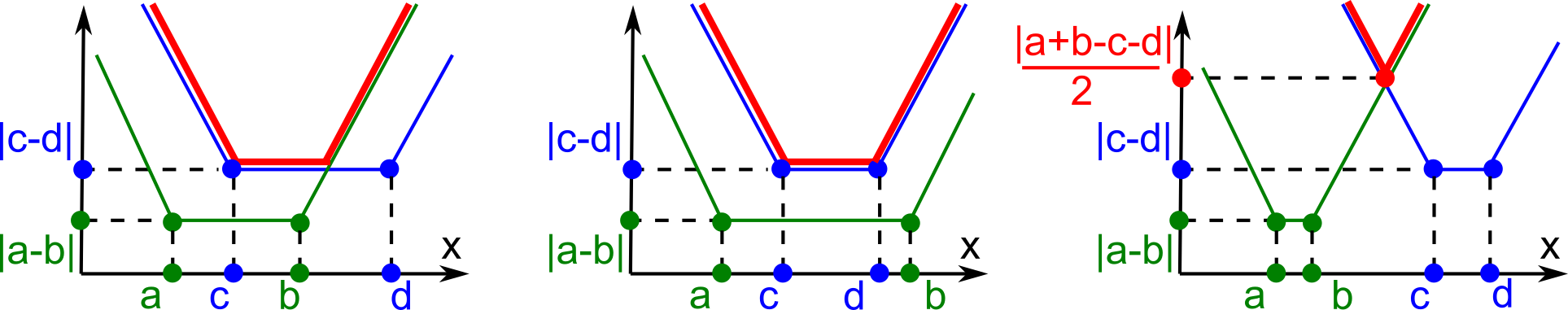}
\caption{Lemma~\ref{lem:max_sum_moduli} finds a minimum value of $S(x)$ for all different positions of $a,b,c,d$.}
\label{fig:max_sum_moduli1}
\end{figure}

If $b\geq c$, then $a\leq c\leq b\leq d$ or $a\leq c\leq d\leq b$, see the first two pictures of Fig.~\ref{fig:max_sum_moduli1}.
In the first case $S(x)$ has the minimum value $\max\{b-a,d-c\}$ for any $x\in[c,b]$.
In the second case $S(x)$ has the minimum value $b-a\geq d-c$ for any $x\in[a,b]$.
In both cases the minimum value coincides with $\max\{|a-b|,|c-d|\}$.
\smallskip

If $b\leq c$, then $a\leq b\leq c\leq d$, see the last picture of Fig.~\ref{fig:max_sum_moduli1}.
Then $S(x)=\max\{2x-a-b,c+d-2x\}$ has a minimum at $x$ such that $2x-a-b=c+d-2x$, so $x=\frac{1}{4}(a+b+c+d)$.
The minimum value is $\MS(a,b,c,d)=\frac{1}{2}|a+b-c-d|\geq\max\{|a-b|,|c-d|\}$ in this case.
The last inequality is reversed in the previous two cases.
So $\MS(a,b,c,d)=\max\{|a-b|,|c-d|,\frac{1}{2}|a+b-c-d|\}$ in all cases.
\ws
\end{proof}
 
If lattices $\La_1,\La_2$ have non-opposite signs, so $\sign(\La_1)\sign(\La_2)\geq 0$, then the metrics $\RM_q^o$ and $\PM_q^o$ from Definition~\ref{dfn:RMo} coincide with the easily computable unoriented metrics $\RM_q,\PM_q$ from Definition~\ref{dfn:RM}.
Hence Propositions~\ref{prop:RMq} and~\ref{prop:PMq} compute $\RM_q^o(\La_1,\La_2)$ and $\PM_q^o(\La_1,\La_2)$ only for lattices of opposite signs.
 
\begin{prop}[root metrics for $q=2,+\infty$]
\label{prop:RMq}
Let $\La_1,\La_2\subset\R^2$ be lattices of opposite signs with $\RI(\La_1)=(r_{12},r_{01},r_{02})$ and 
$\RI(\La_2)=(s_{12},s_{01},s_{02})$.
Then
\smallskip

\noindent
\textbf{(a)}
$\RM_2^o(\La_1,\La_2)$ is the minimum of the Euclidean distances from the point $\RI(\La_1)$ to the three points $(-s_{12},s_{01},s_{02})$, $(s_{01},s_{12},s_{02})$, and $(s_{12},s_{02},s_{01})$ in $\R^3$.  
\smallskip

\noindent
\textbf{(b)}
$\RM_{\infty}^o(\La_1,\La_2)=\min\{ d_0, d_1, d_2 \}$, where \\  
$d_0=\max\{r_{12}+s_{12}, |r_{01}-s_{01}|, |r_{02}-s_{02}|\},$\\
$d_1=\max\{\MS(r_{12},r_{01},s_{12},s_{01}),|r_{02}-s_{02}|\}$,\\
$d_2=\max\{|r_{12}-s_{12}|,\MS(r_{01},r_{02},s_{01},s_{02})\}$, \\
see $\MS(a,b,c,d)=\max\{|a-b|,|c-d|,\frac{1}{2}|a+b-c-d|\}$ in Lemma~\ref{lem:max_sum_moduli}.
\bt
\end{prop}
\begin{proof}
\textbf{(a)}
By Definition~\ref{dfn:RMo} $\RM^o(\La_1,\La_2)$ is the minimum value of $\RM(\La_1,\La_3)+\RM(\La_2,\La_3)$ over mirror-symmetric lattices $\La_3$.
By Lemma~\ref{lem:sign0} the root invariant $\RI(\La_3)$ belongs to one of the boundary sectors of the triangular cone $\TC$.
Let $\RI'$ be the mirror image of $\RI(\La_2)$ in the boundary sector containing $\RI(\La_3)$.
\smallskip

The triangle inequality for the Euclidean distance with $q=2$ implies that
$$\RM(\La_1,\La_3)+\RM(\La_2,\La_3)=||\RI(\La_1)-\RI(\La_3)||_2+||\RI(\La_3)-\RI'||_2$$ achieves minimum value $||\RI(\La_1)-\RI'||_2$ when the point $\RI(\La_3)$ is in the straight line between the points $\RI(\La_1),\RI'$.
The mirror images $\RI'$ of $\RI(\La_2)=(s_{12},s_{01},s_{02})$ in the three boundary sectors $\{s_{12}=0\}$, $\{s_{12}=s_{01}\}$, $\{s_{01}=s_{02}\}$ of the cone $\TC$ are the points $(-s_{12},s_{01},s_{02})$, $(s_{01},s_{12},s_{02})$, $(s_{12},s_{02},s_{01})$, respectively.
So $\RM^o(\La_1,\La_2)$ is the minimum of the Euclidean distances to the points above.
\medskip

\noindent
\textbf{(b)}
For $\RI(\La_1)=(r_{12},r_{01},r_{02})$ and 
$\RI(\La_2)=(s_{12},s_{01},s_{02})$, the required formula $\RM_{\infty}^o(\La_1,\La_2)=\min\{d_0,d_1,d_2\}$ will be proved by minimising the total length 
$D=\RM_{\infty}(\La_1,\La_3)+\RM_{\infty}(\La_2,\La_3)$ of a path from $\La_1$ to $\La_2$ via $\La_3$ whose root invariant $\RI(\La_3)$ can be in one of the three boundary sectors of $\TC$.
\smallskip

\noindent
\emph{Horizontal boundary} : $\RI(\La_3)=(0,t_{01},t_{02})$ for variables $0\leq t_{01}\leq t_{02}$.
Then $$D=\max\{|r_{12}+s_{12}|,|r_{01}-t_{01}|+|t_{01}-s_{01}|,|r_{02}-t_{02}|+|t_{02}-s_{02}|\}$$ has the minimum 
$d_0=\max\{r_{12}+s_{12}, |r_{01}-s_{01}|, |r_{02}-s_{02}|\}$ for  $t_{01}=\frac{1}{2}(r_{01}+s_{01})$,
$t_{02}=\frac{1}{2}(r_{02}+s_{02})$ or any values of $t_{01},t_{02}$ close enough to these averages.
\smallskip

\noindent
\emph{Inclined boundary} : $\RI(\La_3)$ consists of variables  $t_{12}=t_{01}\leq t_{02}$.
By Lemma~\ref{lem:max_sum_moduli}
$$D=\max\{|r_{12}-t_{12}|+|t_{12}-s_{12}|,|r_{01}-t_{12}|+|t_{12}-s_{01}|,|r_{02}-t_{02}|+|t_{02}-s_{02}|\}$$ has the minimum value $d_1=\max\{\MS(r_{12},r_{01},s_{12},s_{01}),|r_{02}-s_{02}|\}$, where $t_{02}$ can be anywhere between $r_{02},s_{02}$, see the formula of $\MS$ in Lemma~\ref{lem:max_sum_moduli}. 
\smallskip

\noindent
\emph{Vertical boundary} : $\RI(\La_3)$ consists of variables $t_{12}\leq t_{01}=t_{02}$.
By Lemma~\ref{lem:max_sum_moduli}
$$D=\max\{|r_{12}-t_{12}|+|t_{12}-s_{12}|,|r_{01}-t_{01}|+|t_{01}-s_{01}|,|r_{02}-t_{01}|+|t_{01}-s_{02}|\}$$ has the minimum value $d_2=\max\{|r_{12}-s_{12}|,\MS(r_{01},r_{02},s_{01},s_{02})\}$, where $t_{12}$ can be anywhere between $r_{12},s_{12}$.
The final distance is $D=\min\{d_0,d_1,d_2\}$. 
\ws
\end{proof}

\begin{prop}[projected metrics for $q=2,+\infty$]
\label{prop:PMq}
Let $\La_1,\La_2$ be lattices with opposite signs and projected invariants $\PI(\La_1)=(x_1,y_1)$, $\PI(\La_2)=(x_2,y_2)$.
\smallskip

\noindent
\textbf{(a)}
$\PM_2^o(\La_1,\La_2)$ is the minimum of the Euclidean distances from $\PI(\La_1)=(x_1,y_1)$ to the three points $(-x_2,y_2),(x_2,-y_2),(1-y_2,1-x_2)$ in $\R^2$.  
\smallskip

\noindent
\textbf{(b)}
For $x_1\leq x_2$, $\PM_{\infty}^o(\La_1,\La_2)=\min\{ d_x, d_y, d_{xy} \}$ for $d_x=\max\{x_2-x_1,y_2+y_1\}$, $d_y=\max\{x_2+x_1,|y_2-y_1|\}$, $d_{xy}=\max\{x_2-x_1,1-x_2-y_2+|1-y_1-x_2|\}$.
\bt
\end{prop}
\begin{proof}
\textbf{(a)}
By Definition~\ref{dfn:RMo} $\PM^o(\La_1,\La_2)$ is the minimum value of $\PM(\La_1,\La_3)+\PM(\La_2,\La_3)$ achieved for a mirror-symmetric lattices $\La_3$.
By Lemma~\ref{lem:sign0} the invariant $\PI(\La_3)$ belongs to one of the sides of the quotient triangle $\QT$.
Let $\PI'$ be the mirror image of $\PI(\La_2)$ with respect to the side of $\QT$ containing $\PI(\La_3)$.
\smallskip

The triangle inequality for the Euclidean distance with $q=2$ implies that 
$$\PM(\La_1,\La_3)+\PM(\La_2,\La_3)=||\PI(\La_1)-\PI(\La_3)||_2+||\PI(\La_3)-\PI'||_2$$ achieves minimum value $||\PI(\La_1)-\PI'||_2$ when the point $\PI(\La_3)$ is in the straight line between the points $\PI(\La_1),\PI'$ in the plane $\R^2$ containing the quotient triangle $\QT$.
The mirror images $\PI'$ of $\PI(\La_2)=(x_2,y_2)$ in the three sides $x_2=0$, $y_2=0$, $x_2+y_2=1$ are the points $(-x_2,y_2)$, $(x_2,-y_2)$, $(1-y_2,1-x_2)$, respectively.
Hence $\PM^o(\La_1,\La_2)$ is the minimum of the Euclidean distances to the points above.
\medskip

\noindent
\textbf{(b)}
$\PM^o(\La_1,\La_2)$ is the minimum value of $\PM(\La_1,\La_3)+\PM(\La_2,\La_3)$ over mirror-symmetric lattices $\La_3$.
The formula $\PM_{\infty}^o(\La_1,\La_2)=\min\{d_x,d_y,d_{xy}\}$ will be proved by minimising the Minkowski length $M_{\infty}$ of a path from $\La_1$ to $\La_2$ via $\La_3$ whose projected invariant $\PI(\La_3)$ can be in one of the three sides of the quotient triangle $\QT$.
For given $\PI(\La_1)=(x_1,y_1)$ and 
$\RI(\La_2)=(x_2,y_2)$ with $x_1\leq x_2$, we minimise $D=\PM_{\infty}(\La_1,\La_3)+\PM_{\infty}(\La_2,\La_3)$ for each side of $\QT$ below.
\smallskip

\noindent
\emph{Horizontal side} : $\PI(\La_3)=(x_3,0)$ for a variable parameter $x_3\geq 0$.
For any $x_3\in[x_1,x_2]$, the distance $D=\max\{(x_3-x_1)+(x_2-x_3),y_2+y_1\}$ equals the simpler function $d_x=\max\{x_2-x_1,y_2+y_1\}$ and can be only larger for any $x_3\not\in[x_1,x_2]$.
\smallskip

\noindent
\emph{Vertical side} : $\PI(\La_3)=(0,y_3)$ for a variable parameter $y_3$.
For any $y_3$ between $y_1,y_2$, the distance $D=\max\{x_1+x_2,|y_1-y_3|+|y_3-y_2|\}$ has the minimum $d_y=\max\{x_1+x_2,|y_1-y_2|\}$ and is larger for any $y_3$ that is not between $y_1,y_2$.
\smallskip

\begin{figure}[h]
\includegraphics[width=1.0\textwidth]{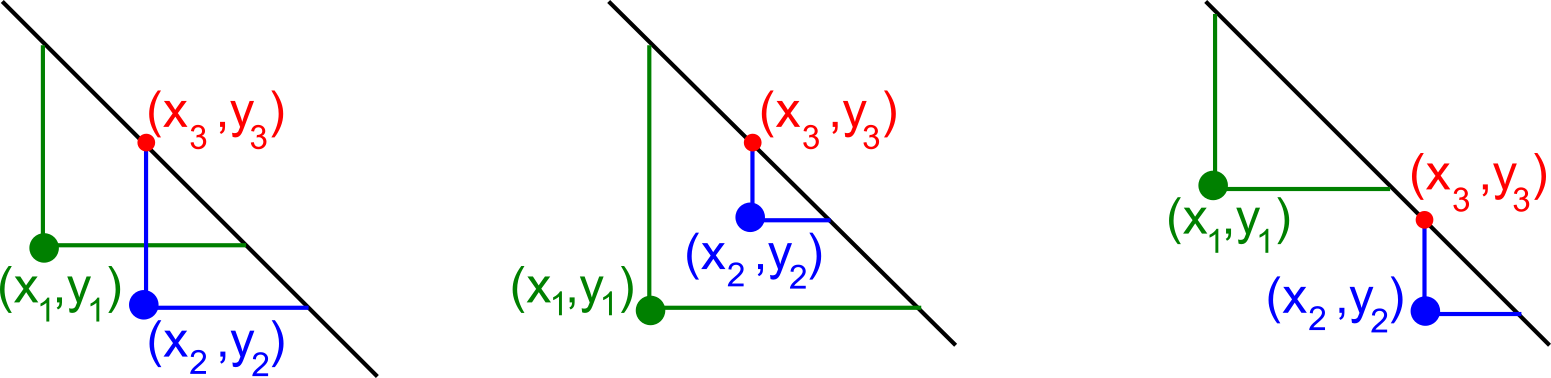}
\caption{Relative positions of $(x_1,y_1),(x_2,y_2)$ for $x_1\leq x_2$ in the proof of Proposition~\ref{prop:PMq}(b).}
\label{fig:max_sum_moduli2}
\end{figure}

\noindent
\emph{Hypotenuse} : $\PI(\La_3)=(x_3,y_3)$ for variables $x_3,y_3$ such that $x_3+y_3=1$.
For $x_1\leq x_2$, we aim to minimise $D=\max\{|x_1-x_3|+|x_3-x_2|,|y_1-y_3|+|y_3-y_2|\}$.
Fig.~\ref{fig:max_sum_moduli2} shows all three different cases how we can find an optimal chiral lattice $\La_3$ with a projected invariant $\PI(\La_3)=(x_3,y_3)$ for given $\PI(\La_1)=(x_1,y_1)$ and $\PI(\La_2)=(x_2,y_2)$.
The sum $|x_1-x_3|+|x_3-x_2|$ has the minimum value $x_2-x_1$ for any $x_3\in[x_1,x_2]$. 
Similarly, $|y_1-y_3|+|y_3-y_2|$ has the minimum value $|y_1-y_2|$ for any $y_3$ between $y_1,y_2$. 
If the point $(x_3,y_3)$ moves along the hypotenuse so that both $x_3,y_3$ are outside their minimum ranges above, then $|x_1-x_3|+|x_3-x_2|=|2x_3-x_1-x_2|$ increases at the same rate as $|y_1-y_3|+|y_3-y_2|=|2y_3-y_1-y_2|$ decreases during this movement because $x_3+y_3=1$.
Hence $x_3$ can be chosen as $x_1$ or $x_2$ to minimise $D$.
In the first two pictures of Fig.~\ref{fig:max_sum_moduli2} we choose $x_3=x_2$.
\smallskip

In the last picture of Fig.~\ref{fig:max_sum_moduli2}, any $(x_3,y_3)$  between the triangles with right-angled vertices at $(x_1,y_1),(x_2,y_2)$ gives the minimum values $x_2-x_1$ and $|y_2-y_1|$.
Hence $x_3=x_2$ always gives the minimum value of the distance
$$d_{xy}=\max\{x_2-x_1,|y_1-y_3|+|y_3-y_2|\} \text{ for } y_3=1-x_2,$$ 
so $d_{xy}=\max\{x_2-x_1,1-x_2-y_2+|1-y_1-x_2|\}$ and $D=\min\{d_x,d_y,d_{xy}\}$.
\ws
\end{proof}

\section{Real-valued chiral distances measure asymmetry of lattices}
\label{sec:chiral_distances}

The classical concept of chirality is a binary property distinguishing mirror images of the same object such as a molecule or a periodic crystal.
Continuous classifications in Theorem~\ref{thm:classification2d} and Corollary~\ref{cor:similar_lattices2d} imply that the binary chirality is discontinuous under almost any perturbations similar to other discrete invariants such as symmetry groups.
To avoid arbitrary thresholds, it makes more sense to continuously quantify a deviation of a lattice from a higher-symmetry neighbour.  
\medskip

The term \emph{chirality} often refers to 3-dimensional molecules or crystal lattices.
One reason is the fact that in $\R^2$ a reflection with respect to a line $L$ is realised by the rotation in $\R^3$ around $L$ through $180^\circ$. 
However, if our ambient space is only $\R^2$, the concepts of isometry and rigid motion differ.
For example, Lemma~\ref{lem:achiral_lattices} described root invariants of all lattices that are related to their mirror images by rigid motion.
Such lattices can be called \emph{achiral}. 
We call them mirror-symmetric. 
\medskip

After consulting with crystallographers, Definition~\ref{dfn:RC} introduces the real-valued $G$-chiral distances of a lattice $\La\subset\R^2$.
Proposition~\ref{prop:cont_chiral_distances} will prove continuity of these functions $\RC[G]:\LIS(\R^2)\to\R$ and $\PC[G]:\LSS(\R^2)\to\R$.
\medskip

Recall that the \emph{crystallographic point group} $G$ of a lattice $\La\subset\R^2$ containing the origin $0$ consists of all symmetry operations that keep $0$ and map $\La$ to itself.
For example, any such group $G$ includes the central symmetry with respect to $0\in\La\subset\R^2$.
If $G$ has no other non-trivial symmetries, we get $G=C_2$ in Schonflies notations.
All 2D lattices split into four crystal families by their point groups: oblique ($C_2$), orthorhombic ($D_2$), tetragonal or square ($D_4$) and hexagonal ($D_6$).
Orthorhombic lattices split into rectangular and centred rectangular, see Fig.~\ref{fig:QT+QS}.

\begin{figure}[h]
\includegraphics[width=1.0\textwidth]{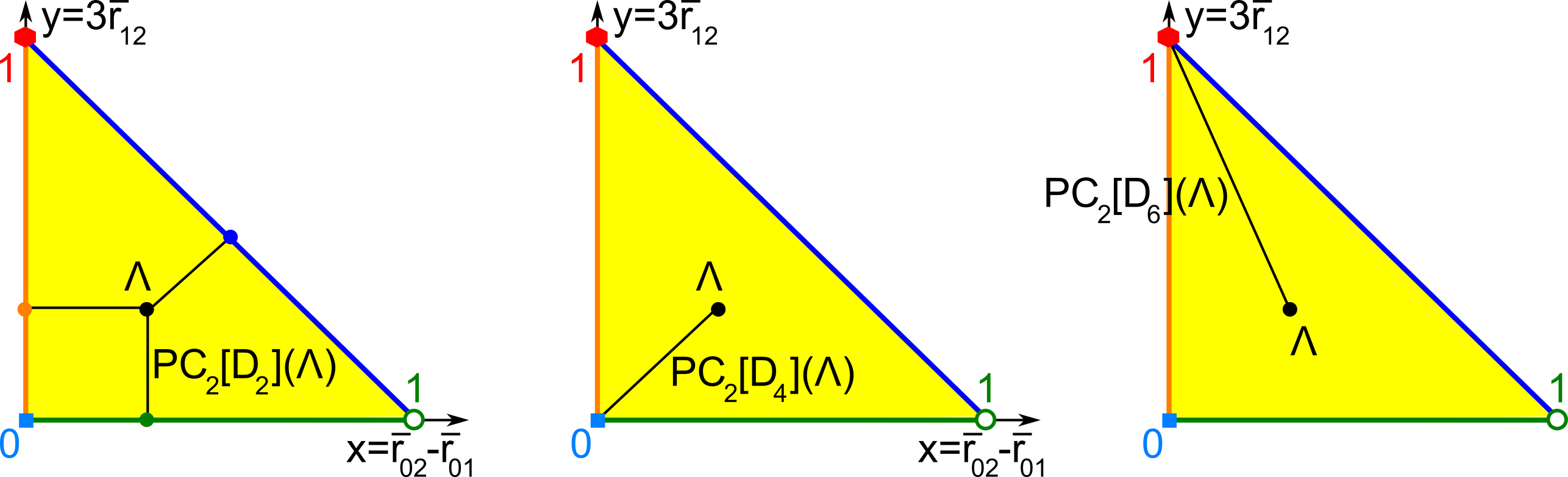}
\caption{\textbf{Left}: by Definition~\ref{dfn:RC}, the projected $D_2$ chiral distance $\PC_2[D_2](\La)$ is the minimum Euclidean distance from $\PI(\La)\in\QT$ to the boundary $\bd\QT$.
\textbf{Middle}: $\PC_2[D_4](\La)$ is the distance from $\PI(\La)$ to $(0,0)$.
\textbf{Right}: $\PC_2[D_6](\La)$ is the distance from $\PI(\La)$ to $(0,1)$.
}
\label{fig:QT+PC}
\end{figure}
 
\begin{dfn}[$G$-chiral distances {$\RC[G]$ and $\PC[G]$}]
\label{dfn:RC}
For any crystallographic point group $G$ in $\R^2$, let $\LIS[G]\subset\LIS(\R^2)$ be the closure of the subspace of all (isometry classes of) lattices that have the crystallographic  point group $G$.
For $G=D_2$ or $G=D_4$ or $G=D_6$, the root and projected \emph{$G$-chiral} distances are 
$$\RC[G](\La)=\min\limits_{\La'\in\LIS[G]}\RM(\La,\La')\geq 0 \text{ and }
\PC[G](\La)=\min\limits_{\La'\in\LIS[G]}\PM(\La,\La')\geq 0,$$
where $\RM$, $\PM$ are any metrics from Definition~\ref{dfn:RM} with a base metric $d$.
If $d=M_q$ for $q\in[1,+\infty]$, denote the $G$-chiral distances by $\RC_q[G]$ and $\PC_q[G]$.
\bs
\end{dfn}

Since any lattice $\La$ is symmetric with respect to the origin $0\in\La$,
the closed subspace $\LIS[C_2]$ coincides with the 3-dimensional Lattice Isometry Space $\LIS(\R^2)$.
The 2-dimensional subspace $\LIS[D_2]$ consists of all mirror-symmetric lattices (rectangular and centred-rectangular) represented by root invariants $\RI$ on the boundary $\bd\TC$ of the triangular cone in Definition~\ref{dfn:TC}, see  Fig.~\ref{fig:TC}.
The 1-dimensional subspaces $\LIS[D_4],\LIS[D_6]\subset\LIS[D_2]$ can be viewed as the blue and orange rays $\{r_{12}=0<r_{01}=r_{02}\}$ and  $\{0<r_{12}=r_{01}=r_{02}\}$, respectively.   
\medskip

The $G$-chiral distance $\RC[G]$ in Definition~\ref{dfn:RC} measures a distance from $\RI(\La)$ to the root invariant of a closest neighbour in the subspace $\LIS[G]$.
Any $\RC[G](\La)$ is invariant up to isometry and measures a distance from $\La$ to its nearest neighbour $\La'\in\LIS[G]$.
The \emph{signed chiral distances} $\sign(\La)\RC(\La)$ and $\sign(\La)\PC(\La)$ are invariant up to rigid motion.
Since $\LIS[G]$ is a closed subspace within $\LIS(\R^2)$, the continuous distances $\RM,\PM$ achieve their minima if their base distances $d$ are continuous.
If $\LIS[D_2]$ was defined as an open subspace of only lattices that have the point group $D_2$ (not $D_4$ or $D_6$), then $\RC[G],\PC[G]$ should be defined via infima instead of simpler minima.
Indeed, any square or hexagonal lattice $\La$ can be approximated by infinitely many closer and closer orthorhombic lattices $\La'$, but the expected distance $\RM(\La,\La')=0$ will not be achieved on an open set.
\medskip

For 
$q=2,+\infty$, the distances $\RC_q,\PC_q$ are computed in Propositions~\ref{prop:RC},~\ref{prop:PC}.

\begin{lem}[properties of chiral distances]
\label{lem:RC}
\textbf{(a)}
A lattice $\La$ is mirror-symmetric if and only if $\RC[D_2](\La)=0$ or, equivalently, $\PC[D_2](\La)=0$.
\medskip

\noindent
\textbf{(b)}
For any crystallographic point group $G$ in $\R^2$, mirror reflections $\La^{\pm}\subset\R^2$ have equal $G$-chiral distances: $\RC[G](\La^+)=\RC[G](\La^-)$,
$\PC[G](\La^+)=\PC[G](\La^-)$.
\bt
\end{lem}
\begin{proof}
\textbf{(a)}
By Definition~\ref{dfn:RC} $\RC[D_2](\La)=\min\limits_{\sign(\La')=0}\{\RM(\La,\La')\}=0$ means that $\RI(\La)=\RI(\La')$ for some mirror-symmetric lattice $\La'$ due to Lemma~\ref{lem:sign0}.
Then $\La$ is isometric to $\La'$ by Theorem~\ref{thm:classification2d} and is mirror-symmetric.
\medskip

\noindent
\textbf{(b)}
The $G$-chiral distances from Definition~\ref{dfn:RC} depend only on a root invariant $\RI(\La)$ and projected invariant $\PI(\La)$, which are the same for mirror images $\La^{\pm}$.
\ws
\end{proof}

\begin{lem}[lower bounds]
\label{lem:RCbounds}
\textbf{(a)}
If lattices $\La_1,\La_2$ have opposite signs, then
$\RM^o(\La_1,\La_2)\geq
\RC[D_2](\La_1)+\RC[D_2](\La_2)$ and \\
$\PM^o(\La_1,\La_2)\geq
\PC[D_2](\La_1)+\PC[D_2](\La_2)$.
\medskip

\noindent
\textbf{(b)}
For the mirror images $\La^{\pm}$ of any lattice $\La$, the lower bounds in part (a) become equalities:
$\RM^o(\La^+,\La^-)=2\RC[D_2](\La)$ and
$\PM^o(\La^+,\La^-)=2\PC[D_2](\La)$.
\bt
\end{lem}
\begin{proof}
\textbf{(a)}
The first lower bound follows from Definitions~\ref{dfn:RC} and~\ref{dfn:RMo}, because the minimisations of the root chiral distances are over separate lattices $\La'_1,\La'_2$.
$$\RM^o(\La_1,\La_2)=\min\limits_{\sign(\La')=0}(\RM(\La_1,\La')+\RM(\La_2,\La'))\geq$$
$$\min\limits_{\sign(\La'_1)=0}\RM(\La_1,\La'_1)+
\min\limits_{\sign(\La'_2)=0}\RM(\La_2,\La'_2)
=\RC[D_2](\La_1)+\RC[D_2](\La_2),$$
where the root chiralities are independent of a sign by Lemma~\ref{lem:RC}(b).
The proof for the projected metric $\PM^o$ is similar to the above arguments. 
\medskip

\noindent
\textbf{(b)}
If $\La_1=\La_2$, the sum $\RM(\La_1,\La')+\RM(\La_2,\La')$ in part \textbf{(a)} consists of equal terms.
Hence the inequality becomes equality giving the double distance. 
\ws
\end{proof}

Lemma~\ref{lem:max_modulus} will help prove
Propositions~\ref{prop:RC},~\ref{prop:PC} to explicitly express the chiral distances $\RC_q[G](\La),\PC_q[G](\La)$ via the forms $\RI(\La)\in\TC$ and $\PI(\La)\in\QT$. 

\begin{lem}[maximum modulus]
\label{lem:max_modulus}
For any fixed points $a,b\in\R$, the function $D_{\infty}(x)=\max\{|a-x|,|x-b|\}$ has the minimum value $\dfrac{1}{2}|a-b|$ for $x=\dfrac{a+b}{2}$.
\bt
\end{lem}
\begin{proof}
Without loss of generality assume that $a\leq b$.
For $x\in[a,b]$, the function $D_{\infty}(x)=\max\{x-a,b-x\}$ has the minimum value $\dfrac{b-a}{2}$ at the mid-point $x$ of the interval $[a,b]$ and takes values larger than $b-a$ for any $x\not\in[a,b]$.
\ws
\end{proof}

\begin{prop}[chiral distances {$\RC_q[G]$} for $q=2,+\infty$]
\label{prop:RC}
Let a lattice $\La\subset\R^2$ have a root invariant $\RI(\La)=(r_{12},r_{01},r_{02})$ with $0\leq r_{12}\leq r_{01}\leq r_{02}$.
Then
\medskip

\hspace*{-7mm}
$\begin{array}{llr}
(\ref{prop:RC}a) &
\RC_2[D_2](\La)=\min\left\{ r_{12}, \dfrac{r_{01}-r_{12}}{\sqrt{2}}, \dfrac{r_{02}-r_{01}}{\sqrt{2}}\right\}; & \\
\smallskip

&
\RC_2[D_4](\La)=\sqrt{r_{12}^2+\frac{1}{4}(r_{02}-r_{01})^2}; & \\
\smallskip

&
\RC_2[D_6](\La)=\sqrt{\frac{2}{3}(r_{12}^2+r_{01}^2+r_{02}^2-r_{12}r_{01}-r_{12}r_{02}-r_{01}r_{02})}; & \\

(\ref{prop:RC}b) &
\RC_{\infty}[D_2](\La)=\min\left\{ r_{12}, \dfrac{r_{01}-r_{12}}{2}, \dfrac{r_{02}-r_{01}}{2}\right\}. & \\

&
\RC_{\infty}[D_4](\La)=\min\{r_{12},\dfrac{r_{02}-r_{01}}{2}\}; & \\
\smallskip

&
\RC_{\infty}[D_6](\La)=\dfrac{r_{02}-r_{12}}{2}. & \blacktriangle
\end{array}$
\end{prop}
\begin{proof}
\textbf{(a)}
The chiral distance $\RC_2[D_2](\La)=\min\limits_{\sign(\La')=0}||\RI(\La)-\RI(\La')||_2$ by Definition~\ref{dfn:RC} is the minimum Euclidean distance from $\RI(\La)$ to the boundary of $\TC$.
This boundary of three triangular sectors consists of root invariants of all mirror-symmetric lattices by Lemma~\ref{lem:sign0}.
Any point $\RI(\La)=(r_{12},r_{01},r_{02})$ with $0\leq r_{12}\leq r_{01}\leq r_{02}$ has Euclidean distances $r_{12}$, $\dfrac{r_{01}-r_{12}}{\sqrt{2}}$, $\dfrac{r_{02}-r_{01}}{\sqrt{2}}$ to the horizontal boundary $r_{12}=0$, inclined boundary $r_{12}=r_{01}$, and the vertical boundary $r_{01}=r_{02}$, respectively, see Fig.~\ref{fig:QS+DC}(right).
Then $\RC_2[D_2](\La)$ is the minimum of the above Euclidean distances to the three boundary sectors of $\TC$.
\medskip

\noindent
$\RC_2[D_4](\La)$ is the Euclidean distance from $\RI(\La)=(r_{12},r_{01},r_{02})$ to a closest root invariant $\RI(\La')=(0,s,s)$ of a square lattice $\La'$.
The square of 
$||\RI(\La)-\RI(\La')||_2$ is $d_4(s)=r_{12}^2+(s-r_{01})^2+(s-r_{02})^2$.
The quadratic function $d_4(s)$ is minimised when $0=d'_4(s)=2(s-r_{01})+2(s-r_{02})$, so $s=\frac{1}{2}(r_{01}+r_{02})$.
Then $d_4(\frac{1}{2}(r_{01}+r_{02}))=r_{12}^2+\frac{1}{4}(r_{02}-r_{01})^2$ and
$\RC_2[D_4](\La)=\sqrt{r_{12}^2+\frac{1}{4}(r_{02}-r_{01})^2}$.
\medskip

\noindent
$\RC_2[D_6](\La)$ is the Euclidean distance from $\RI(\La)=(r_{12},r_{01},r_{02})$ to a closest root invariant $\RI(\La')=(s,s,s)$ of a hexagonal lattice $\La'$.
The square of 
$||\RI(\La)-\RI(\La')||_2$ is $d_6(s)=(s-r_{12})^2+(s-r_{01})^2+(s-r_{02})^2$.
The quadratic function $d_6(s)$ is minimised when $0=d'_6(s)=2(s-r_{12})+2(s-r_{01})+2(s-r_{02})$, so $s=\frac{1}{3}(r_{12}+r_{01}+r_{02})$.
Substituting the minimum point $s$ above, we get
$$9d_6(s)=(2r_{12}-r_{01}-r_{02})^2+(2r_{01}-r_{12}-r_{02})^2
+(2r_{02}-r_{12}-r_{01})^2=$$
$$=6(r_{12}^2+r_{01}^2+r_{02}^2)-(4+4-2)(r_{12}r_{01}+r_{12}r_{02}+r_{01}r_{02}).$$
Then $\RC_2[D_6](\La)=\sqrt{d_6(s)}=\sqrt{\frac{2}{3}(r_{12}^2+r_{01}^2+r_{02}^2-r_{12}r_{01}-r_{12}r_{02}-r_{01}r_{02})}$.
\medskip

\noindent
\textbf{(b)}
$\RC_{\infty}[D_2](\La)=\min\limits_{\sign(\La')=0}||\RI(\La)-\RI(\La')||_{\infty}$ is minimised over mirror-symmetric lattices $\La'$ whose root invariants by Lemma~\ref{lem:sign0} belong to one of the three boundary sectors of the triangular cone $\TC$.
We consider them below one by one.
\medskip

\noindent
\emph{Horizontal boundary} : $\RI(\La')=(0,s_{01},s_{02})$ for 
$0<s_{01}\leq s_{02}$.
Then
$||\RI(\La)-\RI(\La')||_{\infty}=||(r_{12},r_{01},r_{02})-(0,s_{01},s_{02})||_{\infty}=\max\{r_{12},|r_{01}-s_{01}|,|r_{02}-s_{02}|\}$ has the minimum $r_{12}$ for $s_{01}=r_{01}$, $s_{02}=r_{02}$ or any $s_{01},s_{02}$ close to $r_{01},r_{02}$.
\medskip

\noindent
\emph{Inclined boundary} : a variable root invariant $\RI(\La')$ has equal root products $s_{12}=s_{01}$.
By Lemma~\ref{lem:max_modulus} the $M_{\infty}$-distance 
$||\RI(\La)-\RI(\La')||_{\infty}=||(r_{12},r_{01},r_{02})-(s_{01},s_{01},s_{02})||_{\infty}=\max\{|r_{12}-s_{01}|,|r_{01}-s_{01}|,|r_{02}-s_{02}|\}$ has the minimum $\dfrac{r_{01}-r_{12}}{2}$ for $s_{02}=r_{02}$ and $s_{01}$ at the mid-point of the interval $[r_{12},r_{01}]$.
\medskip

\noindent
\emph{Vertical boundary} : $\RI(\La')$ has $s_{01}=s_{02}$.
By Lemma~\ref{lem:max_modulus}
$||\RI(\La)-\RI(\La')||_{\infty}=||(r_{12},r_{01},r_{02})-(s_{12},s_{01},s_{01})||_{\infty}=\max\{|r_{12}-s_{12}|,|r_{01}-s_{01}|,|r_{02}-s_{01}|\}$ has the minimum $\dfrac{r_{02}-r_{01}}{2}$ for $s_{12}=r_{12}$ and $s_{01}$ at the mid-point of the interval $[r_{01},r_{02}]$.
Finally, $\RC_{\infty}[D_2](\La)$ is the minimum of the three $L_{\infty}$ distances.
\medskip

\noindent
$\RC_{\infty}[D_4](\La)$ is the distance $M_{\infty}$ from $\RI(\La)$ to a closest root invariant $\RI(\La')=(0,s,s)$ of a square lattice $\La'$.
Then
$||\RI(\La)-\RI(\La')||_{\infty}=\max\{r_{12},|s-r_{01}|,|s-r_{02}|\}$.
If we ignore $r_{12}$, by Lemma~\ref{lem:max_modulus} the minimum of the largest value among the last two is $\frac{1}{2}(r_{02}-r_{01})$, so $\RC_{\infty}[D_4](\La)=\min\left\{ r_{12}, \frac{r_{02}-r_{01}}{2}\right\}$.
\medskip

\noindent
$\RC_{\infty}[D_6](\La)$ is the distance $M_{\infty}$ from $\RI(\La)=(r_{12},r_{01},r_{02})$ to a closest root invariant $\RI(\La')=(s,s,s)$ of a hexagonal lattice $\La'$.
Then
$||\RI(\La)-\RI(\La')||_{\infty}=\max\{|s-r_{12}|,|s-r_{01}|,|s-r_{02}|\}$.
Since $r_{12}\leq r_{01}\leq r_{02}$, we could ignore  $|s-r_{01}|$ in the maximum.
By Lemma~\ref{lem:max_modulus} the final maximum is
$\frac{r_{02}-r_{01}}{2}=\RC_{\infty}[D_6](\La)$.
\ws
\end{proof}

When considering lattices up to similarity, the subspace $\LSS[D_4]$ consists of a single class of all square lattices, which are all equivalent up to isometry and uniform scaling.
The subspace $\LSS[D_6]$ is also a single point representing all hexagonal lattices.
Then $\PC[D_4]$ and $\PC[D_6]$ are distances to these single points. 

\begin{prop}[chiral distances $\PC_q$ for $q=2,+\infty$]
\label{prop:PC}
Let a lattice $\La$ have a projected invariant $\PI(\La)=(x,y)\in\QT$ so that $x\in[0,1)$, $y\in[0,1]$, $x+y\leq 1$.
\medskip

\hspace*{-7mm}
$\begin{array}{llr}
(\ref{prop:PC}a) &
\PC_2[D_2](\La)=\min\left\{x,y,\dfrac{1-x-y}{\sqrt{2}}\right\}, & \\
\smallskip

&
\PC_q[D_4](\La)=(x^q+y^q)^{1/q} \text{ for any } q\in[1,+\infty), & \\
\smallskip

&
\PC_q[D_6](\La)=(x^q+(1-y)^q)^{1/q} \text{ for any } q\in[1,+\infty); & \\
\smallskip

(\ref{prop:PC}b) &
\PC_{\infty}[D_2](\La)=\min\left\{x,y,\dfrac{1-x-y}{2}\right\}, & \\
\smallskip

&
\PC_{\infty}[D_4](\La)=x, & \\
\smallskip

&
\PC_{\infty}[D_6](\La)=1-y. & 
\end{array}$
\medskip

\noindent
\emph{(\ref{prop:PC}c)}
The upper bounds 
$\PC_{2}[D_2](\La)\leq\frac{1}{2+\sqrt{2}}$,
$\PC_{\infty}[D_2](\La)\leq\frac{1}{4}$ hold for any $\La$, achieved for lattices with $\PI(\La_2)=(\frac{1}{2+\sqrt{2}},\frac{1}{2+\sqrt{2}})$, $\PI(\La_{\infty})=(\frac{1}{4},\frac{1}{4})$, respectively.
For $q\in[1,+\infty]$, the bound $\PC_{q}[D_4](\La)\leq 1$ holds for any $\La$ and is achieved for any hexagonal lattice.
For $q\in[1,+\infty)$, the upper bound $\PC_{q}[D_6](\La)<2^{1/q}$ holds for any $\La$ and is approached but not achieved as $x\to 1$.
The bound $\PC_{\infty}[D_6](\La)\leq 1$ holds for any $\La$ and is achieved for any square and rectangular lattice. 
\bt
\end{prop}
\begin{proof}
\textbf{(a)}
By Definition~\ref{dfn:RC} $\PC_2[D_2](\La)$ is the minimum Euclidean distance from $\PI(\La)$ to the three boundary sides of the quotient triangle $\QT$.
These sides consists of projected invariants of all mirror-symmetric lattices by Lemma~\ref{lem:sign0}.
Any point $\PI(\La)=(x,y)\in\QT$ has distances $x$, $y$, $\dfrac{1-x-y}{\sqrt{2}}$ to the vertical side $x=0$, horizontal side $y=0$, and the hypotenuse $x+y=1$, respectively, see Fig.~\ref{fig:QS+DC}~(left).
Hence $\PC_2(\La)$ is the minimum of the above Euclidean distances.  
\bigskip

\noindent
$\PC_q[D_4](\La)=(x^q+y^q)^{1/q}$ is the Minkowski $M_q$ distance from $\PI(\La)=(x,y)\in\QT$ to the single-point subspace $\LSS[D_4]=(0,0)$ for any $q\in[1,+\infty)$.
\medskip

\noindent
$\PC_q[D_6](\La)=(x^q+(1-y)^q)^{1/q}$
 is the Minkowski $M_q$ distance from $\PI(\La)=(x,y)\in\QT$ to the single-point subspace $\LSS[D_6]=(0,1)$ for any $q\in[1,+\infty)$.
\medskip

\noindent
\textbf{(b)}
$\PC_{\infty}(\La)=\min\limits_{\sign(\La')=0}||\PI(\La)-\PI(\La')||_{\infty}$ is minimised over lattices $\La'$ whose projected invariant is in one of the three sides of the quotient triangle $\QT$.
\medskip

\noindent
\emph{Horizontal side} : $\PI(\La')=(s,0)$ for a variable $s$.
The distance $||\PI(\La)-\PI(\La')||_{\infty}=||(x,y)-(s,0)||_{\infty}=\max\{|x-s|,y\}$ has the minimum value $y$ for $s=x$.
\medskip

\noindent
\emph{Vertical side} : $\PI(\La')=(0,t)$ for a variable $t$.
The distance $||\PI(\La)-\PI(\La')||_{\infty}=||(x,y)-(0,t)||_{\infty}=\max\{x,|y-t|\}$ has the minimum $x$ for $t=y$.
\medskip

\noindent
\emph{Hypotenuse} : $\PI(\La')=(s,t)$ for variables $s,t\geq 0$ such that $s+t=1$.
To compute the $M_{\infty}$ distance from $(x,y)$ to $(s,t)$, first assume that $s\geq x$, $t\geq y$.
Then $||(x,y)-(s,t)||_{\infty}=\max\{s-x,t-y\}$ is minimised when $s-x=t-y$.
Substituting $t=1-s$, we get $s-x=1-s-y$, $s=\dfrac{1+x-y}{2}$,
$t=\dfrac{1-x+y}{2}$.
One can check that $s+t=1$ and $s\geq x$, $t\geq y$  due to $x+y\leq 1$ as $(x,y)\in\QT$.
Then $s-x=t-y=\dfrac{1-x-y}{2}$.
It remains to show that the minimum $M=\dfrac{1-x-y}{2}$ of the distance from $(x,y)$ to $(s,t)$ cannot have a smaller value for $s\leq x$ or $t\leq y$.
\medskip

If $s\leq x$, then $||(x,y)-(s,t)||_{\infty}=\max\{x-s,t-y\}=\max\{x-1+t,t-y\}\leq t-y$, whose minimum value $(1-x)-y$ is not less than $M=\dfrac{1-x-y}{2}$ as $x+y\leq 1$.  
\medskip

If $t\leq y$, then $||(x,y)-(s,t)||_{\infty}=\max\{s-x,y-t\}=\max\{s-x,y-1+s\}\leq s-x$, whose minimum value $(1-y)-x$ is not less than $M=\dfrac{1-x-y}{2}$ as $x+y\leq 1$.  
\medskip

\noindent
So $\PC_{\infty}(\La)$ is the minimum of the above three $M_{\infty}$ distances.  
\medskip

\noindent
$\PC_{\infty}[D_4](\La)=x$ is the distance  $M_{\infty}$ from $(x,y)$ to  $\LSS[D_4]=(0,0)$.
\medskip

\noindent
$\PC_{\infty}[D_6](\La)=1-y$ is the distance $M_{\infty}$ from $(x,y)$ to  $\LSS[D_6]=(0,1)$.
\medskip

\noindent
\textbf{(c)}
The quotient triangle $\QT$ in Fig.~\ref{fig:QT+QS}~(left) is parameterised by $0\leq x<1$ and $0\leq y\leq 1$ such that $x+y\leq 1$,
By Theorem~\ref{prop:PC}(a) the distance $\PC_{2}[D_2](\La)=\min\{x,y,\dfrac{1-x-y}{\sqrt{2}}\}$ is maximal when $x=y=\dfrac{1-x-y}{\sqrt{2}}$, so $x=y=\dfrac{1}{2+\sqrt{2}}$.
By Theorem~\ref{prop:PC}(b)
$\PC_{\infty}[D_2](\La)=\min\{x,y,\dfrac{1-x-y}{2}\}$ is maximal when $x=y=\dfrac{1-x-y}{2}$, $x=y=\dfrac{1}{4}$.
Then $x^q\leq x$, $y^q\leq y$ and $(x^q+y^q)^{1/q}\leq (x+y)^{1/q}\leq 1$ for any $q\in[1,+\infty)$.
Hence the upper bound $\PC_{q}[D_4](\La)\leq 1$ holds for any $q\in[1,+\infty]$ and is achieved for any hexagonal lattice with $\PI=(0,1)$.
Similarly, the upper bound $\PC_{q}[D_6](\La)<2^{1/q}$ holds for any $q\in[1,+\infty)$ and is approached but not achieved as $x\to 1$, $y=0$.
The bound $\PC_{\infty}[D_6](\La)=1-y\leq 1$ holds for any $\La$ and is achieved for any square and rectangular lattice with $y=0$. 
\ws
\end{proof}

\begin{exa}[distances $\RC_q,\PC_q$]
\label{exa:RC}
Table~\ref{tab:RC} shows the chiral distances computed by Propositions~\ref{prop:RC},~\ref{prop:PC} for the prominent lattices $L_2^{\pm}$, $L_{\infty}^{\pm}$ in
Example~\ref{exa:inverse_design}.
\bs
\end{exa}

\begin{table}[h]
\caption{Chiral distances $\PC_q,\RC_q$ for the lattices $L_2^{\pm},L_{\infty}^{\pm}$ in Fig.~\ref{fig:QS+DC} and~\ref{fig:DT}, see Example~\ref{exa:RC}. }
\medskip

\begin{tabular}{|l|cc|}  
\hline    
$\La$ & $L_{\infty}$ & $L_2$ \\

$\PI(\La)$ & 
$(\frac{1}{4},\frac{1}{4})$ & 
$(\frac{1}{2+\sqrt{2}},\frac{1}{2+\sqrt{2}})$ \\
\hline
 
$\PC_2[D_2]$ & 
$\frac{1}{4}$ & 
$\frac{1}{2+\sqrt{2}}$ \\

$\PC_2[D_4]$  & 
$\frac{\sqrt{2}}{4}$ & 
$\sqrt{2}-1$ \\

$\PC_2[D_6]$ & 
$\frac{\sqrt{10}}{4}$ & 
$\sqrt{2-\sqrt{2}}$ \\
\hline

$\PC_{\infty}[D_2]$  & 
$\frac{1}{4}$ &
$\frac{1}{2+\sqrt{2}}$ \\

$\PC_{\infty}[D_4]$ & 
$\frac{1}{4}$ &
$\frac{1}{2+\sqrt{2}}$ \\

$\PC_{\infty}[D_6]$  & 
$\frac{3}{4}$ &
$\frac{1}{\sqrt{2}}$ \\
\hline
\end{tabular}
\begin{tabular}{|l|cc|}  
\hline    
$\La$ & $L_{\infty}$ & $L_2$ \\

$\RI(\La)$ & 
$(1,4,7)$ & 
$(2-\sqrt{2},2\sqrt{2}-1,5-\sqrt{2})$ \\
\hline
 
$\RC_2[D_2]$ & 
$1$ & 
$2-\sqrt{2}$ \\

$\RC_2[D_4]$  & 
$\frac{\sqrt{13}}{2}$ & 
$(2-\sqrt{2})\frac{\sqrt{13}}{2}$ \\

$\RC_2[D_6]$ & 
$3\sqrt{2}$ & 
$\sqrt{2(13-3\sqrt{2})}$ \\
\hline

$\RC_{\infty}[D_2]$  & 
$1$ & 
$2-\sqrt{2}$ \\

$\RC_{\infty}[D_4]$ & 
$1$ & 
$2-\sqrt{2}$ \\

$\RC_{\infty}[D_6]$  & 
$3$ &
$\frac{3}{2}$ \\
\hline
\end{tabular}
\label{tab:RC}
\end{table}

\begin{exa}[metrics $\RM_q^o,\PM_q^o$]
\label{exa:RMo}
Table~\ref{tab:RMo} has $\RM_q^o,\RM_q^o$ for $q=2,+\infty$ and
 the prominent lattices $L_2^{\pm}$, $L_{\infty}^{\pm}$, which were inversely designed in 
Example~\ref{exa:inverse_design}.
\medskip

If lattices have the same sign, then $\RM^o,\PM^o$ coincide with their unoriented versions by Definition~\ref{dfn:RMo}.
For example, $\PM_q^o(L_2^+,L_{\infty}^+)$ is the distance $M_q$ between the invariants $\PI(L_{\infty})=(\frac{1}{4},\frac{1}{4})$ and $\PI(L_2)=(\frac{1}{2+\sqrt{2}},\frac{1}{2+\sqrt{2}})=(1-\frac{1}{\sqrt{2}},1-\frac{1}{\sqrt{2}})$, so
$\PM_{\infty}^o(L_2^+,L_{\infty}^+)=\frac{3}{4}-\frac{1}{\sqrt{2}}\approx 0.04$ and $\PM_{2}^o(L_2^+,L_{\infty}^+)=\frac{3}{4}\sqrt{2}-1\approx 0.06$.
\medskip

Similarly, 
$\RM_q^o(L_2^+,L_{\infty}^+)$ is the $M_q$ distance between the root invariants $\PI(L_{\infty})=(1,4,7)$ and $\RI(L_2)=(2-\sqrt{2},2\sqrt{2}-1,5-\sqrt{2})$, so
$\RM_{\infty}^o(L_2^+,L_{\infty}^+)=\max\{\sqrt{2}-1,5-2\sqrt{2},2+\sqrt{2}\}=2+\sqrt{2}\approx 3.41$ and $\RM_{2}^o(L_2^+,L_{\infty}^+)=\sqrt{6(7-3\sqrt{2})}\approx 4.1$.
\medskip

By Lemma~\ref{lem:RCbounds}(b) the distance between mirror images of the same lattice equals the doubled $D_2$-chiral distance.
For example, $\PM_q^o(L_{\infty}^+,L_{\infty}^-)=2\PC_q[D_2](L_{\infty})=\frac{1}{2}$ and 
$\PM_q^o(L_2^+,L_2^-)=2\PC_q[D_2](L_2)=\frac{2}{2+\sqrt{2}}=2-\sqrt{2}\approx 0.59$ for $q=2,+\infty$.
\medskip

Lemma~\ref{lem:RCbounds}(b) and Table~\ref{tab:RC} also give $\RM_q^o(L_{\infty}^+,L_{\infty}^-)=2\RC_q[D_2](L_{\infty})=2$ and 
$\RM_q^o(L_2^+,L_2^-)=2\RC_q[D_2](L_2)=2(2-\sqrt{2})\approx 1.17$ for $q=2,+\infty$.
\medskip

Lemma~\ref{lem:reversed_signs} says that $\RM^o(L_2^+,L_{\infty}^-)=\RM^o(L_2^-,L_{\infty}^+)$ and $\PM^o(L_2^+,L_{\infty}^-)=\PM^o(L_2^-,L_{\infty}^+)$.
Using the above properties, it remains to find four distances. 
\medskip

\begin{table}[h]
\caption{Metrics $\PM_q^o$ and $\RM_q^o$ for the lattices given by their forms in Table~\ref{tab:RC}, see Fig.~\ref{fig:QS+DC}. }
\medskip

\begin{tabular}{|l|cccc|}
\hline      
$\PM_2^o$ & $L_{\infty}^+$ & $L_{\infty}^-$ & $L_2^+$ & $L_2^-$\\
\hline

$L_{\infty}^+$ & 0 & $\frac{1}{2}$ & $\frac{3}{4}\sqrt{2}-1\approx 0.06$ & $\frac{\sqrt{25-16\sqrt{2}}}{2\sqrt{2}}\approx 0.54$ \\

$L_{\infty}^-$  & $\frac{1}{2}$ & 0 & $\frac{\sqrt{25-16\sqrt{2}}}{2\sqrt{2}}\approx 0.54$ & $\frac{3}{4}\sqrt{2}-1\approx 0.06$  \\

$L_2^+$ & $\frac{3}{4}\sqrt{2}-1\approx 0.06$
& $\frac{\sqrt{25-16\sqrt{2}}}{2\sqrt{2}}\approx 0.54$ & 0 & $2-\sqrt{2}\approx 0.59$ \\

$L_2^-$ & $\frac{\sqrt{25-16\sqrt{2}}}{2\sqrt{2}}\approx 0.54$
 & $\frac{3}{4}\sqrt{2}-1\approx 0.06$ & $2-\sqrt{2}\approx 0.59$ & 0 \\
\hline
\end{tabular}
\medskip

\begin{tabular}{|l|cccc|}
\hline      
$\PM_{\infty}^o$ & $L_{\infty}^+$ & $L_{\infty}^-$ & $L_2^+$ & $L_2^-$\\
\hline

$L_{\infty}^+$ & 0 & $\frac{1}{2}$ & $\frac{3}{4}-\frac{1}{\sqrt{2}}\approx 0.04$ & $\frac{5}{4}-\frac{1}{\sqrt{2}}\approx 0.54$ \\

$L_{\infty}^-$  & 
$\frac{1}{2}$ & 
0 & 
$\frac{5}{4}-\frac{1}{\sqrt{2}}\approx 0.54$ & 
$\frac{3}{4}-\frac{1}{\sqrt{2}}\approx 0.04$  \\

$L_2^+$ & 
$\frac{3}{4}-\frac{1}{\sqrt{2}}\approx 0.04$ & 
$\frac{5}{4}-\frac{1}{\sqrt{2}}\approx 0.54$ & 
0 & 
$2-\sqrt{2}\approx 0.59$ \\

$L_2^-$ & 
$\frac{5}{4}-\frac{1}{\sqrt{2}}\approx 0.54$ & 
$\frac{3}{4}-\frac{1}{\sqrt{2}}\approx 0.04$ & 
$2-\sqrt{2}\approx 0.59$ & 0 \\
\hline
\end{tabular}
\medskip

\begin{tabular}{|l|cccc|}
\hline      
$\RM_{2}^o$ & $L_{\infty}^+$ & $L_{\infty}^-$ & $L_2^+$ & $L_2^-$\\
\hline

$L_{\infty}^+$ & 0 & 
2 & 
$\sqrt{6(7-3\sqrt{2})}\approx 4.1$ & 
$\sqrt{50-22\sqrt{2}}$ \\

$L_{\infty}^-$  & 2 & 0 & $\sqrt{50-22\sqrt{2}}\approx 4.3$ & $\sqrt{6(7-3\sqrt{2})}$  \\

$L_2^+$ & $\sqrt{6(7-3\sqrt{2})}$
& $\sqrt{50-22\sqrt{2}}\approx 4.3$ & 
0 & 
$2(2-\sqrt{2})$ \\

$L_2^-$ & 
$\sqrt{50-22\sqrt{2}}$ & 
$\sqrt{6(7-3\sqrt{2})}\approx 4.1$ & 
$2(2-\sqrt{2})\approx 1.17$ & 0 \\
\hline
\end{tabular}
\medskip

\begin{tabular}{|l|cccc|}
\hline      
$\RM_{\infty}^o$ & $L_{\infty}^+$ & $L_{\infty}^-$ & $L_2^+$ & $L_2^-$\\
\hline

$L_{\infty}^+$ & 0 & 2 & $2+\sqrt{2}\approx 3.41$ & 3 \\

$L_{\infty}^-$  & 2 & 0 & 3 & $2+\sqrt{2}\approx 3.41$  \\

$L_2^+$ & 
$2+\sqrt{2}\approx 3.41$ & 
3 & 
0 & 
$2(2-\sqrt{2})\approx 1.17$ \\

$L_2^-$ & 
3 & 
$2+\sqrt{2}\approx 3.41$ & 
$2(2-\sqrt{2})\approx 1.17$ & 0 \\
\hline
\end{tabular}
\label{tab:RMo}
\end{table}

Proposition~\ref{prop:PMq}(a) finds
$\PM_2^o(L_2^+,L_{\infty}^-)$ as the minimum of the Euclidean distances from $\PI(L_2)=(\frac{1}{2+\sqrt{2}},\frac{1}{2+\sqrt{2}})=(1-\frac{1}{\sqrt{2}},1-\frac{1}{\sqrt{2}})$ to the three points $(-\frac{1}{4},\frac{1}{4})$, $(-\frac{1}{4},\frac{1}{4})$, $(\frac{3}{4},\frac{3}{4})$ obtained from $\PI(L_{\infty})=(\frac{1}{4},\frac{1}{4})$ by reflections in the edges of $\QT$.
The first two distances equal to $\frac{\sqrt{25-16\sqrt{2}}}{2\sqrt{2}}\approx 0.54$ are larger than the third.
\medskip

Given $\PI(L_2)=(x_1,y_1)=(1-\frac{1}{\sqrt{2}},1-\frac{1}{\sqrt{2}})$ and $\PI(L_{\infty})=(x_2,y_2)=(\frac{1}{4},\frac{1}{4})$, Proposition~\ref{prop:PMq}(b) computes
$\PM_{\infty}^o(L_2^+,L_{\infty}^-)$ for as the minimum of
$d_x=\max\{x_2-x_1,y_2+y_1\}=\frac{5}{4}-\frac{1}{\sqrt{2}}$, 
$d_y=\max\{x_2+x_1,|y_2-y_1|\}=\frac{5}{4}-\frac{1}{\sqrt{2}}$, 
$d_{xy}=\max\{x_2-x_1,1-x_2-y_2+|1-y_1-x_2|\}=\frac{1}{4}+\frac{1}{\sqrt{2}}$, so $\PM_{\infty}^o(L_2^+,L_{\infty}^-)=\frac{5}{4}-\frac{1}{\sqrt{2}}\approx 0.54$
\medskip

Proposition~\ref{prop:RMq}(a) computes
$\RM_2^o(L_2^+,L_{\infty}^-)$ as the minimum of the Euclidean distances from $\RI(L_2)=(2-\sqrt{2},2\sqrt{2}-1,5-\sqrt{2})$ to the three points $(-1,4,7)$, $(4,1,7)$, $(1,7,4)$ obtained from $\RI(L_{\infty})=(1,4,7)$ by reflections in the boundaries of $\TC$.
The first distance is the smallest, so $\RM_2^o(L_2^+,L_{\infty}^-)=\sqrt{50-22\sqrt{2}}\approx 4.3$. 
\medskip

Given $\RI(L_2)=(r_{12},r_{01},r_{02})=(2-\sqrt{2},2\sqrt{2}-1,5-\sqrt{2})$ and
$\RI(L_{\infty})=(s_{12},s_{01},s_{02})=(1,4,7)$,
by Proposition~\ref{prop:RMq}(b) 
$\RM_{\infty}^o(L_2^+,L_{\infty}^-)=\min\{d_0,d_1,d_2\}$.
Using $\MS(a,b,c,d)=\max\{|a-b|,|c-d|,\frac{1}{2}|a+b-c-d|\}$, we compute 
\medskip

\noindent
$d_0=\max\{r_{12}+s_{12}, |r_{01}-s_{01}|, |r_{02}-s_{02}|\}$ \\
$=\max\{3-\sqrt{2},5-2\sqrt{2},2+\sqrt{2}\}
=2+\sqrt{2}\approx 3.4,$ 
\medskip

\noindent
$d_1=
\max\{\MS(r_{12},r_{01},s_{12},s_{01}),|r_{02}-s_{02}|\}=$\\
$=\max\{\MS(r_{12},r_{01},s_{12},s_{01}),2+\sqrt{2}\}=$ \\
$=\max\{\MS(2-\sqrt{2},2\sqrt{2}-1,1,4),2+\sqrt{2}\}=$ \\
$=\max\{\max\{3(\sqrt{2}-1),3,2-\frac{1}{\sqrt{2}}\}, 2+\sqrt{2}\}=\max\{3,2+\sqrt{2}\}=2+\sqrt{2},$
\medskip

\noindent
$d_2=\max\{|r_{12}-s_{12}|,\MS(r_{01},r_{02},s_{01},s_{02})\}=$ \\
$=\max\{\sqrt{2}-1, \MS(2\sqrt{2}-1,5-\sqrt{2},4,7)\}=$ \\
$=\max\{\sqrt{2}-1, \max\{6-3\sqrt{2},3,\frac{7-\sqrt{2}}{2})\}=3$, hence $\RM_{\infty}^o(L_2^+,L_{\infty}^-)=3$.
\bs
\end{exa}

\section{Bi-continuity of the map from obtuse superbases to root invariants}
\label{sec:continuity}

Theorems~\ref{thm:isometric_superbases} and~\ref{thm:classification2d} established the bijections $\LIS\to\OSI\to\RIS$, $\La\mapsto B\to\RI(B)=\RI(\La)$, mapping
any lattice $\La\subset\R^2$ to its (unique up to isometry) obtuse superbase $B$ and then to the complete invariant $\RI(\La)$.
Hence the Lattice Isometry Space $\LIS(\R^2)$ having a root metric $\RM$ can be identified with the Root Invariant Space $\RIS=(\TC,d)$ having any metric $d$ on the triangular cone $\TC\subset\R^3$.
\medskip

This section studies continuity of the bijection $B\mapsto\La(B)$, where an obtuse superbase $B$ and its lattice $\La(B)$ are considered up to isometry, rigid motion or two types of similarity.
To state continuity results, we need metrics on lattices and superbases.
Up to each of the four equivalence relations, a lattice $\La$ will be identified with its complete invariant with a relevant metric.
For example, up to isometry, the space $\LIS(\R^2)$ is represented by root invariants $\RI$ with the root metric $\RM$.
Now we define natural metrics on obtuse superbases in any $\R^n$.

\begin{dfn}[space $\OSI$ of obtuse superbases up to isometry]
\label{dfn:OSI}
\textbf{(a)}
Let \\ $B=\{v_i\}_{i=0}^n$ and $B'=\{u_i\}_{i=0}^n$ be any obtuse superbases in $\R^2$.
The \emph{Superbase Isometry Metric} $\SIM_{\infty}(B,B')=\min\limits_{f\in\Or(\R^n)}\max\limits_{i=0,\dots,n}|f(u_i)-v_i|$  minimises vector differences over orthogonal maps $f$ from the group $\Or(\R^n)$.
Let $\OSI(\R^n)$ denote the space of all \emph{obtuse superbases up to isometry} with the metric $\SIM_{\infty}$.
\medskip

\noindent
\textbf{(b)}
Define the space $\OSI^o(\R^n)$ of \emph{obtuse superbases up to rigid motion} (orientation-preserving isometry) with the metric $\SIM_{\infty}^o(B,B')=\min\limits_{f\in\SO(\R^n)}\max\limits_{i=0,\dots,n}|f(u_i)-v_i|$.

\noindent
\textbf{(c)}
Define the spaces $\OSS(\R^n),\OSS^o(\R^n)$ of \emph{obtuse superbases up to similarity} and orientation-preserving similarity with the \emph{Superbase Similarity Metrics} $\SSM_{\infty},\SSM_{\infty}^o$ minimising basis differences over $\Or(\R^n)\times\R_+,\SO(\R^n)\times\R_+$. 
\bs
\end{dfn}

Since any continuous function over a compact domain achieves its minimum value and $\SO(\R^n),\Or(\R^n)$ are compact, the minima in Definition~\ref{dfn:OSI} are achievable.

\begin{lem}[metric axioms for $\SIM_{\infty}$]
\label{lem:OSI_metric_axioms}
The metrics from Definition~\ref{dfn:OSI} on $\OSI(\R^n),\OSI^o(\R^n),\OSS(\R^n),\OSS^o(\R^n)$ satisfy all metric axioms in~(\ref{pro:map}c).
\bt
\end{lem}
\begin{proof}
Let $B_j=\{v_{j0},\dots,v_{jn}\}$, $j=1,2,3$, be any obtuse superbases in $\R^n$.
The first axiom: if $0=\SIM_{\infty}(B_1,B_2)=\min\limits_{f\in\Or(\R^n)}\max\limits_{i=0,\dots,n}|f(v_{1i})-v_{2i}|$, there is an isometry $f\in\Or(\R^n)$ such that $f(B_1)=B_2$.
The superbases $B_1,B_2$ are isometric.
\medskip

Since any isometry $f$ preserves Euclidean distance, we get 
$|f(v_{1i})-v_{2i}|=|f^{-1}(f(v_{1i})-v_{2i})|=|v_{1i}-f^{-1}(v_{2i})|$ and the second axiom: 
$\SIM_{\infty}(B_1,B_2)=\min\limits_{f\in\Or(\R^n)}\max\limits_{i=0,\dots,n}|f(v_{1i})-v_{2i}|=\min\limits_{f^{-1}\in\Or(\R^n)}\max\limits_{i=0,\dots,n}|v_{1i}-f^{-1}(v_{2i})|=\SIM_{\infty}(B_2,B_1)$.
\smallskip

To prove the triangle inequality in the third axiom for $\SIM_{\infty}$, let $f,g\in\Or(\R^n)$ be rotations that minimise the distances $\SIM_{\infty}(B_1,B_2)=\max\limits_{i=0,\dots,n}|f(v_{1i})-v_{2i}|$ and
$\SIM_{\infty}(B_2,B_3)=\max\limits_{i=0,\dots,n}|g(v_{2i})-v_{3i}|$, respectively.
Then we get $\SIM_{\infty}(B_1,B_3)\leq\max\limits_{i=0,\dots,n}|g(f(v_{1i}))-v_{3i}|\leq \max\limits_{i=0,\dots,n}|g(f(v_{1i}))-g(v_{2i})|+\max\limits_{i=0,\dots,n}|g(v_{2i})-v_{3i}|
=$ \\ $\max\limits_{i=0,\dots,n}|f(v_{1i})-v_{2i}|+\SIM_{\infty}(B_2,B_3)
=\SIM_{\infty}(B_1,B_2)+\SIM_{\infty}(B_2,B_3)$.
The proof for other spaces is identical
after replacing $\Or(\R^n)$ with the relevant groups. 
\ws 
\end{proof}


\begin{lem}[bounds for root products]
\label{lem:bounds_products}
Let vectors $u_1,u_2,v_1,v_2\in\R^n$ have a maximum length $l$, have non-positive scalar products $u_1\cdot u_2,v_1\cdot v_2\leq 0$, and be $\de$-close in the Euclidean distance so that 
$|u_i-v_i|\leq\de$ for $i=1,2$.
Then 
$$|u_1\cdot u_2-v_1\cdot v_2|\leq 2l\de \qquad\text{ and }\qquad
|\sqrt{-u_1\cdot u_2}-\sqrt{-v_1\cdot v_2}|\leq\sqrt{2l\de}.
\eqno{\blacktriangle}$$
\end{lem}
\begin{proof}
If $\sqrt{-u_1\cdot u_2}+\sqrt{-v_1\cdot v_2}\leq\sqrt{2l\de}$, the difference of square roots is at most $\sqrt{2l\de}$.
Assuming that $\sqrt{-u_1\cdot u_2}+\sqrt{-v_1\cdot v_2}\geq\sqrt{2l\de}$, it remains to prove that 
$$|u_1\cdot u_2-v_1\cdot v_2|=|\sqrt{-u_1\cdot u_2}-\sqrt{-v_1\cdot v_2}|(\sqrt{-u_1\cdot u_2}+\sqrt{-v_1\cdot v_2})\leq 2l\de.$$
We estimate the scalar product $|u\cdot v|\leq |u|\cdot|v|$ by using Euclidean lengths.
Then we apply the triangle inequality for scalars and replace vector lengths by $l$ as follows: 

\noindent
$|u_1\cdot u_2-v_1\cdot v_2|=|(u_1-v_1)\cdot u_2+v_1\cdot(u_2-v_2)|\leq |(u_1-v_1)\cdot u_2|+|v_1\cdot(u_2-v_2)|\leq$

\noindent
$\leq |u_1-v_1|\cdot |u_2|+|v_1|\cdot|u_2-v_2|\leq \de(|u_2|+|v_1|)\leq 2l\de$ as required.
\end{proof}

\begin{lem}[a lower bound of the size]
\label{lem:size_bound}
If all vectors of an obtuse superbase $B=\{v_0,v_1,v_2\}$ of a lattice $\La\subset\R^2$ have a maximum length $l$, the size $\si(\La)=r_{12}+r_{01}+r_{02}$ of the lattice $\La$ has the lower bound $l\leq \si(\La)$.
\bt
\end{lem}
\begin{proof}
Let $|v_1|=l$. 
Formula~\ref{dfn:vonorms}(a) gives $p_{12}+p_{01}=v_1^2=l^2$ and $r_{12}+r_{01}+r_{02}=\sqrt{(r_{12}+r_{01}+r_{02})^2}\geq\sqrt{r_{12}^2+r_{01}^2+r_{02}^2}=\sqrt{p_{12}+p_{01}+p_{02}}\geq l$.
\ws
\end{proof}

For $q=+\infty$, both $2^{1/q},3^{1/q}$ are interpreted as their limit 1 when $q\to+\infty$.

\begin{thm}[continuity of $\OSI\to\LIS=\RIS$]
\label{thm:OSI->RIS}
\textbf{(a)}
Let $q\in[1,+\infty]$ and lattices $\La,\La'\subset\R^2$ have obtuse superbases $B$ and $B'$ whose vectors have a maximum length $l$.
If $\SIM_{\infty}(B,B')=\de\geq 0$,
then $\RM_q(\La,\La')\leq 3^{1/q}\sqrt{2l\de}$.
Hence the bijection $\OSI(\R^2)\to\LIS(\R^2)$ is continuous in the metrics $\SIM_{\infty}$ and $\RM_q$.
\medskip

\noindent
\textbf{(b)}
In the conditions above, the projected metric satisfies
 $\PM_q(\La,\La')\leq 2^{1/q}3\sqrt{2\de/l}$, so the bijection $\OSS(\R^2)\to\LSS(\R^2)$ is continuous in the metrics $\SIM_{\infty},\PM_q$.
\medskip

\noindent
\textbf{(c)}
In the oriented case, if $\de\to 0$ then $\RM_{q}^o(\La,\La')\to 0$ and $\PM_{q}^o(\La,\La')\to 0$, so both $\OSI^o(\R^2)\to\LIS^o(\R^2)$ and $\OSS^o(\R^2)\to\LSS^o(\R^2)$ are continuous.
\bt
\end{thm}
\begin{proof}
\textbf{(a)}
One can assume that given obtuse superbases $B=(v_0,v_1,v_2)$ and $B'=(u_0,u_1,u_2)$  satisfy $|u_i-v_i|\leq\de$ for $i=0,1,2$ after applying a suitable isometry to $B'$ by Definition~\ref{dfn:OSI}.
Lemma~\ref{lem:bounds_products} implies that the root products $r_{ij}=\sqrt{-v_i\cdot v_j}$ and $\sqrt{-u_i\cdot u_j}$ differ by at most $2l\de$ for any pair $(i,j)$ of indices.
The $M_q$-norm of the vector difference in $\R^3$ is
$\RM_q(\La,\La')\leq 3^{1/q}\sqrt{2l\de}$, $q\in[1,+\infty]$.
\medskip

\noindent
\textbf{(b)}
The coordinates $x=\bar r_{02}-\bar r_{01}$ and $y=3\bar r_{12}$ have an error bound that is at most three times larger than the error for $\bar r_{ij}$.
In Definition~\ref{dfn:PF} each $\bar r_{ij}$ is obtained by dividing the root product $r_{ij}$ by the sizes with the lower bound $l\leq \si=r_{12}+r_{01}+r_{02}$ by Lemma~\ref{lem:size_bound}.
The above error bound $\sqrt{2l\de}$ for $r_{ij}$ gives the error bound $3\sqrt{2\de/l}$ for $x,y$.
Then $\PM_q(\La,\La')\leq 2^{1/q}3\sqrt{2\de/l}$, $q\in[1,+\infty]$.
\medskip

\noindent
\textbf{(c)}
In the oriented case, if $\sign(\La)>0$, then $\RI(\La)$ is strictly inside the triangular cone $\TC$.
The continuity of $\RM$ implies that under any continuous motion of $\de$-close superbases $B\to B'$, if $\de$ is sufficiently small, then all intermediate lattices have their unoriented root invariants inside $\TC$, so their signs are positive and $\sign(\La')>0$.
Hence $\sign(\La),\sign(\La')$ coincide or one of them is 0.
In all cases by Definition~\ref{dfn:RMo} the metric $\RM^o(\La,\La')$ coincides with $\RM(\La,\La')$ whose convergence to $0$ as $\de\to 0$ was proved above.
The proof of $\PM^o(\La,\La')\to 0$ is similar.
\ws
\end{proof}

Theorem~\ref{thm:OSI->RIS} is proved for the metrics $\RM_q,\PM_q$ only to give explicit upper bounds.
A similar argument proves continuity for any metrics $\RM,\PM$ in Definition~\ref{dfn:RM} based on a metric $d$ satisfying $d(u,v)\to 0$ when $u\to v$ coordinate-wise. 
All Minkowski norms in $\R^n$ are topologically equivalent due to the bounds  $||v||_q\leq ||v||_{r}\leq n^{\frac{1}{q}-\frac{1}{r}}||v||_q$ for any $1\leq q\leq r$ \cite{norms}, hence continuity for one value of $q$ is enough.
Theorem~\ref{thm:OSI->RIS} implies continuity of $\OSI^o\to\RIS^o$, because closeness of superbases up to rigid motion is a stronger condition than up to isometry.
\medskip

Example~\ref{exa:RF_deformation} illustrates Theorem~\ref{thm:OSI->RIS} and shows that the root invariant changes continuously for a deformation when a reduced basis changes discontinuously.

\begin{exa}[continuity of root invariants]
\label{exa:RF_deformation}
Let the basis vectors $v_1(t)=(1,0),v_2(t)=(t,1)$ vary for $t\in[0,1]$.
The resulting lattice $\La(t)$ starts its deformation from the unit square lattice $\La_4=\La(0)$, deforms to the centred rectangular lattice $\La=\La(\frac{1}{2})$ in the second picture of Fig.~\ref{fig:reduced_bases} and then returns to $\La_4=\La(1)$.
The initial basis $\{v_1(t),v_2(t)\}$ by Definition~\ref{dfn:reduced_cell} remains reduced for $t\in[0,\frac{1}{2}]$ and then discontinuously jumps to $v_1(t)=(1,0)$ and $v'_2=(t-1,1)$ for $t\in(\frac{1}{2},1]$.
\medskip

The root invariant changes continuously as follows.
The superbase $v_1=(1,0)$, $v'_2=(t-1,1)$, $v_0=(-t,-1)$ of $\La(t)$ remains obtuse for $t\in[0,1]$ and is unique up to isometry by Theorem~\ref{thm:isometric_superbases}.
Then $r_{12}=\sqrt{1-t}$, $r_{01}=\sqrt{t}$, $r_{02}=\sqrt{1-t+t^2}$.
Since $1-t+t^2\geq\max\{t,1-t\}$ for $t\in[0,1]$, the root invariant can be written as 
$$\RI(\La(t))=\left\{\begin{array}{ll}
(\sqrt{t},\sqrt{1-t},\sqrt{1-t+t^2}) \text{ for } t\in[0,\frac{1}{2}],\\
(\sqrt{1-t},\sqrt{t},\sqrt{1-t+t^2}) \text{ for } t\in[\frac{1}{2},1].
\end{array} \right.$$

\begin{figure}[h]
\includegraphics[width=\textwidth]{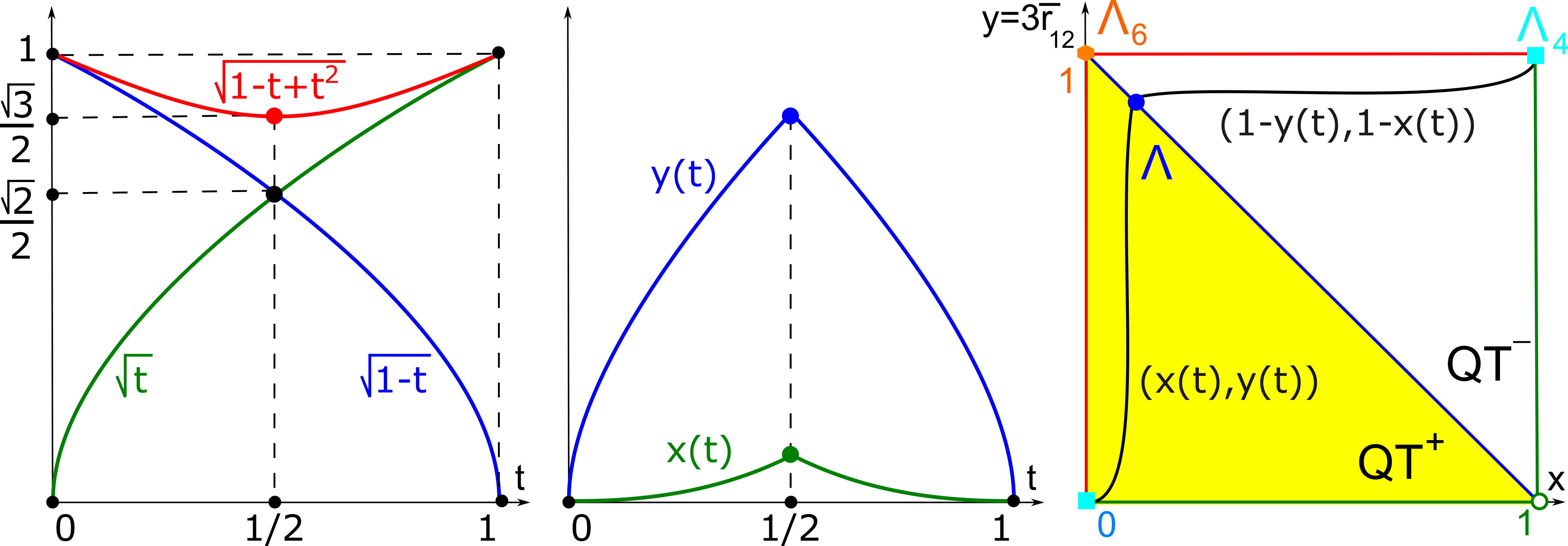}
\caption{
\textbf{Left}: graphs of root products in $\RI(\La(t))$, see Example~\ref{exa:RF_deformation}.
\textbf{Middle}: graphs of the components in $\PI(\La(t))$.
\textbf{Right}: the continuous path of $\PI(\La(t))$ in the quotient square $QS$.}
\label{fig:RF_deformation}
\end{figure}

By Definition~\ref{dfn:PF} the size is $\si(\La(t))=r_{12}+r_{01}+r_{02}=\sqrt{t}+\sqrt{1-t}+\sqrt{1-t+t^2}$.
The projected invariant is $\PI(\La(t))=(x(t),y(t))$, see Fig.~\ref{fig:RF_deformation}. where
$$x(t)=\frac{\sqrt{1-t+t^2}-\max\{\sqrt{t},\sqrt{1-t}\}}{\sqrt{t}+\sqrt{1-t}+\sqrt{1-t+t^2}},\quad
y(t)=\frac{3\min\{\sqrt{t},\sqrt{1-t}\}}{\sqrt{t}+\sqrt{1-t}+\sqrt{1-t+t^2}}.$$
If $t=\frac{1}{2}$, then $\sqrt{t}=\sqrt{1-t}=\frac{\sqrt{2}}{2}$, $\sqrt{1-t+t^2}=\frac{\sqrt{3}}{2}$, $\si(\La(\frac{1}{2}))=\sqrt{2}+\frac{\sqrt{3}}{2}$.
So
$$\RI\left(\La\Big(\frac{1}{2}\Big)\right)=\left(\frac{\sqrt{2}}{2},\frac{\sqrt{2}}{2},\frac{\sqrt{3}}{2}\right),\; 
\PI\left(\La\Big(\frac{1}{2}\Big)\right)=\left(\frac{\sqrt{3}-\sqrt{2}}{\sqrt{3}+2\sqrt{2}},\frac{3\sqrt{2}}{\sqrt{3}+2\sqrt{2}}\right).$$
The last point is approximately $(0.07,0.93)$ in the diagonal $x+y=1$ of the quotient square $\QS$.
Under the symmetry $t\lra 1-t$, all the functions above remain invariant and $\La(t)$ changes its sign.
The path $\RI(\La(t))$ and its projection $\PI(\La(t))\in\QS$ are continuous everywhere, while the reduced basis is discontinuous at $t=\frac{1}{2}$.
\bs
\end{exa}

Theorem~\ref{thm:RIS->OSI} below proves the inverse continuity 
of $\RIS\to\OSI$ and a weaker claim in the oriented case saying that we can choose an obtuse superbase $B'$ of a perturbed lattice arbitrarily close to a given superbase $B$ of an original lattice.

\begin{thm}[continuity of $\LIS\to\OSI$]
\label{thm:RIS->OSI}
\textbf{(a)}
For $q\in[1,+\infty]$, let lattices $\La,\La'$ in $\R^2$ satisfy $\RM_q(\La,\La')\leq\de$.
For any obtuse superbase $B$ of $\La$, there is an obtuse superbase $B'$ of $\La'$ such that $\SIM_{\infty}(B,B')\leq\SIM_{\infty}^o(B,B')\to 0$ as $\de\to 0$.
\medskip

\noindent
\textbf{(b)}
The bijection $\LIS(\R^2)\to\OSI(\R^2)$ is continuous in the metrics $\RM_q,\SIM_{\infty}$.
\medskip

\noindent
\textbf{(c)}
The above conclusions hold for lattices and superbases up to similarity.
\bt
\end{thm}
\begin{proof}
\textbf{(a)}
The obtuse superbase $B=(v_0,v_1,v_2)$ of $\La$ is already given.
Let $B'=(u_0,u_1,u_2)$ be any obtuse superbase of $\La'$ found from $\RI(\La')$ by Lemma~\ref{lem:superbase_reconstruction}.
\medskip

Up to rigid motion in $\R^2$, one can assume that $\La,\La'$ share the origin and the first vectors $v_0,u_0$ lie in the positive $x$-axis.
Let $r_{ij},s_{ij}$ be the root products of $B,B'$, respectively.
Formulae~(\ref{dfn:vonorms}a) imply that $v_i^2=r_{ij}^2+r_{ik}^2$ and $u_i^2=s_{ij}^2+s_{ik}^2$ for distinct indices $i,j,k\in\{0,1,2\}$, for example if $i=0$ then $j=1$, $k=2$. 
\medskip

For any given continuous transformation of root invariants from $\RI(\La)$ to $\RI(\La')$, all root products have a finite upper bound $M$, which we use to estimate 
$$|v_i^2-u_i^2|=|(r_{ij}^2+r_{ik}^2)-(s_{ij}^2+s_{ik}^2)|\leq
|r_{ij}^2-s_{ij}^2|+|r_{ik}^2-s_{ik}^2|=$$
$$(r_{ij}+s_{ij})|r_{ij}-s_{ij}|+(r_{ik}+s_{ik})|r_{ik}-s_{ik}|\leq 
(r_{ij}+s_{ij})\de+(r_{ik}+s_{ik})\de\leq 4M\de.$$

Since at least two continuously changing conorms are strictly positive to guarantee positive lengths of basis vectors by formula~(\ref{dfn:vonorms}a), all basis vectors reconstructed by Lemma~\ref{lem:superbase_reconstruction} have a minimum length $a>0$.
Then $||v_i|-|u_i||\leq\dfrac{4M\de}{|v_i|+|u_i|}\leq\dfrac{2M}{a}\de$.
Since the vectors $v_0,u_0$ lie in the positive horizontal axis, the lengths can be replaced by vectors: $|v_0-u_0|\leq\dfrac{2M}{a}\de$, so $|v_0-u_0|\to 0$ as $\de\to 0$.
\medskip

If the superbases $B,B'$ have opposite signs, apply to $B'$ the reflection  with respect the fixed $x$-axis. 
To conclude that $\SIM_{\infty}^o(B,B')\to 0$, we show below that the basis vectors $v_i,u_i$ from both superbases have close angles $\al_i,\be_i$ measured anticlockwise from the positive $x$-axis for $i=1,2$.
To estimate the small difference $\al_i-\be_i$, we first express the angles via the root products by Lemma~\ref{lem:superbase_reconstruction}:
$$\al_i=\arccos\dfrac{v_0\cdot v_i}{|v_0|\cdot|v_i|}
=\arccos\dfrac{-r_{0i}^2}{\sqrt{r_{01}^2+r_{02}^2}\sqrt{r_{ij}^2+r_{ik}^2}},$$
$$
\be_i=\arccos\dfrac{u_0\cdot u_i}{|u_0|\cdot|u_i|}
=\arccos\dfrac{-s_{0i}^2}{\sqrt{s_{01}^2+s_{02}^2}\sqrt{s_{ij}^2+s_{ik}^2}}
,$$
where $j\neq k$ differ from $i=1,2$.
If $\de\to 0$, then $s_{ij}\to r_{ij}$ and $\al_i-\be_i\to 0$ for all indices because all functions above are continuous for $|u_j|,|v_j|\geq a$, $j=0,1,2$.
\medskip

We estimate the squared length of the difference by using the scalar product:
$$|v_i-u_i|^2=v_i^2+u_i^2-2u_iv_i
=(|v_i|^2-2|u_i|\cdot |v_i|+|u_i|^2)+2|u_i|\cdot |v_i|-2|u_i|\cdot |v_i|\cos(\al_i-\be_i)$$
$$
=(|v_i|-|u_i|)^2+2|u_i|\cdot |v_i|(1-\cos(\al_i-\be_i))
=(|v_i|-|u_i|)^2+|u_i|\cdot |v_i|4\sin^2\dfrac{\al_i-\be_i}{2}
\leq$$
$$\leq (|v_i|-|u_i|)^2+|u_i|\cdot |v_i|4\left(\dfrac{\al_i-\be_i}{2}\right)^2
=(|v_i|-|u_i|)^2+|u_i|\cdot |v_i|(\al_i-\be_i)^2,$$
where we used that $|\sin x|\leq|x|$ for $x\in\R$. 
The upper bound $M$ of root products guarantees a fixed upper bound for the lengths $|u_i|,|v_i|$.
If $\de\to 0$, then $|v_i|-|u_i|\to 0$ and $\al_i-\be_i\to 0$ as proved above, so
 $v_i-u_i\to 0$ and $\SIM_{\infty}^o(B,B')\to 0$.
\medskip

\noindent
\textbf{(b)}
Since the metric $\SIM_{\infty}$ from Definition~\ref{dfn:OSI} is minimised over the larger group $\Or(\R^2)$ in comparison with $\SO(\R^2)$, we have the inequality $\SIM_{\infty}(B,B')\leq\SIM_{\infty}^o(B,B')$, hence $\SIM_{\infty}(B,B')\to 0$ as $\de\to 0$ by part (a).
Up to isometry, the obtuse superbases $B,B'$ of $\La,\La'$ are unique by Theorem~\ref{thm:isometric_superbases}.
Since we can start with any obtuse superbase $B$ and can also apply a reflection to $B'$, the above convergence $\SIM_{\infty}(B,B')\to 0$ for any $B,B'$ proves continuity of $\LIS\to\OSI$.
\medskip

\noindent
\textbf{(c)}
To extend part (a) to the similarity equivalence, we use the size $\si=r_{12}+r_{01}+r_{02}$ of the given superbase $B$ to reconstruct an obtuse superbase $B'$ of $\La'$ from $\PI(\La')$ with the same size $\si$ by Proposition~\ref{prop:inverse_design}.
By formula~(\ref{prop:inverse_design}) the given condition $\PM_q(\La,\La')\to 0$ implies that $\RM_q(\La,\La')\to 0$, hence part (a) implies the required conclusions for the smaller metrics $\SSM_{\infty}\leq\SIM_{\infty}$ and $\SSM_{\infty}^o\leq\SIM_{\infty}^o$.
The same argument extends part~(b) to the similarity case.
\ws
\end{proof}

Lemma~\ref{lem:SIM_bound} proves a non-trivial lower bound needed for Corollary~\ref{cor:rect_discontinuity} later.

\begin{lem}[lower bound for $\SIM_{\infty}^o$]
\label{lem:SIM_bound}
For any ordered obtuse superbases $B,B'$ of lattices in $\R^2$ with coforms $\CF(B)=(p_{12},p_{01},p_{02})$ and $\CF(B')=(p'_{12},p'_{01},p'_{02})$, let $\CM_{\infty}(B,B')=\min\limits_{\zeta\in A_3} \max\limits_{i\neq j}\{|p_{ij}-p'_{\zeta(i)\zeta(j)}|\}$ be minimised is over three cyclic permutations $\zeta\in A_3$ of indices $0,1,2$.
Let $l$ be a maximum length of all vectors from $B,B'$.
Then we have the lower bound $\SIM_{\infty}^o(B,B')\geq \CM_{\infty}(B,B')/l$.
\bt
\end{lem}
\begin{proof}
For the superbase $B=\{v_0,v_1,v_2\}$, by Definition~\ref{dfn:OSI}find an optimal rotation around the origin so that the resulting image $\{v'_0,v'_1,v'_2\}$ of $B'$ satisfies $|u_i-v_i|\leq\SIM_{\infty}^o(B,B')$, $i=0,1,2$.
Lemma~\ref{lem:bounds_products} implies that $|p_{ij}-p'_{ij}|\leq 2l\cdot\SIM_{\infty}^o(B,B')$ for all distinct $i,j\in\{0,1,2\}$.
The above rotation might have cyclically shifted the coforms of $B'$, but $\CM_{\infty}(B,B')$ is minimised over cyclic permutations.
Then $\CM_{\infty}(B,B')\leq \max\limits_{i\neq j}\{|p_{ij}-p'_{ij}|\}\leq2l\cdot\SIM_{\infty}^o(B,B')$ gives the lower bound.
\ws
\end{proof}

One can prove that min-max distance in Lemma~\ref{lem:SIM_bound} satisfies metric axioms in $\R^3$.
Corollary~\ref{cor:rect_discontinuity} shows that Theorem~\ref{thm:RIS->OSI}(b) is the strongest possible continuity in the oriented case.
In $\R^3$, a similar discontinuity around high-symmetry lattices will be much harder to resolve for continuous invariants even up to isometry \cite{kurlin2022complete}.

\begin{cor}[partial discontinuity of $\RIS^o\to\OSI^o$]
\label{cor:rect_discontinuity}
The bijection $\LIS^o\to\OSI^o$ is discontinuous in the metrics $\RM_{\infty},\SIM_{\infty}^o$ at any rectangular lattice.
\bt
\end{cor}
\begin{proof}
For any $0\leq 3\de<a<b$, start from any rectangular lattice with a unit cell $a\times b$ and consider the lattices $\La^{\pm}(\de)\subset\R^2$ with the obtuse superbases 
\smallskip

\noindent
\begin{tabular}{llll}
$B^+(\de)$ : & 
$v_1=(a,0)$ & 
$v_2^{+}(\de)=(-\de,b)$ & 
$v_0^+(\de)=(\de-a,-b)$ \\
$B^-(\de)$ : & $v_1=(a,0)$ & 
$v_2^-(\de)=(\de-a,b)$ & 
$v_0^{-}(\de)=(-\de,-b)$
\end{tabular}
\medskip

Notice that the vectors in both superbases are ordered anticlockwise around $0$.
The initial lattice $\La^{\pm}(0)$ has two superbases $v_1=(a,0)$, $v_2^\pm(0)=(0,\pm b)$, $v_0=(-a,\mp b)$ related by reflection, not by rigid motion, see Fig.~\ref{fig:hexagonal_graphene}~(right).
\medskip

Keeping the anticlockwise order above, write the ordered coforms below.
\smallskip

\hspace*{-8mm}
\begin{tabular}{llll}
$\CF(B^+(\de))$ & 
$-v_1\cdot v_2^+=\de a$ & 
$-v_0\cdot v_1^+=a^2-\de a$ & 
$-v_0\cdot v_2^+=b^2-\de a+\de^2$\\
$\CF(B^-(\de))$ & 
$-v_1\cdot v_2^-=a^2-\de a$ & 
$-v_0\cdot v_1^-=\de a$ & 
$-v_0\cdot v_2^-=b^2-\de a+\de^2$
\end{tabular}
\medskip

The above coforms differ by the transposition of the first two conorms.
The maximum difference of all corresponding conorms in $\CF(B^\pm(\de))$ is $a^2-2\de a$.
If we cyclically shift $\CF(B^-(\de))$ to the left, the maximum difference becomes $b^2-a^2+\de^2$.
If we cyclically shift $\CF(B^-(\de))$ to the right, the maximum difference becomes $b^2-2\de a+\de^2$.
By Lemma~\ref{lem:SIM_bound}, the cyclic metric between the above coforms is 
$$\CM_{\infty}(B^+(\de),B^-(\de))=\min\{\, a^2-2\de a,\,  b^2-a^2+\de^2,\, b^2-2\de a+\de^2 \, \}\geq\min\{\frac{a^2}{3},b^2-a^2\}$$ 
due to $\de<\frac{a}{3}$.
Since the maximum length of vectors from $B^\pm(\de)$ is $l\leq\sqrt{a^2+b^2}$, we get $\SIM_{\infty}^o(B^+(\de),B^-(\de))\geq \CM_{\infty}(B^+(\de),B^-(\de))/l\geq\min\{\frac{a^2}{3},b^2-a^2\}/\sqrt{a^2+b^2}$.
This lower bound shows that, for any $0<\de<\frac{a}{3}$, the (unique up to rigid motion) obtuse superbases $B^\pm(\de)$ of $\La^\pm(\de)$ are not close in the metric $\SIM_{\infty}^o$.
\medskip

The lattices $\La^\pm(\de)$ have 
$\RI(\La^\pm(\de))$ consisting of $\sqrt{\de a},\sqrt{a^2-\de a},\sqrt{b^2-a^2+\de^2}$, which might need to be ordered. 
Since the lattices $\La^\pm(\de)$ are related by reflection,
Lemma~\ref{lem:RCbounds}(b) computes $\RM(\La^+(\de),\La^-(\de))$ as the double distance $2\RC[D_2](\La(\de))$ depending only on the root invariant above without signs.
For the Minkowski parameter $q=+\infty$, Proposition~\ref{prop:RC}(a) computes the required distance as follows:
$$\RC_{\infty}[D_2](\La(\de))=\min\{\de a,\frac{a^2-2\de a}{2}, \frac{b^2-2a^2+\de a+\de^2}{2}\}\leq\de a\to 0 \text{ as }\de\to 0.$$  
Hence the lattices $\La^\pm(\de)$ have close root invariants with $\RM_{\infty}(\La^+(\de),\La^-(\de))\to 0$ as $\de\to 0$, but their obtuse superbases have a constant lower bound for the metric 
$\SIM_{\infty}^o$ independent of $\de$.
The discontinuity conclusion holds for all $q\in[1,+\infty)$, because all Minkowski distances $M_q$ are topologically equivalent \cite{norms}.
\ws
\end{proof}

Corollary~\ref{cor:rect_discontinuity} should be positively interpreted in the sense that we need to study lattices up to rigid motion by their complete oriented root invariants in the continuous space $\LIS^o(\R^2)$ rather than in terms of reduced bases (or, equivalently, obtuse superbases due to Proposition~\ref{prop:reduced_bases}b), which are inevitably discontinuous.
\medskip

Proposition~\ref{prop:cont_chiral_distances} shows that all $G$-chiral distances $\RC[G]:\LIS(\R^2)\to\R$ and $\PC[G]:\LSS(\R^2)\to\R$ are continuous in any metrics $\RM,\PM$ from Definition~\ref{dfn:RM}.

\begin{prop}[continuous chiral distances]
\label{prop:cont_chiral_distances}
For a crystallographic point group $G$ and lattices $\La_1,\La_2$ in $\R^2$, we have 
$|\RC[G](\La_1)-\RC[G](\La_2)|\leq \RM(\La_1,\La_2)$
and
$|\PC[G](\La_1)-\PC[G](\La_2)|\leq \PM(\La_1,\La_2)$
for any metrics $\RM$ and $\PM$. 
\bt
\end{prop}
\begin{proof}
In Definition~\ref{dfn:RC} let $\La_1,\La_2\in\LIS[G]$ be lattices that minimise $\RC[G](\La_1)=\RM(\La_1,\La_1')$ and $\RC[G](\La_2)=\RM(\La_2,\La_2')$.
The triangle inequality implies that
$$\RC[G](\La_1)\leq \RM(\La_1,\La'_2)\leq \RM(\La_1,\La_2)+\RM(\La_2,\La'_2)=\RM(\La_1,\La_2)+\RC[G](\La_2)$$
and $\RC[G](\La_1)-\RC[G](\La_2)\leq \RM(\La_1,\La_2)$.
Swapping indices $1\lra 2$, we similarly get 
$\RC[G](\La_2)-\RC[G](\La_1)\leq \RM(\La_1,\La_2)$.
Hence we get the required upper bound 
$|\RC[G](\La_1)-\RC[G](\La_2)|\leq \RM(\La_2,\La_1)$.
The proof for $\PC[G]$ is similar.
\ws
\end{proof}

\section{New mathematical structures on lattices, conclusions, and discussion}
\label{sec:conclusions}

This section first connects the recent invariants of more general periodic point sets with the complete invariants of lattices.
Then we discuss linear operations, scalar products, CAT(0) property of $\LIS(\R^2)$ and finally describe the future work.
\medskip

Below we prove that other continuous isometry invariants  $\AMD$ (average minimum distances) and $\PDD$ (pointwise distance distribution) are complete for lattices, though they make sense for  general periodic and finite point sets \cite{widdowson2021pointwise,widdowson2022average}.
Other isometry invariants such as persistent homology turned out to be weaker than expected, see infinite families of sets that have identical persistence in \cite{smith2022families}.

\begin{dfn}[RSD invariant]
\label{dfn:RSD}
For any lattice $\La\subset\R^n$, both AMD and PDD invariants reduce to the sequence of distances $(d_1,d_1,d_2,d_2,d_3,d_3,\dots)$ from the origin $0\in\La$ to its $k$-th nearest neighbour in $\La$ for $k\geq 1$.
Since any $\La$ is symmetric with respect to $0$, define the \emph{Reduced Sequence of Distances} $\RSD(\La)=(d_1,d_2,d_3,\dots)$ containing one distance from each pair of equal distances above.
\bs
\end{dfn}

\begin{figure}[h]
\includegraphics[width=1.0\textwidth]{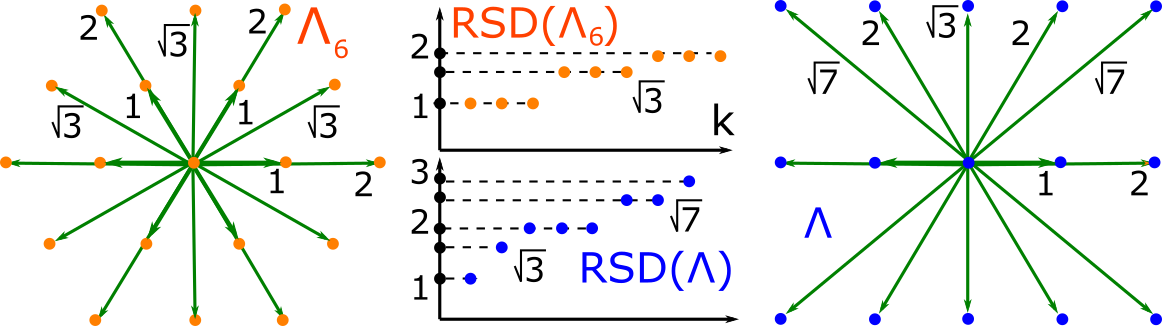}
\caption{\textbf{Left and Right}: first neighbours of the origin $0$ in the hexagonal lattice $\La_6$ and rectangular lattice $\La$ with unit cell $1\times\sqrt{3}$.
\textbf{Middle}: $\RSD(\La_6)$ and $\RSD(\La)$, see Definition~\ref{dfn:RSD}.
}
\label{fig:hex_rect_RSD}
\end{figure}

In 1938 Delone reduced $\RSD(\La)$ even further and considered only distinct increasing distances \cite[p.~163]{delone1938geometry}.  
He proved that the resulting weaker invariant (of only the first four distinct distances) is complete for all lattices $\La\subset\R^2$ except the two lattices $\La_6,\La$ in Fig.~\ref{fig:hex_rect_RSD}, which are distinguished by 
the stronger $\RSD$:
$\begin{array}{l}
\RSD(\La_6)=(1,1,1,\sqrt{3},\sqrt{3},\sqrt{3},2,2,2,\sqrt{7},\sqrt{7},\sqrt{7},\sqrt{7},\sqrt{7},\sqrt{7},3,3,3,\dots),\\
\RSD(\La)=(1,\sqrt{3},2,2,2,\sqrt{7},\sqrt{7},3,\dots).
\end{array}$

Any quadratic form $Q(x,y)$ is uniquely determined (up to a linear change of variables) by the sets of its values with the only exception of $Q_6=x^2+xy+y^2$ and $Q=x^2+3y^2$ corresponding to the lattices in Fig.~\ref{fig:hex_rect_RSD}, see references in \cite{watson1980determination}.
\medskip

For any lattice $\La\subset\R^n$, the `halved' sequence $\RSD$ contains the same information as $\AMD$ and $\PDD$.
We conjecture that $\PDD$ is complete for all finite and periodic points sets in $\R^2$.
Proposition~\ref{prop:RSD} implies completeness for lattices in $\R^2$.
\medskip

For any lattice $\La\subset\R^n$, $\RSD(\La)$ can be represented as the \emph{theta series} $\Theta_{\La}(q)=\sum\limits_{v\in\La}q^{|v|^2}=1+2\sum\limits_{k\geq 1}q^{\RSD_k^2(\La)}$, $q\in\C$.
Lecture~2 in \cite[p.~45]{conway2013sphere} mentions that $\Theta_{\La(q)}$ determines the shape of any lattice $\La\subset\R^2$, which is proved below via $\RSD(\La)$.

\begin{prop}[$\RSD$ completeness]
\label{prop:RSD}
Any lattices $\La,\La'\subset\R^2$ are isometric if and only if $\RSD(\La)=\RSD(\La')$.
\bt
\end{prop}
\begin{proof}
The lengths $|v_1|\leq|v_2|\leq|v_0|$ of shortest Voronoi vectors from Fig.~\ref{fig:Voronoi2D} are not necessarily the first three distances in $\RSD$.
For example, if $\La$ has the basis $v_1=(1,0)$, $v_2=(0,3)$, then $\RSD(\La)=(1,2,3,\dots)$, where $2=2|v_1|\neq|v_2|=3$.
We will extract $|v_1|\leq|v_2|\leq|v_0|$ from $\RSD(\La)$, which proves completeness.
\medskip

For any integer $k>1$, the shortest Voronoi vector $v_1\in\La$ of length $d_1=|v_1|$ generates a single integer multiple $kv_1\in\La$ of length $kd_1$, which can be removed from $\RSD(\La)$.
The second shortest Voronoi vector $v_2\in\La$ may accidentally have the same length $|v_2|=d_1$.
If not, the next distance $d_2$ in the resulting sequence equals $|v_2|$ as in $\RSD(\La)=\{1,\sqrt{3},\dots\}$ for the lattice $\La$ in Fig.~\ref{fig:hex_rect_RSD}~(right).
If yes, we recognise the repeated value of $d_1$ as $d_2=|v_2|$ and do not confuse $d_2$ with any multiple $kd_1$ for $k>1$.
For example, if $\La$ has the basis $v_1=(1,0)$, $v_2=(0,2)$, then $\RSD(\La)=\{1,2,2,\dots\}$, where one distance $2$ is $2d_1$, another distance $2$ is $d_2$, so we remove only one distance $2=2d_1$.
Again for any $k>1$, we remove one multiple $kd_2$ and find the next distance $d_3$ equal to the length of the third shortest Voronoi vector $v_0=-v_1-v_2$ as in Fig.~\ref{fig:Voronoi2D}.
Since $(v_1^2,v_2^2,v_0^2)=(d_1^2,d_2^2,d_3^2)$ is  a complete invariant by Theorem~\ref{thm:classification2d} and Lemma~\ref{lem:VF}, then so is $\RSD(\La)$.
\ws
\end{proof}

The invariant $\RSD(\La)$ can be made complete for lattices up to rigid motion by adding $\sign(\La)$ and up to similarity after dividing all distances by the first $d_1$.
\medskip

Now Remark~\ref{rem:structures} summarises a wide range of rich mathematical structures that can be considered on the lattice spaces in addition to continuous metrics.

\begin{rem}[linear structure, scalar product on lattices]
\label{rem:structures}
Since the triangular cone $\TC$ in Fig.~\ref{fig:TC} is convex, we can consider any convex linear combination of root invariants $t\RI(\La_1)+(1-t)\RI(\La_2)\in\TC$, $t\in[0,1]$.
The resulting root invariant determines (an isometry class of) the new lattice that can be denoted by $t\La_1+(1-t)\La_2$.
The average of the square and hexagonal lattices with $\RI(\La_4)=(0,1,1)$, $\RI(\La_6)=(1,1,1)$ has $\RI=(\frac{1}{2},1,1)$.
The new lattice $\frac{1}{2}(\La_4+\La_6)$ is centred rectangular and has the basis $v_1=(\sqrt{\frac{3}{2}},0)$ and $v_2=(-\frac{1}{9}\sqrt{\frac{3}{2}},\frac{4}{9}\sqrt{\frac{15}{2}})$.
We can define similar sums in $\LSS(\R^2)$ due to convexity of the triangle $\QT$. 
\medskip

The usual scalar product of vectors in $\R^3$ defines the positive product of root invariants: $\RI(\La_4)\cdot\RI(\La_6)=(0,1,1)\cdot(1,1,1)=2$.
The lattice spaces $\LIS^o(\R^2)$ and $\LSS^o(\R^2)$ up to rigid motion and orientation-preserving similarity are geodesic metric spaces, even CAT(0) spaces and flat manifolds (locally Euclidean).
\bs
\end{rem}

In conclusion, Problem~\ref{pro:map} was resolved by the new invariants $\RI,\RI^o,\PI,\PI^o$ classifying all 2D lattices up to four equivalences, see a summary in Table~\ref{tab:classifications}.
\smallskip

\noindent
(\ref{pro:map}a) Invariants :
Definitions~\ref{dfn:RI} and~\ref{dfn:sign}, 
Lemma~\ref{lem:RI_invariance} and Lemma~\ref{lem:lattice_invariants}(a).
\smallskip

\noindent
(\ref{pro:map}b) Completeness of invariants : 
Theorem~\ref{thm:classification2d} and Corollary~\ref{cor:similar_lattices2d}.
\smallskip

\noindent
(\ref{pro:map}c) Continuous metrics : 
Definitions~\ref{dfn:RM} and~\ref{dfn:RMo}, Theorems~\ref{thm:OSI->RIS} and~\ref{thm:RIS->OSI}.
\smallskip

\noindent
(\ref{pro:map}d) Computability of metrics : 
Propositions~\ref{prop:RMq},\ref{prop:PMq} with Examples~\ref{exa:RM},~\ref{exa:RC}. 
\smallskip

\noindent
(\ref{pro:map}e) Inverse design : 
Lemma~\ref{lem:superbase_reconstruction} and
Proposition~\ref{prop:inverse_design} with Example~\ref{exa:inverse_design}.
\medskip

\begin{table}[h]
\caption{A summary of classifications of all lattices $\La\subset\R^2$ up to four equivalence relations.}

\hspace*{-4mm}
\begin{tabular}{c|cccc}      
equivalence  & complete invariant & configuration space & continuous metric & visual results \\
\hline

any   & root invariant  & $\LIS(\R^2)\lra\TC$ & root metric &  
Theorem~\ref{thm:classification2d} \\

isometry & $\RI(\La)$ & triangular cone & $\RM$ & Fig.~\ref{fig:TC} (left) \\
\hline

rigid & oriented  & $\LIS^o(\R^2)\lra\DC$ & oriented & Theorem~\ref{thm:classification2d} \\

motion & invariant $\RI^o(\La)$ & doubled cone & metric $\RM^o$ & Fig.~\ref{fig:QS+DC} (right) \\
\hline

any & projected  & $\LSS(\R^2)\lra\QT$ & projected & Corollary~\ref{cor:similar_lattices2d} 
\\

similarity & invariant $\PI(\La)$ & quotient triangle & metric $\PM$ & Fig.~\ref{fig:QT+QS} (left) \\
\hline

orientation- & oriented & $\LSS^o(\R^2)\lra\QS$ & oriented & Corollary~\ref{cor:similar_lattices2d} \\

preserving & projected & quotient & projected & Fig.~\ref{fig:QT+QS} (right) \\

similarity & invariant $\PI^o(\La)$ & square & metric $\PM^o$ & Fig.~\ref{fig:DT}
\end{tabular}
\label{tab:classifications}
\end{table}
\medskip

The key contributions are the easily computable metrics in Definitions~\ref{dfn:RM},\ref{dfn:RMo}, which led to continuous real-valued deviations of lattices from their higher symmetry neighbours.
The chiral distances in Definition~\ref{dfn:RC} continuously extend the classical binary chirality and have explicit formulae in Propositions~\ref{prop:RC},\ref{prop:PC}.
\medskip

The discontinuity of reduction in \cite[Theorem~15]{widdowson2022average} was proved with a simple metric on bases without isometry. 
When we consider obtuse superbases up to isometry, continuity holds in Theorem~\ref{thm:RIS->OSI} without orientation.
If orientation should be preserved, Corollary~\ref{cor:rect_discontinuity} proves discontinuity at any rectangular lattice in $\R^2$.
\medskip

The structures in Remark~\ref{rem:structures} help treat lattices as vectors in a meaningful way (independent of a basis), for example, as inputs or outputs in machine learning algorithms.
Future work \cite{kurlin2022complete,bright2021welcome} extends key results to 3D lattices.
The author thanks any reviewers for their valuable time and helpful suggestions in advance. 

\begin{acknowledgements}
This research was supported by the £3.5M EPSRC grant `Application-driven Topological Data Analysis' (2018-2023), the £10M Leverhulme Research Centre for Functional Materials Design (2016-2026) and the Royal Academy of Engineering Fellowship `Data Science for Next Generation Engineering of Solid Crystalline Materials' (2021-2023).
\end{acknowledgements}


\bibliographystyle{spmpsci}      
\bibliography{lattices2Dmaths}   

\renewcommand{\thesection}{\Alph{section}}
\setcounter{section}{0}
\section{Appendix~A: detailed proofs of past results by Conway-Sloane}
\label{sec:proofs}

To make the paper self-contained, the appendix includes detailed proofs of past results whose outlines in Delone \cite{delone1975bravais}, Conway, Sloane \cite{conway1992low} contained a few typos.
\medskip

Lemma~\ref{lem:Voronoi_classification} (probably due to Voronoi) was mentioned in \cite[section~2.3]{delone1975bravais}.

\begin{lem}[lattices $\lra$ Voronoi domains]
\label{lem:Voronoi_classification}
Lattices $\La,\La'\subset\R^n$
are related by an isometry $f$
if and only if Voronoi domains $V(\La),V(\La')$ are related by $f$.
\bt  
\end{lem}
\begin{proof}
Since any isometry $f$ preserves distances and the Voronoi domain is defined in terms of distances, if $f$ maps a lattice $\La$ to $\La'$, then $f$ restricts to an isometry of Voronoi domains: $V(\La)\to V(\La')$.
Conversely, the whole space $\R^n$ is covered by the lattice translates $V(\La)+\La=\{V(\La)+v \mid v\in\La\}$, which have disjoint interiors.
Hence any isometry $f:V(\La)\to V(\La')$ gives rise to an isometry of $\R^n$.
\ws
\end{proof}

\begin{proof}[of Lemma~\ref{lem:Voronoi_vectors}]
We prove the second part for strict Voronoi vectors with all strict inequalities.
The first part follows by making all inequalities non-strict. 
\medskip

Assume that $v\in\La$ is a Voronoi vector of a lattice $\La\subset\R^2$ but there is a shorter vector $w\in v+2\La$.
Then the point $\frac{1}{2}v$ is closer to $\frac{1}{2}(v+w)\in\La$ than to $v$, because $\frac{1}{2}|w|<\frac{1}{2}|v|$, which contradicts the assumption that $v$ is a Voronoi vector.
\medskip

Conversely, let $v$ be a shortest vector in its $2\La$-class.
If $v$ is not a Voronoi vector, there is another vector $w\in\La$ whose bisector hyperspace separates $\frac{1}{2}v$ from $0$.
Then $\frac{1}{2}v\cdot w>\frac{1}{2}|w|^2$, 
$|v^2|>|v-2w|^2$, so $v-2w$ is shorter than $v$.
\ws
\end{proof}

\begin{lem}
\label{lem:squared_norm}
For any basis $v_1,\dots,v_n$ in $\R^n$, let $p_{ij}=-v_i\cdot v_j$ be the conorms of a superbase $v_0,v_1,\dots,v_n$ with $v_0=-\sum\limits_{i=1}^n v_i$.
The squared norm $v^2$ of any vector $v=\sum\limits_{i=1}^n c_i v_i$ equals 
$N(v)=\sum\limits_{i=1}^n c_i^2 p_{0i}+\sum\limits_{1\leq i<j\leq n}(c_i-c_j)^2 p_{ij}$.
If we decompose $v=\sum\limits_{i=0}^n c_i v_i$ over the full superbase, then
$N(v)=\sum\limits_{0\leq i<j\leq n}(c_i-c_j)^2 p_{ij}$.
\bt
\end{lem}
\begin{proof}
In the right hand side of the required formula, we substitute the conorms in terms of scalar products of basis vectors as follows: 
$p_{ij}=-v_i\cdot v_j$ for $i,j\in\{1,\dots,n\}$.
Then $p_{0i}=-v_0\cdot v_i=v_i\cdot\sum\limits_{j=1}^n v_j=v_i^2+\sum\limits_{j\neq i}^n v_i\cdot v_j$ and
$$\sum\limits_{i=1}^n c_i^2 p_{0i}+\sum\limits_{1\leq i<j\leq n}(c_i-c_j)^2 p_{ij}
=\sum\limits_{i=1}^n c_i^2 (v_i^2+ \sum\limits_{j\neq i}^n v_i\cdot v_j)-\sum\limits_{1\leq i<j\leq n}(c_i^2-2c_ic_j+c_j^2) v_i\cdot v_j$$
$$=\sum\limits_{i=1}^n c_i^2 v_i^2+\sum\limits_{i=1}^n c_i^2 \sum\limits_{j\neq i}^n v_i\cdot v_j-\sum\limits_{1\leq i<j\leq n}(c_i^2+c_j^2)(v_i\cdot v_j)+\sum\limits_{1\leq i<j\leq n}2c_ic_j (v_i\cdot v_j)=$$
$=\mathlarger{\sum}\limits_{i=1}^n c_i^2 v_i^2+2\mathlarger{\sum}\limits_{1\leq i<j\leq n}c_ic_j (v_i\cdot v_j)=\left(\mathlarger{\sum}\limits_{i=1}^n c_i v_i\right)^2=N(v)$ as required.
The last formula for $n+1$ vectors follows by replacing $c_i$ with $c_i-c_0$ for $i=1,\dots,n$.
\ws
\end{proof}

\begin{proof}[of Lemma~\ref{lem:partial_sums}]
By Lemma~\ref{lem:Voronoi_vectors} Voronoi vectors have smallest squared norms $N(v)$ in their $2\La$-classes.
The $2\La$-class of any vector $v=\sum\limits_{i=0}^n c_i v_i\in\La$ with $c_i\in\Z$ remains invariant of any coefficient $c_i$ keeps its parity modulo 2.
Within the $2\La$-class, the squared norm $N(v)=\sum\limits_{0\leq i<j\leq n}(c_i-c_j)^2 p_{ij}$ computed in Lemma~\ref{lem:squared_norm} is minimal if all even $c_i$ are replaced by 0 and all odd $c_i$ are replaced by 1. 
The resulting shortest vectors with coefficients $0,1$ are all $2^n-1$ symmetric pairs of partial sums $\pm v_S$ for a proper subset $S\subset\{0,1,\dots,n\}$.
If all conorms $p_{ij}>0$, to guarantee a minimum value of $N(v)$, every difference $|c_i-c_j|$ should be 0 or 1, hence there are no other Voronoi vectors apart from the partial sums above. 
\ws
\end{proof}

\begin{proof}[of Theorem~\ref{thm:reduction} for $n=2$]
For any lattice $\La\subset\R^2$, permuting vectors of a superbase $B=(v_0,v_1,v_2)$ allows us to order the conorms: $p_{12}\leq p_{01}\leq p_{02}$.
Our aim is to reduce $B$ so that all $p_{ij}\geq 0$.
Assuming that $p_{12}=-v_1\cdot v_2=-\ep<0$, we change the superbase: $u_1=-v_1$, $u_2=v_2$, $u_0=v_1-v_2$ so that $u_0+u_1+u_2=0$.
\medskip

Two vonorms remain the same: $u_1^2=v_1^2$, $u_2^2=v_2^2$.
The third vonorm decreases by $4\ep>0$ as follows: $u_0^2=(v_1-v_2)^2=(v_1+v_2)^2-4v_1v_2=v_0^2-4\ep$.
One conorm changes its sign: $q_{12}=-u_1\cdot u_2=-p_{12}=\ep>0$.
The two other conorms decrease:
$$q_{01}=-u_0\cdot u_1=-(v_1-v_2)\cdot(-v_1)
=-(-v_1-v_2)v_1-2v_1\cdot v_2=p_{01}-2\ep,$$  
$$q_{02}=-u_0\cdot u_2=-(v_1-v_2)\cdot v_2
=-(-v_1-v_2)v_2-2v_1\cdot v_2=p_{02}-2\ep.$$  
If one of the new conorms becomes negative, we apply the above reduction again.
\medskip

To prove that all conorms eventually become non-negative, note that every reduction can make superbase vectors only shorter, but not shorter than a minimum distance between points of $\La$.
The angle between $v_i,v_j$ can have only finitely many values when lengths of $v_i,v_j$ are bounded.
Then the scalar product $\ep=v_i\cdot v_j>0$ cannot converge to 0.
Since every reduction makes one superbase vectors shorter by a positive constant, the reductions will finish in finitely many steps.
\ws
\end{proof}

\end{document}